\documentclass[11pt]{amsart}
\usepackage[english]{babel}
\usepackage[T1]{fontenc}
\usepackage[utf8]{inputenc} 

\usepackage{mathtools}
\usepackage{tikz}

\usepackage{xpatch}

\usepackage{parskip}
\usepackage{blindtext}
\usepackage{adjustbox}
\usepackage{stmaryrd}
\usepackage{amsmath}
\usepackage{amssymb}
\usepackage{amsthm}

\newtheoremstyle{note}
{2pt}
{2pt}
{\itshape}
{0pt}
{\bfseries}
{.}
{.5em}
{}
\theoremstyle{note}

\usepackage[all,cmtip]{xy}
\usepackage{tikz-cd}
\usepackage{setspace}
\usepackage{enumitem}
\setlist[itemize]{noitemsep, topsep=0pt}
\usepackage{mathrsfs}
\usepackage{subfiles}
\usepackage{csquotes}
\usepackage{booktabs}
\usepackage{xcolor}
\usepackage[margin=.8in]{geometry}
\usepackage{fancyhdr}
\usepackage{braket}
\usepackage{amscd}
\usepackage{extarrows}
\usepackage{dsfont}
\usepackage{scalerel}
\usepackage{scalerel,stackengine}
\usepackage{xpatch}

\usetikzlibrary{arrows,calc,matrix,decorations}

\newsavebox{\pullback}
\sbox\pullback{%
\begin{tikzpicture}%
\draw (0,0) -- (1ex,0ex);%
\draw (1ex,0ex) -- (1ex,1ex);%
\end{tikzpicture}}
\usepackage[colorlinks,citecolor=blue,menucolor=magenta,urlcolor=blue,bookmarks=false,hypertexnames=true]{hyperref}

\usepackage{diagbox}
\usepackage{float}
\definecolor{codegreen}{rgb}{0,0.6,0}
\definecolor{codegray}{rgb}{0.5,0.5,0.5}
\definecolor{codepurple}{rgb}{0.58,0,0.82}
\definecolor{backcolour}{rgb}{0.95,0.95,0.92}
\makeatletter

\makeatletter
\xpatchcmd{\@thm}{\thm@headpunct{.}}{\thm@headpunct{ :}}{}{}
\makeatother

\newtheorem{pro}{Proposition}[subsection]

\newtheorem{rk}[pro]{Remark}

\usepackage{xparse}
\let\oldrk\rk
\RenewDocumentCommand{\rk}{o}{%
  \IfNoValueTF{#1}
    {\oldrk}
    {\oldrk[#1]}%
  \normalfont
}

\newtheorem{ex}[pro]{Example}
\newtheorem{lemma}[pro]{Lemma}
\newtheorem{foreshadow}[pro]{Proof sketch}
\newtheorem{defi}[pro]{Definition}

\newtheorem{corro}[pro]{Corollary}
\newtheorem{thm}[pro]{Theorem}

\newtheorem{prodef}[pro]{Proposition--Definition}

\newtheorem*{pro_no}{Proposition}

\newcommand{\lpar}{\left(}
\newcommand{\rpar}{\right)}

\renewcommand{\part}[1]{\addcontentsline{toc}{part}{#1}}

\newcommand\sublsection[1]{\Alph{subsection}}

\newcommand{\A}{\ensuremath{\mathbb{A}}} 
\newcommand{\cA}{\ensuremath{\mathcal{A}}} 
\newcommand{\cB}{\ensuremath{\mathcal{B}}} 
\newcommand{\bB}{\ensuremath{\mathbb{B}}} 
\newcommand{\C}{\ensuremath{\mathbb{C}}} 
\newcommand{\cC}{\ensuremath{\mathcal{C}}} 
\newcommand{\D}{\ensuremath{\mathcal{D}}} 
\newcommand{\cD}{\ensuremath{\mathcal{D}}} 
\newcommand{\e}{\ensuremath{\Acute{e}}} 
\newcommand{\et}{\ensuremath{\Acute{e}t}}
\newcommand{\F}{\ensuremath{\mathbb{F}}} 
\newcommand{\cF}{\ensuremath{\mathcal{F}}} 
\newcommand{\sF}{\ensuremath{\mathscr{F}}}
\newcommand{\bG}{\ensuremath{\mathds{G}}} 
\newcommand{\G}{\ensuremath{\mathcal{G}}} 
\newcommand{\cG}{\ensuremath{\mathcal{G}}} 
\renewcommand{\L}{\ensuremath{\mathcal{L}}}  
\newcommand{\cL}{\ensuremath{\mathcal{L}}} 
\newcommand{\m}{\ensuremath{\mathfrak{m}}} 
\newcommand{\N}{\ensuremath{\mathbb{N}}} 
\renewcommand{\O}{\ensuremath{\mathcal{O}}} 
\newcommand{\cP}{\ensuremath{\mathcal{P}}} 
\renewcommand{\P}{\ensuremath{\mathbb{P}}} 
\newcommand{\Q}{\ensuremath{\mathbb{Q}}} 
\newcommand{\R}{\ensuremath{\mathbb{R}}} 
\newcommand{\T}{\ensuremath{\mathbb{T}}} 
\newcommand{\bT}{\ensuremath{\mathbb{T}}} 

\newcommand{\Z}{\ensuremath{\mathbb{Z}}} 
\renewcommand{\bar}[1]{\overline{#1}} 

\newcommand{\rlim}{\ensuremath{\varinjlim}}
\newcommand{\llim}{\ensuremath{\varprojlim}}
\renewcommand{\=}{\ensuremath{\overset{a}{=}}}

\newcommand{\acong}{\ensuremath{\overset{a}{\cong}}}

\newcommand{\Ab}{\operatorname{Ab}}

\newcommand{\Hom}{\operatorname{Hom}} 
 
\newcommand{\Gal}{\operatorname{Gal}} 
\newcommand{\Fil}{\operatorname{Fil}} 
\newcommand{\Rep}{\operatorname{Rep}} 
\newcommand{\Isoc}{\operatorname{Isoc}} 
\newcommand{\Grad}{\operatorname{Grad}}

\newcommand{\Spa}{\operatorname{Spa}} 
\newcommand{\Spec}{\operatorname{Spec}}

\newcommand{\Dbs}{\cD^a_{\square}}

\newcommand{\hotimes}{\widehat{\otimes}}

\newcommand{\la}{\left\langle}
\newcommand{\ra}{\right\rangle}

\renewcommand{\1}{\mathds{1}}

\newcommand{\doublefaktor}[3]{%
    {\textstyle #1}
    \mkern-4mu\scalebox{1.5}{$\diagdown$}\mkern-5mu^{\textstyle #2}%
    \mkern-4mu\scalebox{1.5}{$\diagup$}\mkern-5mu{\textstyle #3} }

\makeatletter
\newcommand*{\relrelbarsep}{.386ex}
\newcommand*{\relrelbar}{%
  \mathrel{%
    \mathpalette\@relrelbar\relrelbarsep
  }%
}

\newcommand*{\@relrelbar}[2]{%
  \raise#2\hbox to 0pt{$\m@th#1\relbar$\hss}%
  \lower#2\hbox{$\m@th#1\relbar$}%
}
\providecommand*{\rightrightarrowsfill@}{%
  \arrowfill@\relrelbar\relrelbar\rightrightarrows
}
\providecommand*{\leftleftarrowsfill@}{%
  \arrowfill@\leftleftarrows\relrelbar\relrelbar
}
\providecommand*{\xrightrightarrows}[2][]{%
  \ext@arrow 0359\rightrightarrowsfill@{#1}{#2}%
}
\providecommand*{\xleftleftarrows}[2][]{%
  \ext@arrow 3095\leftleftarrowsfill@{#1}{#2}%
}
\makeatother

\stackMath
\newcommand\reallywidehat[1]{%
\savestack{\tmpbox}{\stretchto{%
  \scaleto{%
    \scalerel*[\widthof{\ensuremath{#1}}]{\kern-.6pt\bigwedge\kern-.6pt}%
    {\rule[-\textheight/2]{1ex}{\textheight}}
  }{\textheight}%
}{0.5ex}}%
\stackon[1pt]{#1}{\tmpbox}%
}

\makeatletter
\renewcommand\@listI{%
\leftmargin=25pt
\rightmargin=0pt
\labelsep=10pt
\labelwidth=20pt
\itemindent=0pt
\listparindent=2pt
\topsep=5pt plus 2pt minus 4pt
\partopsep=0pt plus 1pt minus 1pt
\parsep=0pt plus 1pt
\itemsep=\parsep}

\makeatletter
\renewcommand\paragraph{\@startsection{paragraph}{5}{\z@}%
  {3.25ex \@plus1ex \@minus.2ex}%
  {-1em}%
  {\normalfont\normalsize\bfseries}}
\renewcommand\subparagraph{\@startsection{subparagraph}{6}{\parindent}%
  {3.25ex \@plus1ex \@minus .2ex}%
  {-1em}%
  {\normalfont\normalsize\bfseries}}
\def\toclevel@subsubsubsection{4}
\def\toclevel@paragraph{5}
\def\toclevel@paragraph{6}
\def\l@subsubsubsection{\@dottedtocline{4}{7em}{4em}}
\def\l@paragraph{\@dottedtocline{5}{10em}{5em}}
\def\l@subparagraph{\@dottedtocline{6}{14em}{6em}}
\makeatother

\usetikzlibrary{decorations.markings}

\makeatletter
\tikzcdset{
  open/.code     = {\tikzcdset{hook, circled};},
  closed/.code   = {\tikzcdset{hook, slashed};},
  open'/.code    = {\tikzcdset{hook', circled};},
  closed'/.code  = {\tikzcdset{hook', slashed};},
  circled/.code  = {\tikzcdset{markwith = {\draw (0,0) circle (.375ex);}};},
  slashed/.code  = {\tikzcdset{markwith = {\draw[-] (-.4ex,-.4ex) -- (.4ex,.4ex);}};},
  markwith/.code ={
    \pgfutil@ifundefined%
    {tikz@library@decorations.markings@loaded}%
    {\pgfutil@packageerror{tikz-cd}{You need to say %
      \string\usetikzlibrary{decorations.markings} to use arrows with markings}{}}{}%
    \pgfkeysalso{/tikz/postaction = {
      /tikz/decorate,
      /tikz/decoration={markings, mark = at position 0.5 with {#1}}}
    }
  },
}
\makeatother

\usepackage{subfiles}

\title{Mod $p$ Poincaré duality for $p$-adic period domains}
\author{Guillaume PIGNON-YWANNE}

\newcommand{\GL}{\mathbb{GL}}

\renewcommand{\Rep}{\operatorname{Rep}}
\renewcommand{\L}{\mathbb{L}}

\usepackage[backend=biber,style=alphabetic,maxbibnames=99, maxnames = 99]{biblatex}
\AtEndDocument{
    \printbibliography}

\begin{document}

\begin{abstract}
    In this article, we introduce a new class of smooth partially proper rigid analytic varieties over a $p$-adic field that satisfy Poincaré duality for étale cohomology with mod $p$-coefficients : the varieties satisfying "primitive comparison with compact support". We show that almost proper varieties, as well as p-adic (weakly admissible) period domains in the sense of Rappoport-Zink \cite{RappZink}, belong to this class. In particular, we recover Poincaré duality for almost proper varieties as first established by Li-Reinecke-Zavyalov \cite{Zavyalov&co}, and we compute the étale cohomology with $\F_p$-coefficients of p-adic period domains, generalizing a computation of Colmez-Dospinescu-Niziol \cite{CDN_IntegralDrinfeld} for Drinfeld's symmetric spaces. The arguments used in this paper rely crucially on Mann's six functors formalism for solid $\O^{+,a}/\pi$ coefficients as developed in \cite{mann2022-6F}.
\end{abstract}
\maketitle
\setcounter{tocdepth}{2}
{\setlength{\parskip}{0pt}
\tableofcontents}

\section{Introduction}


\subsection{Main results}
Let $K$ be a nonarchimedean field of mixed characteristic $(0, p)$, $C$ be a complete algebraically closed extension of $K$, and $\Lambda = \Z/n\Z$ be a fixed ring of coefficients for some $n \in \N_{>1}$. Let $X$ be an analytic variety over $C$, in an sense that we voluntarily keep imprecise for now.

\smallskip
There are fully developed theories of étale cohomology for sheaves of $\Lambda$-modules on analytic varieties over $K$ in many incarnations of analytic geometry, constructed by Berkovich \cite{Berkovich_ECohD}, Huber \cite{Huber_Ecoh}, or Scholze \cite{Scholze_ECohD}. When $n$ is coprime to $p$, étale cohomology admits similar properties to the ones satisfied by étale cohomology of schemes in characteristic zero, as developed by Grothendieck and Verdier in \cite{SGA4}. Indeed, the derived\footnote{One needs to pass to left-completion in the setup of spatial diamonds.} category of étale sheaves of $\Lambda$-modules over some $X$ : $\cD_{\e t}(X, \Lambda)$ admits a six functor formalism ; in particular, it satisfies Poincaré duality. However, the Poincaré duality for rigid analytic spaces fails badly when $p \mid n$ - famously, it does not hold for the open unit ball over $C$ (see, for example, \cite[Rk 6.4.11]{Zavyalov&co}).

\smallskip
Nonetheless, there have been recent advances in establishing Poincaré duality with $\F_p$-coefficients for some class of \textit{nice} analytic spaces. It has been proven for smooth and proper rigid-analytic varieties over $C$ simultaneously by L.Mann in \cite{mann2022-6F} and B.Zavyalov in \cite{zavyalov2024modp} (following unpublished ideas of O.Gabber). This has then been extended to Zariski-open subsets of such smooth proper rigid analytic varieties in \cite{Zavyalov&co}, extending upon the approach of Zavyalov.

\smallskip
In this paper, we extend Poincaré duality to smooth partially proper rigid analytic varieties satisfying an analogue of the primitive comparison theorem for compactly supported étale cohomology, expanding on Mann's approach in \cite{mann2022-6F}.

Our main result is as follows : 

\begin{thm} (cf. \ref{PCCS_PD}) 
    \label{PrimtiveComparisonhasPD}
    Let $K$ be a complete extension of $\Q_p$, $\overline{K}$ be its algebraic closure, and $C$ be the completion of $\overline{K}$. Let $X$ be a partially proper proper rigid analytic variety over $K$ that is smooth of pure dimension $d$, and let $X_C$ be its base change to $C$. 

    Let $\cL$ be an overconvergent\footnote{i.e. the stalks of $\cL$ are constant on specializations (cf. \ref{Overconv} for a precise definition). This includes local systems and Zariski-constructible sheaves in the sense of \cite[Defi 3.1]{BH_Sixfunctors_Zariskiconstr}.} étale sheaf of $\F_p$-vector spaces on $X$ and $\cL_C$ denote its pullback to $X_C$. 
    
    Assume furthermore that, for any $i \in \N$, the natural morphism : 
    \begin{equation*}
        H^i_{\e t, c}(X_C, \cL_C) \otimes \O_C/p \to H^i_{\e t, c}(X_C, \cL_C \otimes \O_{X_C}^+/p)
    \end{equation*} is an almost isomorphism of $\O_C/p$-modules. When this is the case, we say that $X_C$ satisfies primitive \textbf{comparison with compact support} with respect to $\cL_C$.
    
    Then, for all $0 \le k \le 2d$, there exists a natural $\operatorname{Gal}(\overline{K}/K)$-equivariant isomorphism : $$H^k_{\e t}(X_C, \cL_C^\vee(d)) \cong \operatorname{Hom}_{\F_p}(H^{2d-k}_{\e t, c}(X_C, \cL_C), \F_p),$$
    where $\cL_C^\vee$ denotes the dual sheaf $\cL_C^\vee = \underline{\operatorname{Hom}}_{\F_p}(\cL_C, \F_p)$, and $(d)$ denotes a Tate twist.
\end{thm}

A proof sketch of the above will be laid out in \ref{Proof_subseq}.

\smallskip
Note that the Poincaré duality in the above theorem does \textit{not} always come from a perfect pairing : $$H^k_{\e t}(X_C, \cL_C^\vee) \otimes_{\F_p} H^{2d-k}_{\e t, c}(X_C, \cL_C) \to \F_p(-d)$$ as cohomology groups may not be finite dimensional. In particular, one cannot a priori recover compactly supported cohomology groups $H^\bullet_{\e t, c}$ from $H^\bullet_{\e t}$.

If we furthermore assume $X$ to be proper and $\cL$ to be a local system, primitive comparison with compact support follows from Scholze's primitive comparison theorem \cite[Thm 5.1]{Scholze_padHTforRAV}, so that the above theorem recovers Poincaré duality as in \cite[Thm 1.1.1]{mann2022-6F}. 

When $X$ is not necessarily proper, primitive comparison with compact support does not hold in general, even for $\cL = \F_p$ - a direct computation shows that the open unit disk over $C$ does not satisfy it.

We also establish primitive comparison with compact support for so-called "almost proper" varieties, i.e. rigid analytic varieties that can be written as a Zariski open subset of a smooth proper rigid analytic variety. In particular, we recover a version of \cite[Corro 1.1.2]{Zavyalov&co} : 

\begin{corro} (cf. \ref{Almost_proper})
    Let $K$ be a complete extension of $\Q_p$, $\overline{K}$ be its algebraic closure, and $C$ be the completion of $\overline{K}$. Let $X$ be a smooth proper rigid analytic variety over $K$ of pure dimension $d$, and $U \subset X$ be a Zariski open subset. Let $\cL$ be an étale local system of $\F_p$-vector spaces on $U$.
    
    Let $U_C$ denote the base change of $U$ to $C$, and $\cL_C$ be the pullback of $\cL$ to $U_C$. Then $U_C$ satisfies primitive comparison with compact support with respect to $\cL_C$, so that, for any $ 0 \le k \le 2d$, there is a natural $\operatorname{Gal}(\overline{K}/K)$-isomorphism : $$H^k_{\e t}(U_C,  \cL_C^\vee(d)) \cong \operatorname{Hom}_{\F_p}(H^{2d-k}_{\e t, c}(U_C, \cL_C), \F_p).$$
\end{corro}

This corollary follows straightforwardly from Thm. \ref{PrimtiveComparisonhasPD}, together with a version of the primitive comparison theorem for Zariski-constructible coefficients, as considered in \cite{BH_Sixfunctors_Zariskiconstr}.

\smallskip
In the second part of the article, we prove that $p$-adic period domains, in the sense of  Rappoport-Zink \cite[Defi 1.36]{RappZink} satisfy primitive comparison with compact support. Let us first state our result for Drinfeld symmetric spaces, which are important examples of such period domains :

\begin{pro}
    Let $C$ be the completion of an algebraic closure of $\Q_p$, and $d \ge 1$ be a fixed integer. Let $\Omega^d_{\Q_p}$ be the Drinfeld symmetric space of dimension $d$ over $\Q_p$, and let $\Omega^d_{C}$ be its base change to $C$.
    
    Then, $\Omega^d_{C}$ satisfies primitive comparison with compact support with $\F_p$-coefficients.
    
    In particular, $\Omega^d_{C}$ satisfies Poincaré duality with $\F_p$-coefficients, so that, for any $0 \le k \le 2d$, there is a natural $\operatorname{Gal}(\overline{\Q}_p/\Q_p) \times \operatorname{GL_{d+1}(\Q_p)}$-equivariant isomorphism : $$ H^{k}_{\e t}(\Omega^d_{C}, \F_p(d)) \cong \operatorname{Hom}_{\F_p}(H_{\e t,c}^{2d-k}(\Omega^d_{C}, \F_p), \F_p).$$
\end{pro}

It now follows from the computation of compactly supported cohomology of Drinfeld symmetric spaces in \cite[Corro 6.8]{cdhn} that the étale cohomology of Drinfeld symmetric spaces with $\F_p$-coefficients is given by the Tate-twisted dual of a generalised Steinberg representation, as defined in \cite[Section 4.1]{CDN_IntegralDrinfeld} for all $0 \le k \le d$ :\footnote{Note that the Steinberg representation is sometimes denoted $\operatorname{Sp}$ (for SPecial representation) in the litterature.} $$H_{\e t}^{k}(\Omega^d_{C}, \F_p) \cong \operatorname{St}_{k}^*(\F_p)(k).$$ This recovers the computation of \cite[Thm 1.1]{CDN_IntegralDrinfeld}. Note that the approach of loc.cit. does not stem from a Poincaré duality result, but rather an explicit computation using formal models and $A_{inf}$-cohomology, in the sense of \cite{BMS1}.

The above result stems from a more general statement about primitive comparison with compact support for p-adic period domains, as follows (cf. section \ref{PerDom} for a more detailed introduction) :

Let $\Breve{\Q}_p$ denote the completion of the maximal unramified extension of $\Q_p$. Consider the datum $(G, [b], \{ \mu \})$ consisting of : 
    \begin{itemize}
        \item A connected and quasi-split reductive algebraic group $G$ over $\Q_p\, $.
        \item An element $[b]$ of the Kottwitz set $B(G)$ of $\sigma$-conjugacy classes in $G(\Breve{\Q}_p)$, associated to some $b \in G(\Breve{\Q}_p)$ that is assumed to be basic and of $s$-decent\footnote{The decency hypothesis is not a constraint, as any conjugacy class contains a decent element.} for some integer $s \ge 1$.
        \item A conjugacy class $\{\mu\}$ of geometric cocharacters of $G$.
    \end{itemize}

    Let $E$ be the field of definition of $\{ \mu \}$, also known as its reflex field. The datum of $(G, \{ \mu\})$ induces a flag variety $\sF(G, \{ \mu \})$, defined over $E$, classifying par-equivalence classes of cocharacters in $\{\mu \}$.
    
    Let $\Breve{E} = E \cdot \Breve{\Q}_p$, and $J_b$ be the inner form of $G$ defined as the $\sigma$-centralizer of $b$. The period domain associated to $(G, [b], \{ \mu\})$, denoted $\sF^{wa}(G, [b], \{\mu\})$, is the admissible partially proper open subset of $\sF(G, \{\mu\}) \otimes_E \Breve{E}$ parameterizing weakly admissible filtrations of type $\{ \mu\}$ on $N_b$. It admits a canonical model over $E_s := E \cdot \Q_{p^s}$, still denoted $\sF^{wa}(G, [b], \{\mu\})$, and is endowed with a natural action of $J_b(\Q_p)$.

We may now state our result :

\begin{thm} \label{PCCS_Period_Dom_Intro}(cf. \ref{PCCS_Period_Domains})
    Let $(G, b, [\mu])$ be as above, and $\sF^{wa}(G, [b], \{\mu\})$ be the associated period domain over $E_s$. Let $C = \overline{E_s}$, and let $\sF^{wa}_C$ be the base change of $\sF^{wa}(G, [b], \{\mu\})$ to $C$.
    
    Then, $\sF^{wa}_{C}$ satisfies primitive comparison with compact support with $\F_p$-coefficients.

    It hence satisfies Poincaré duality, in the sense that, letting $d$ be its dimension, for all $0 \le k \le 2d$, there is a natural $\operatorname{Gal}(\overline{E_s}/E_s) \times J_b(\Q_p)$-equivariant isomorphism : $$ H^{k}_{\e t}(\sF^{wa}_C, \F_p(d)) \cong \operatorname{Hom}_{\F_p}(H_{\e t,c}^{2d-k}(\sF^{wa}_C, \F_p), \F_p).$$
\end{thm}

The proof of the above theorem relies on adapting the computation of étale cohomology with compact support of $\sF^{wa}_C$, done in \cite{cdhn} with $\F_p$-coefficients, to $\O^+/p$ coefficients. This relies on a stratification of the complement of the inclusion of $\sF_C^{wa}$ inside the flag variety $\sF_C$ due to \cite{Orlik_local}. Using geometric arguments, we reduce to establishing a comparison result for Schubert cells, for which Scholze's primitive comparison theorem holds.

\smallskip
In particular, using the computation of the compactly supported cohomology of period domains in \cite[Thm 1.2]{cdhn}, the above theorem allows us to compute the étale cohomology of $\sF^{wa}_C$ with $\F_p$-coefficients. We now state the result after introducing the necessary notation.

\smallskip
    Let $\G_{E_s} := \operatorname{Gal}(\overline{E_s}/E_s)$ denote the absolute Galois group of $E_s$. 
    
    Let $T$ be a maximum torus of $G$ such that $\mu$ factors through $T$, and $W = N(T)/T$ be the absolute Weyl group of $G$ with respect to $T$, which acts naturally on $X_*(T)$. 
    
    Let $W^\mu$ be the set of Kostant representatives with respect to $W/\operatorname{Stab}(\mu)$. The action of $\G_{E_s}$ on $W$ preserves $W^\mu$. To each $\G_{E_s}$-orbit $[w] \in W^\mu/\G_{E_s}$, we associate : 
\begin{itemize}
    \item Its length, an integer denoted by $l_{[w]} \,$, defined as the length of any element of the orbit.
    \item A $\F_p$-module $\rho_{[w]}(\F_p)$ of $\F_p$-valued functions on $[w]$, equipped with the induced $\G_{E_s}$-action twisted by $(-l_{[\omega]})$, and $\rho_{[w]}(\F_p)^*$ denotes its $\F_p$-linear dual.
\end{itemize}

Consider a maximal $\Q_p$-split torus $S$ of the derived group $J_{b, der}$ of $J_b$ contained in $T$, and a maximal $\Q_p$-parabolic subgroup $P_0$ of $J$ containing $S$. Let $\Delta \subset X^*(S)$ be the associated set of relative simple roots. For any $I \subset \Delta$, we let $P_I$ be the corresponding standard parabolic subgroup of $J_b$, and $X_I =J_b(\Q_p)/P_I(\Q_p)$, viewed as a profinite topological space. We let :
 $$v_{P_I}^{J_b}(\F_p) = \operatorname{LC}(X_I, \F_p)/\sum_{I \subsetneq I'} \operatorname{LC}(X_{I'}, 
\F_p)$$ be the corresponding generalized Steinberg representation of $J_b(\Q_p)$, where $\operatorname{LC}(-, \F_p)$ denotes the module of locally constant $\F_p$-valued functions. We let $v_{P_I}^{J_b}(\F_p)^* = \operatorname{Hom}(v_{P_I}^{J_b}(\F_p), \F_p)$ be the dual Steinberg representation, equipped with trivial Galois action. 

Fix an invariant inner product $(-,-)$ on $G$. For each $\G_{E_s}$-orbit $[w]$ of some $w \in W^\mu$, we define $$ I_{[w]} = \{\alpha \in \Delta, (w \mu - v,\omega_\alpha) \le 0 \}, \, \, \, P_{[w]} := P_{I_{[w]}},$$
where $(\omega_\alpha)_{\alpha \in \Delta} \in X_*(S) \otimes \Q$ forms the dual basis of $\Delta$.

\begin{corro}
    Let $\sF^{wa}_C$ be a period domain as in the above setup over $C = \overline{E_s}$, and let $d$ be its dimension.\footnote{Note that $d$ can be explicited as $d = (2\rho, \nu_b)$ where $\rho$ denotes the half sum of the positive roots of $\Delta$, $\nu_b$ is the slope morphism associated to $b$ as defined in \ref{Slope_defi}, and $\langle -, -\rangle$ is the natural pairing between characters and cocharacters.} Assume furthermore that $p \ge 5$ and that $\sF^{wa}_C$ is nonempty.\footnote{By \cite[Thm 9.5.10]{dat_orlik_rapoport_2010}, the non-emptiness is equivalent to $b$ being acceptable for $\{ \mu \}$, in the sense of \cite[Section 2.2]{rapoport2014theorylocalshimuravarieties}.}
    
    Then, there is a $J_b(\Q_p) \times \G_{E_s}$-equivariant isomorphism : $$ H^\bullet_{\e t}(\sF^{wa}_C, \F_p(d)) \simeq \bigoplus_{[w] \in W^\mu/\G_{E_s}} v^{J_b}_{P_{[w]}}(\F_p)^* \otimes \rho_{[w]}(\F_p)^* \, [- 2d + n_{[w]}].$$
\end{corro}

This follows directly from Poincaré duality \ref{PCCS_Period_Dom_Intro} and the computation of compactly supported cohomology in \cite[Thm 1.2]{cdhn}. 

The condition "$p \ge 5$" in the above theorem stems from a technical statement  \cite[Thm 1.8]{cdhn} where it is needed to establish the vanishing on some $\operatorname{Ext}^1$ groups between generalized Steinberg representations. It is expected from the authors of loc. cit. that this hypothesis is likely not needed.

\smallskip
Let us now briefly explain our strategy for the proof of the main theorem \ref{PrimtiveComparisonhasPD}.

\subsection{The proof of the main result}
\label{Proof_subseq}

The proof of  \ref{PrimtiveComparisonhasPD} relies heavily on Mann's six functor formalism with $\O^+/p$-coefficients. 

In his thesis \cite{mann2022-6F}, L.Mann attaches, to every untilted small v-stack $X$ over $\Q_p$, a category $\Dbs(X, \O_X^+/p)$ of solid quasi-coherent sheaves of almost $\O^+_X/p$-modules, which admits a six functor formalism.

He also constructs a fully faithful "Riemann-Hilbert functor" : $$(-) \otimes \O^{+ a}_X/p : \cD_{\e t}(X, \F_p)^{oc} \to \Dbs(X, \O^+_X/p),$$ relating his category to the left-completed derived  category of overconvergent étale sheaves of $\F_p$-vector spaces on $X$. From the six-functor formalism, it is relatively straightforward to establish Poincaré duality on the essential image of the Riemann-Hilbert functor, in the following sense :

\begin{thm}(cf. \ref{Corro_PCCS_PD})
    \label{Intro_RealPCCSimpliesPD}
    Let $C$ be a complete algebraically closed extension of $\Q_p$, and $\O_C$ be its ring of integers. Let $f : X \to \operatorname{Spa}(C, \O_C)$ be a smooth locally noetherian analytic adic space that is pure of dimension $d$. 
    
    Let $\cL$ be an overconvergent étale sheaf of $\F_p$-vector spaces on $X$.  Assume that there exists a complex $\cF_\cL \in \cD_{\e t}(\operatorname{Spa}(C, \O_C), \F_p)$ together with an isomorphism in $\Dbs(\operatorname{Spa}(C, \O_C), \O^+_{\operatorname{Spa}(C, \O_C)}/p)$ : $$f^{Mann}_! (\cL \otimes\O^{+ a}_X/p) \simeq \cF_\cL \otimes \O^{+a}_{\operatorname{Spa}(C, \O_C)}/p.$$
    
    Then, for any $0 \le k \le 2d$, there exists a natural isomorphism : $$H^k_{\e t}(X, \cL^\vee)(d) \simeq \Hom_{\F_p}(\pi_{k - 2d}(\cF_\cL), \F_p),$$ where $\pi_{k - 2d}(\cF_\cL)$ denotes the $k-2d$-th homology group of the complex $\cF_{\cL}$. 
\end{thm}

Our key technical statement is that, for a partially proper $f$ as above, the complex $f_!^{Mann}(\cL \otimes \O_X^{+a}/p)$ computes the compactly supported cohomology $R\Gamma_{\e t,c}(X, \cL \otimes \O^+_X/p)$ up to an almost isomorphism, so that it suffices to verify primitive comparison with compact support in the sense of Theorem \ref{PrimtiveComparisonhasPD} - in which case the condition of \ref{Intro_RealPCCSimpliesPD} is satisfied with $\cF_{\cL} = R\Gamma_{\e t, c}(X, \cL)$. 

The precise statement is as follows :

\begin{thm} (cf. \ref{Identification_lower_shriek})
    \label{Intro_RealPCCSequivtoNaive} 
     Let $C$ be a complete algebraically closed extension of $\Q_p$, and $\O_C$ be its ring of integers. Let $f : X \to \operatorname{Spa}(C, \O_C)$ be a smooth and partially proper rigid analytic variety over $C$, and $\cL$ be an overconvergent étale sheaf of $\F_p$-vector spaces on $X$.

    Then $f_!^{Mann}(\cL \otimes \O_X^+/p)$ is a discrete object of $\Dbs(\operatorname{Spa}(C, \O_C), \O_{\operatorname{Spa}(C, \O_C)}^+/p)$.
    
    Moreover, there is an natural isomorphism of almost modules : $$f_!^{Mann} (\cL \otimes \O^+_X/p) \cong R\Gamma_{\e t, c}(X, \cL \otimes \O^+_X/p)^a,$$
    
    where we implicitly use the identification of \cite[Prop 3.3.16]{mann2022-6F}.
\end{thm}

The "discreteness" is to be taken in the sense of condensed mathematics - discrete objects are the ones with trivial condensed structure.

Hence the main theorem \ref{PrimtiveComparisonhasPD} directly follows from \ref{Intro_RealPCCSimpliesPD} and \ref{Intro_RealPCCSequivtoNaive}. The proof of \ref{Intro_RealPCCSequivtoNaive} is the technical heart of this paper, and occupies most of section \ref{Sec4}. The gist of the proof consists in establishing localization triangles of the form : $$R\Gamma_c(X, \cL \otimes \O^+_X/p) \to R\Gamma(X, \cL \otimes \O^+_X/p) \to \rlim_{U \in \Phi_X} R\Gamma(X \setminus U, i_{X \setminus  U}^* \cL \otimes \O^+_{X \setminus U}/p), $$ where the colimit is taken along quasicompact open subsets $U \subset X$, in both Mann's six functor formalism and an appropriate almost quasi-coherent version of the formalism of standard étale cohomology. 

Using the five lemma, we reduce to a variant of the second part of \ref{Intro_RealPCCSequivtoNaive} for $f_*$ instead of $f_!$, for which we establish the following :

\begin{pro} (cf. \ref{Disc_OC_over_field})
    In the setup of \ref{Intro_RealPCCSequivtoNaive}, there is a natural isomorphism $$\left(f_*^{Mann}(\cL \otimes \O^{+,a}/p)\right)_\omega \cong R\Gamma_{\e t}(X, \cL \otimes \O_X^+/p)^{a}$$
    through the natural identification of \cite[Prop 3.3.16]{mann2022-6F}, where $(-)_\omega$ denotes the discretization functor of \cite[Defi 3.2.20]{mann2022-6F}.
\end{pro}

This is a variant of \cite[Example 3.3.17]{mann2022-6F}. Note that the functor $f^{Mann}_*$ does not preserve discrete objects in this setup.

\subsection{Outline of the paper}

In section \ref{Sec2}, after a recollection about the sheaf $\O^+/p$ on locally spatial diamonds, we briefly present Mann's six functors formalism and introduce the results that will be needed later. We also extract from the formalism a version of Poincaré duality on the essential image of the Riemann-Hilbert functor, as stated in \ref{Intro_RealPCCSimpliesPD}.

In section \ref{Sec3}, we provide technical topological results on locally spatial diamonds and partially proper analytic adic spaces, that will be used in section \ref{Sec4}. We also establish a localisation sequence in classical étale cohomology.

Section \ref{Sec4} is the technical heart of the paper, and is dedicated to the proof of theorem \ref{Intro_RealPCCSequivtoNaive}, following the path exposed in section \ref{Proof_subseq}.

In section \ref{Sec5}, we deduce Poincaré duality for almost proper spaces, after a recollection on the primitive comparison theorem.

Section \ref{Sec6} is an introduction to the theory of period domains and their geometry, and, in section \ref{Sec7}, we establish primitive comparison with compact support for such spaces.

\subsection{Conventions and notations}

In what follows, $p$ denotes a fixed prime number, and $\N$ denotes the set of nonnegative integers.

For an abelian category $\mathcal{A}$ use the upright $D(\mathcal{A})$ to denote its $1$-derived category, and $\cD(\mathcal{A})$ to denote its associated $\infty$-derived category. We work with homotopical notation, and, for $\cC \in \cD(\mathcal{A})$, we let $\pi_k(\cC)$ denotes its $k$-th homology group.

Recall that, in a stable $\infty$-category, a square :

\begin{center}
\begin{tikzcd}
 & X \ar[r] \ar[d] & Y \ar[d] \\ 
 & 0 \ar[r] & Z
\end{tikzcd}
\end{center}

is cartesian if and only if it is cocartesian. When this is the case, we say that $X \to Y \to Z$ forms an exact sequence (this corresponds to a distinguished triangle on the homotopy category). Moreover, a reduced functor between stable infinity categories preserves fiber sequences if and only if it preserves cofiber sequences, and, when this is the case, we say that the functor is exact.

Using the natural embeddings, we canonicaly identify rigid analytic varieties with their associated analytic adic space (cf. \cite[1.1.11.(d)]{Huber_Ecoh}), analytic adic spaces with their associated locally spatial diamond (cf. \cite[Defi 15.5]{Scholze_ECohD}), and locally spatial diamonds with their associated small v-stack (cf. \cite[Remark 12.2]{Scholze_ECohD}).

We use the subscript of superscript "M" to denote the functors $f^*, f^M, \dots$, as constructed by Mann in \cite{mann2022-6F} (that was denoted "Mann" in the introduction). We use the letter "H" to denote functors as defined by Huber in \cite{Huber_Ecoh}, and "S" to denote their improvements to the setup of diamonds as constructed in \cite{Scholze_ECohD}. 

Note that Mann's functors are always implicitely derived, whereas, the functors $f^*_{H}, f_!^{H}, f^*_{S}, f_!^{S}, \dots$ are not, and we denote $Rf^*_{H}, Rf_!^{H}, Rf^*_{S}, Rf_!^{S}, \dots$ their derived version.

\subsection{Acknowledgements}

This paper was written during the PhD thesis of the author, under the supervision of Gabriel Dospinescu, Arthur-César Le Bras, and Frédéric Déglise, for which I am very thankful. I also thank Lucas Mann and Nataniel Marquis for some useful discussions.

\section{A Crash Course in Mann's six functor formalism}
\label{Sec2}

In this section, we review some of the main ideas and results of Mann's six functor formalism \cite{mann2022-6F}, focusing on the comparison results with more standard objects, and we prove Poincaré duality on the essential image of the Riemann-Hilbert functor (see \ref{PCCSPD}). All categories considered will be $\infty$-categories. 

\subsection{The structure sheaf $\O^+/p$ for locally spatial diamonds} Note that we could replace $p$ with a fixed pseudo-uniformizer $\pi$ throughout this section. For simplicity, we stick with $\pi = p$.

Recall that any locally spatial diamond $X$ admits both an étale site $X_{\e t}$ and a quasi pro-étale site $X_{qpro\e t}$, defined in \cite[Def 14.1]{Scholze_ECohD}. There is a natural morphism of sites $\nu_X : X_{qpro\e t} \to X_{\e t}$, which induces a geometric morphism of topos given by the pushforward $\nu^X_{*} : X_{qpro\e t}^\sim \to X_{\e t}^\sim$ and a pullback $\nu_X^* : X_{\e t}^\sim \to X_{qpro\e t}^\sim$.

Recall that any analytic adic space $X$ over $\Z_p$ can naturally be viewed as a locally spatial diamond $X^\diamond$ that admits an homeomorphic underlying topological space, as well as an isomorphic étale topos to the one of $X$ (cf. \cite[Lemma 15.6]{Scholze_ECohD}).

A process described in \cite[Section 2]{MannWerner2020localsystemsdiamondspadic} endows every diamond $X$ over $\operatorname{Spd}(\Z_p)$ with quasi-pro-étale structure sheaves, $\widehat{\O}_X$ and $\widehat{\O}_X^+$, and, hence every diamond $X$ over $\operatorname{Spd}(\Q_p)$ with a pro-étale sheaf $\widehat{\O}_X^+/p$, viewed as the cokernel of the multiplication by $p$. This sheaf can be constructed as the pro-étale sheafification of the presheaf defined on the basis of perfectoid spaces over $X$ by $S \mapsto \O^+_{S^\sharp}(S^\sharp)/p$, where $S^\sharp$ denotes the untilt of $S$ associated to the composition $S \to X \to \operatorname{Spd}(\Q_p)$. Note that it is important to remember the structure morphism to $\operatorname{Spd}(\Q_p)$, as, otherwise, a lot of information is lost.

As is shown in \cite[Lemma 2.7]{MannWerner2020localsystemsdiamondspadic}, if $X$ is a locally noetherian analytic space over $\operatorname{Spa}(\Q_p, \Z_p)$, the restriction $\nu_X^* \hat{\O}_X^+/p$ to the étale site simply coincides with the classical étale structure sheaf $\O^+_X/p$, viewed as the mod $p$ reduction of the étale sheaf $Y \in X_{\e t} \mapsto \O^+_Y(Y)$. We also refer the reader to \cite[Section 6]{zavyalov2024coherentmodulescoherentsheaves} for a detailed analysis the of étale, pro-étale, and v-structure sheaves in the setup of pre-adic spaces.

In the following, we show that, if $X \to \operatorname{Spd}(\Q_p)$ is a locally spatial diamond, the pro-étale sheaf $\hat{\O}_X^+/p$ comes from the pullback of an étale sheaf, that we may then reasonably call $\O_X^+/p$.\footnote{That this is claimed in the proof of \cite[Theorem 25.1]{Scholze_ECohD}, but, as far as I know, the proof of this statement does not appear explicitly in the litterature.}

Recall the following lemma : 

\begin{lemma}
    \label{pro_et_pullback_is_restr}
    Let $f : X \to Y$ be a quasi-pro-étale morphism of locally spatial diamonds, and $\cF \in Y_{qpro\e t}$. Then, the pullback $f^* \cF$ is given by the restriction of $\cF$ to $X_{qpro\e t}$, i.e, for any quasi-pro-étale $U \to Y$, $f^* \cF(U) = \cF(U \to X)$, where we consider the quasi-pro-étale composition $U \to Y \to X$.
\end{lemma}

The same result holds by replacing "quasi-pro-étale" by "étale" everywhere.

\begin{proof}
    By the general formalism (cf. \cite[\href{https://stacks.math.columbia.edu/tag/00VC}{Section 00VC}]{stacks-project}, or \cite[\href{https://stacks.math.columbia.edu/tag/03PZ}{Section 03PZ}]{stacks-project} in the étale case), the inverse image $f^* \cF$ is given by the sheafification of the following : $$ U/X \mapsto \rlim_{(V, \varphi) \in (I_U)^{op}} \cF(V/Y),$$ where $I_U$ denotes the category of pairs $(V, \varphi)$ defined by : 
    \begin{itemize}[itemsep = .1pt]
        \item $\operatorname{Ob}(I_U)= \{ (V, \varphi) : V \in Y_{qpro\e t}, \varphi : U \to V \times_X Y\}$
        \item $\operatorname{Hom}((V, \varphi), (V', \varphi')) = \{ g : V \to V' \text{ quasi-pro-étale }, g_Y \circ \varphi = \varphi'\}$ where $g_Y$ denotes the base change $V \times_X Y \xrightarrow{g \times_X id_{Y}} V' \times_X Y$ of $g$.
    \end{itemize}
    Since $f$ is quasi-pro-étale, we may view any $U/X$ as quasi-pro-étale over $Y$, so that we may view $(U, \varphi_U : U \to U \times_X Y)$ as an object of $I_U$. Moreover, for any $(V, \varphi_V) \in I_U$, there is a natural morphism $g : U \to V$ given by the composition $U \to  V \times_X Y \to V$, that is quasi-pro-étale by the two-out-of-three property \cite[Prop 10.4]{Scholze_ECohD} applied to the composition $U \to V \to Y$.
    
    It is straightforward to check that this induces a morphism $(U, \varphi_U) \to (V, \varphi_V)$ in $I_U$, so that $(U, \varphi_U)$ is an initial object of $I_U$, i.e. a final object of $(I_U)^{op}$. Hence, the colimit alongside $(I_U)^{op}$ simply computes $\cF(U \to Y)$, which clearly defines a sheaf, so that there is no need to sheafify. This concludes.
\end{proof}

\begin{pro}
    \label{iset}
    Let $X$ a locally spatial diamond over $\operatorname{Spd}(\Q_p)$. 

    Then, there exists an étale sheaf on $X_{\e t}^\sim$, denoted $\O^+_X/p$, such that $\widehat{\O}^+_X/p \cong \nu_X^* \O^+_X/p$. 
\end{pro}

\begin{proof}
    Since this property is Zariski-local over $X$, we may assume that $X$ is spatial, so that there exists a surjective quasi-pro-étale morphism $f : Y \to X$ from a strictly totally disconnected perfectoid space, by \cite[Pro 11.24]{Scholze_ECohD}. 
    
    By definition, $\hat{\O}_X^+$ is such that, for any morphism $Z \to X$ from a totally disconnected affinoid perfectoid space, $\hat{\O}_X^+(Z) = \O^+_{Z^\sharp}(Z^\sharp)$, where $Z^\sharp$ is the untilt of $Z$ viewed through the composition $Z \to X \to \operatorname{Spd}(\Z_p)$ - there is no need to sheafify, using \cite[Theorem 8.7]{Scholze_ECohD}. This construction is clearly compatible with restrictions, so that, by the above lemma, $f^* \hat{\O}_X^+$ is simply $\hat{\O}^+_Y$. Since pullbacks commute with small colimits, the same holds for $\hat{\O}/p$. 
    
    Since $Y$ is perfectoid, it follows from the proof of \cite[Lemma 2.7]{MannWerner2020localsystemsdiamondspadic} that the pro-étale sheaf $\hat{\O}^+_Y$ actually comes from the pullback of an étale sheaf, and so does $\hat{\O}^+_Y/p$. Hence, by \cite[Thm 14.12]{Scholze_ECohD}, the sheaf $\hat{\O}_X^+/p$ on $X$ comes from the pullback of an étale sheaf, that we denote $\O_X^+/p$.
    \end{proof}

We'll also use the following lemma.

\begin{lemma}
    \label{Pullback_qproet}
    Let $f : X \to Y$ be a quasi-pro-étale morphism of locally spatial diamonds over $\operatorname{Spd}(\Q_p)$. 
    
    Then $\O^+_Y/p \cong f_{\e t}^* \O^+_{X}/p$, where $f_{\e t}^*$ denotes the pullback of étale sheaves along $f$.
\end{lemma}

\begin{proof} We start by establishing the commutativity of the following diagram :
    \begin{center}
        \begin{tikzcd}
            &X_{\e t}^\sim \ar[d, swap, "\nu_X^*"] 
            & Y_{\e t}^{\sim} \ar[d, "\nu_Y^*"] \ar[l, "f_{\e t}^*", swap] \\ 
            & X_{qpro\e t}^\sim 
            & Y_{qpro\e t}^{\sim} \ar[l,"f_{qpro\e t}^*", swap ]
        \end{tikzcd}            
    \end{center}

    All of the four arrows in this diagram are induced from morphisms of sites, which are themselves induced by continuous functors satisfying the conditions of \cite[\href{https://stacks.math.columbia.edu/tag/00X6}{Proposition 00X6}]{stacks-project}. Hence, if suffices to consider the following diagram :\footnote{Here, the arrows are continuous functors, and the associated morphisms of site are contravariant.}
    
    \begin{center}
    \begin{tikzcd}[column sep = 2cm]
        & X_{\e t}  \ar[d, "id", swap] 
        & Y_{\e t}  \ar[d, "id"] \arrow[dl, dashed] \ar[l,swap, "Z \mapsto Z \times_X Y"]\\ 
        & X_{qpro\e t} 
        & Y_{qpro\e t} \ar[l, "Z \mapsto Z \times_X Y"]
    \end{tikzcd}            
    \end{center}

    Here, the vertical arrows are obtained by viewing any étale cover as a pro-étale one, so that the diagram obviously commutes. It follows from functoriality of morphisms of sites \cite[\href{https://stacks.math.columbia.edu/tag/03CB}{Lemma 03CB}]{stacks-project} that both compositions $\nu_Y^* \circ f_{\e t}^*$ and $f_{qpro\e t}^* \circ \nu_X^*$ are equal to the pullback along the morphism of site induced by the diagonal $X_{\e t} \to Y_{qproet}$, so that the first diagram commutes.

    \medskip
    Let us now prove that $ f^*_{\e t} \O^+_Y/p \cong \O^+_X/p$. From the lemma \ref{pro_et_pullback_is_restr} and the proof of \ref{iset}, it follows that $ f^*_{pro \e t} \hat{\O}^+_Y/p \cong \hat{\O}^+_X/p$.
    By commutation of the above diagram, we have $$\nu_X^* \O^+_X/p \cong \hat{\O}^+_X/p \cong f^*_{qpro\e t} \nu_Y^* \hat{\O}^+_Y/p \cong \nu_X^* f^*_{\e t} \O^+_Y/p$$
    
    Finally, by \cite[Prop 14.8]{Scholze_ECohD}, $\nu^X_* \nu_X^* \cong id$, so that $\O^+_X/p \cong f^*_{\e t} \O^+_Y/p$. 
\end{proof}

\subsection{Mann's functors $\Dbs(\rule{.2cm}{.5pt}, \O^+/\pi)$ and $\Dbs(\rule{.2cm}{.5pt}, \O^+/\pi)^\varphi$}

\label{Manndef}

Mann associates to any untilted small v-stack $X$, with a fixed pseudouniformizer $\pi$, an $\infty$-category $\Dbs(X, \O^+_X/\pi)$ consisting of \textit{complexes} of \textit{almost, quasi-coherent, solid} sheaves of $\O^+_X/\pi$-modules on $X$. 

For an analytic adic space $X$, objects of $\Dbs(X, \O^+_X/\pi)$ are not explicit - rather, this category is constructed by gluing from an adequate basis of the pro-étale topology, given by (a loosely improved class of) totally disconnected perfectoid spaces, as follows : 

\begin{defi}
    \label{Dbsperf}
    Let $X = \Spa(A, A^+)$ be an affinoid perfectoid space that is of weakly perfectly finite type over some totally disconnected space (for example : the spectrum of a perfectoid field), and let $\pi$ be a pseudo-uniformizer of $A$.
    
    Then, let $\Dbs(X, \O^+_X/\pi) := \Dbs(A^+/\pi)$ be the derived $\infty$-category of solid almost $A^+/\pi$-modules, constructed as follows : 
    \begin{itemize}
        \item Let $\cD(A^+/\pi)^{cond}$ denote the $\infty$-derived category of the category of condensed $A^+/\pi$-modules.
        \item Let $\cD_\square(A^+/\pi) \subset \cD(A^+/\pi)^{cond}$ be the full subcategory formed by solid objects, with respect to $A^+/\pi$. 
        \item We say that a morphism $f : M \to N$ in $\cD_\square(A^+/\pi)$ is an almost isomorphism if, in $\operatorname{fib}(f)$, the multiplication by any pseudo-uniformizer is zero. 
        
        Finally, let $\cD_\square^a(A^+/\pi)$ be the localization of $\cD_\square(A^+/\pi)$ with respect to almost isomorphisms.
    \end{itemize}
\end{defi}

An object of $\cD_\square^a(A^+/\pi)$ can be thought of a complex of solid (i.e. "complete" for some topology) $A^+/\pi$-modules, up to an almost isomorphism, in the sense of Faltings. 

\begin{rk}
More generally, Mann defines and studies analogue over animated rings. In this more general setup, one replaces the derived category by the category of "modules" over the animated ring, in the sense of Lurie's higher algebra \cite[Def 7.1.1.2]{Lurie_HA}. In this paper, we always consider static (i.e. not animated) rings, so that this is simply a derived category by \cite[Remark 7.1.1.6]{Lurie_HA}. 
\end{rk}


Recall that small v-stacks, defined in \cite{Scholze_ECohD}, form a generalization of analytic adic spaces, using a stacky approach, so that they are best used in descent problems. In order to properly endow a small v-stack $X$ with a sheaf $\O^+_X/\pi$, we need to chose an untilt\footnote{Note that, contrary to what the notation might suggest, the sheaf $\O^+_X/\pi$ depends on the choice of the untilt $X^\sharp$, and should maybe rather be denoted $\O^+_{X^\sharp}/\pi$.} $X^\sharp$ of $X$, and a pseudo-uniformizer $\pi$. Let $\operatorname{vStack}^\sharp_\pi$ denote the associated category. 

\begin{pro} (\cite[Thm 1.2.1]{mann2022-6F})
    \label{descentM}
    There exists a unique hypercomplete v-sheaf of $\infty$-categories : $$(X, X^\sharp, \pi) \in (\operatorname{vStack}^\sharp_\pi)^{op} \mapsto \Dbs(X, \O^+_X/\pi) \in \operatorname{Cat}_\infty$$ 

    extending the map defined in \ref{Dbsperf}. This construction admits a six-functor formalism. 
\end{pro}

\begin{rk}
    If $X = \operatorname{Spa}(C, \O_C)$ for some complete algebraically closed field $C$ with ring of integers $\O_C$, we simply let $\Dbs(\operatorname{Spa}(C, \O_C), \O^+_C/p) := \Dbs(\operatorname{Spa}(C, \O_C), \O^+_{\operatorname{Spa}(C, \O_C)}/p)$.
\end{rk}

While we do not precisely define what a six functors formalism \textit{is}\footnote{A precise definition was coined in \cite[Definition A.5.7]{mann2022-6F}. See also lectures notes by Scholze \cite{Scholze_6FF}.}, let us present the main characteristic properties, as in \cite[Theorem 1.2.4]{mann2022-6F}.

\begin{pro} 
    \label{6ffprop}
    The mapping\footnote{For simplicity, we simply denote $(X, X^\sharp, \pi)$ as $X$.} $X \mapsto \Dbs(X, \O^+_X/\pi)$ has the following properties :
    \begin{enumerate}[itemsep=.1cm]
        \item For any $X$, the category $\Dbs(X, \O^+_X/\pi)$ is closed symmetric monoidal, i.e. admits a  tensor product $- {\otimes}^M - $, which admits a right adjoint, given by an internal Hom functor $\underline{\Hom}^M$. \\
        For any $X$, the unit of the tensor product over $X$ is $\O^+_X/\pi$.

        \item Any morphism $f : X \to Y$ induces adjoint functors : $f^M_* : \Dbs(X, \O^+_X/\pi) \rightleftarrows \Dbs(Y, \O^+_Y/\pi) : f_M^* $, where the left adjoint $f_M^*$ is called a pullback functor, and $f^M_*$ is a pushforward. Both are compatible with composition.

        \item Any bdcs\footnote{Bdcs morphisms are defined as locally compactifiable (i.e written locally as the composition of an étale and a proper map) and "p-cohomologically bounded". This is a technical assumption that is satisfied by every morphism of rigid analytic varieties, and by any morphism of analytic adic spaces that is locally weakly of finite type, cf. \cite[Prop 3.5.14]{mann2022-6F}.} morphism $f : X \to Y$ induces an adjoint pair : $f^M_! : \Dbs(X, \O^+_X/\pi) \rightleftarrows \Dbs(Y, \O^+_Y/\pi) : f^!_M$, where the left adjoint $f^M_!$ is called a proper pushforward, and $f^M_*$ is an exceptional pullback. Both are compatible with composition.Both are compatible with composition of bdcs morphisms.

        \item Let $f : X \to Y$ be a bdcs morphism. If $f$ is proper, then $f^M_! = f^M_*$. If $f$ is étale, then $f_M^! = f_M^*$.

        \item (Projection formula) Let $f : X \to Y$ be a bdcs morphism, $\mathcal{M} \in \Dbs(Y, \O_Y/\pi)$ and $\mathcal{N} \in \Dbs(X, \O_X/\pi)$. Then there exists a natural isomorphism : $$ f^M_! (\mathcal{N} \otimes f_M^*(\mathcal M) ) \cong f^M_! \mathcal{N} \otimes \mathcal M.$$

        \item (Proper base change) Consider the following cartesian square : 

            \begin{center}
            \begin{tikzcd}
                X' \arrow[r, "f'"] \arrow[d, "g'"] & Y \arrow[d, "g"]
                \\ X \arrow[r, "f"] & Y'
            \end{tikzcd}                            
            \end{center}

        Assume that $f$ is bdcs. Then, there is a natural isomorphism of functors : $$ g_M^* \circ f^M_! \cong f'^{M}_! \circ g_M'^{*}.$$

        \item \label{Smoothness} Let $f : X \to Y$ be a smooth morphism of locally noetherian analytic spaces adic over $\Q_p$ of equidimension $d$. Then, there is a canonical isomorphism $f_M^! \O^+_Y/\pi \cong f_M^* \O^+_Y/\pi (d)[2d] $.
    \end{enumerate}
\end{pro}

As is classical whenever we have a six functor formalism, we may define a notion of cohomology (resp. cohomology with compact support) of some $X$ as the pushforward $f^M_*$ (resp $f^M_!$) of the unit object $\O^+_X/\pi$, where $f : X \to \operatorname{Spa}(C, \O_C)$ is a structural morphism to a point. This should be interpreted as a form of "almost étale $\O_X^+/\pi$-cohomology", that should, at least in some cases, coincide, up to an almost isomorphism, with the standard cohomology groups of the étale sheaf $\O_X^+/\pi$ (cf. \cite[Ex. 3.3.17]{mann2022-6F} and \ref{Identification_lower_shriek} for precise such statements). From the properties above, we may formally deduce some nice properties for this version of cohomology, such as the Künneth formula (for cohomology with compact support) and Poincaré duality (see \cite[Ex. 1.0.1 and 1.0.2]{Mann_Heyer_Smooth}).

In order to descend from $\O^+/\pi$-cohomology groups to some $\F_p$-cohomology groups, we need to take some form of Frobenius invariants. Assuming that $\pi \mid p$, a reasonable approach is to start from the Artin-Schreier exact sequence for sheaves on $X$ : $$0 \to \F_p \to \O^+_X/\pi \xrightarrow{x \mapsto x^p - x} \O^+_X/\pi \to 0$$

And consider some form of associated long exact sequence in cohomology (or, rather, a derived refinement).

The key technical difficulty is Mann's category only captures the information of the $\O^+/\pi$-cohomology up to an \textit{almost} isomorphism, which need not commute with taking Frobenius invariants. This is not an artifact of Mann's formalism, as, in practice, most computations of $\O^+$ or $\O^+/\pi$-cohomology rely on the fundamental result of \textit{almost} vanishing of cohomology affinoid perfectoid spaces \cite[Lemma 8.8]{Scholze_ECohD}, which only holds in the almost sense. 

This idea can be adapted in a less naive way. Whenever $\pi \mid p$, the unit object $\O^+_X/\pi \in \Dbs(X, \O^+_X/\pi)$ is naturally endowed with a Frobenius endomorphism $\varphi$, which endows the whole category $\Dbs(X, \O^+_X/\pi)$ with a Frobenius endomorphism $\varphi^*_M$. 


\begin{defi}
    A $\varphi$-module in $\Dbs(X, \O_X^+/\pi)$ is the datum of an object $C \in \Dbs(X, \O^+_X/\pi)$, together with a specified isomorphism $\varphi_M^* C \simeq C$. We let

\vspace{-.2cm}
$$\Dbs(X, \O_X^+/\pi)^\varphi = eq\left(
\begin{tikzcd}
\Dbs(X, \O_X^+/\pi) \ar[r,shift left=.75ex,"\varphi_M^*"]
  \ar[r,shift right=.75ex,swap,"id"]
&
\Dbs(X, \O_X^+/\pi) 
        \end{tikzcd}\right)$$

be the category of $\varphi$-modules in $\Dbs(X, \O^+_X/\pi)^\varphi$.
\end{defi}

\begin{rk}
    \label{Forget_phiinv}
    Note that $\Dbs(X, \O_X^+/\pi)^\varphi$ is very far from being a full subcategory of $\Dbs(X, \O_X^+/\pi)$. Indeed, morphisms in $\Dbs(X, \O_X^+/\pi)^\varphi$ should be thought of as \textit{Frobenius-equivariant} morphisms.
    
    There is a canonical "forgetful" functor $\Dbs(X, \O^{+}_X/\pi)^\varphi \to \Dbs(X, \O^{+}_X/\pi)$, that forgets the $\varphi$-module structure.
\end{rk}

 Moreover, the assignment $X \mapsto \Dbs(X, \O^+_X/p)^\varphi$ is particularly well behaved, since it also satisfies a six-functor formalism, and all properties of \ref{6ffprop} are satisfied with $\Dbs(X, \O^+_X/p)^\varphi$ in place of $\Dbs(X, \O^+_X/p)$.



 
\subsection{Classical $\F_p$-étale sheaves and the Riemann-Hilbert adjunction}

\label{Etale_Shf_on_LSD}

Let us start by recalling some facts about the classical theory of étale sheaves in analytic geometry, as considered by Huber in \cite{Huber_Ecoh}, and refined by Scholze in \cite{Scholze_ECohD}.

Fix $S$ a locally noetherian analytic adic space, and some $n \in \N$. To any locally noetherian adic space $X$ over $S$, Huber associates in \cite{Huber_Ecoh} a 1-category $D_{\e t}(X, \Z/n\Z)$, as the derived category of étale $\Z/n\Z$-sheaves on $X$. He shows that, whenever $n$ is invertible in $\O_S^+$ (and when restricting to spaces of finite type over $S$, the setup $D(X_{\e t}, \Z/n\Z)$ admits a six functor formalism.

\begin{rk}
While Huber works 1-categorically, and the notion of "six functor formalism" was not properly formalized yet, his constructions can be translated in an $\infty$-categorical setup, and can be shown to be an \textit{actual} six-functor formalism. This has been checked thoroughly by Zavyalov in \cite[Section 8]{Zavyalov_Adicnotes} (see also \cite[Section 6.1]{Zavyalov_Abstract_PD}). In particular, it satisfies all the properties of \ref{6ffprop}. In what follows, we always work $\infty$-categorically.
\end{rk}

When $n$ is not invertible in $\O_S^+$, most of the constructions still apply, but some key properties are lacking. In particular, there is no general proper base change (see \cite[4]{mann2022-6F} for a crucial counter-example), nor a general projection formula. 

\begin{rk}
    \label{Scholze_vs_H}
    Most of Huber's theory (at least in the $\ell \neq p$ setup) has been generalized by Scholze in \cite{Scholze_ECohD} in the generality of small v-stacks, and the constructions coincide for locally noetherian analytic adic spaces. We'll use both theories, usually preferring the more modern notations of \cite{Scholze_ECohD}.
\end{rk}


Let $X$ be a locally spatial diamond, and $X_{\e t}$ the associated étale site. Following \cite[Prop 14.15]{Scholze_ECohD}, we let $\cD(X_{\e t}, \F_p)$ be the $\infty$-derived category of sheaves of $\F_p$-modules over $X_{\e t}$. This category needs not be left-complete in general\footnote{Recall that a stable $\infty$-category $\cC$ with a t-structure is said to be left-complete if, for all objects $A$ of $\cC$, the natural morphism $A \to \llim_{n \in \Z} \tau^{\ge n} A$ is an isomorphism, where $\tau^{\ge n}$ denotes the truncation. The left completion of $\cC$ is then defined by $\hat{\cC} = \llim_{n \in \Z} \tau^{\ge n} \cC$, which remains a stable $\infty$-category.}, and we denote its left completion $\cD_{\e t}(X, \F_p)$.

\begin{rk}
    If $X$ is of $p$-cohomologically finite dimension, the category $\cD(X_{\e t}, \F_p)$ is already left-complete by \cite[Prop 20.17]{Scholze_ECohD}.
\end{rk}

This is functorial, in the sense that for any morphism $f : X \to Y$ of locally spatial diamonds, there is an adjoint pair of functors, constructed by Scholze\footnote{If $X$ and $Y$ come from locally noetherian analytic adic spaces, they coincide with the ones defined by Huber : $R f_*^{S} = R f_*^{H}$ and $R f_{S}^* = R f_{H}^*$.} : $$ R f_*^{S} : \cD_{\e t}(X, \F_p) \leftrightarrows \cD_{\e t}(Y, \F_p) : R f_{S}^*$$

Moreover, for any $X$, the category $\cD_{\e t}(X, \F_p)$ admits a well behaved derived tensor product $\otimes^\L_{\F_p}$, which admits a right adjoint, a derived internal Hom, denoted\footnote{This is neither the notation of \cite{Huber_Ecoh} nor \cite{Scholze_ECohD}, but mimics the one of Mann.} $\underline{\operatorname{RHom}}_{\F_p}$. 

Let us now introduce overconvergent sheaves.

\begin{defi}
\label{Overconv}
An étale complex $\cF \in \cD_{\e t}(X, \F_p)$ is said to be \textbf{overconvergent} if, for any quasi-pro-étale morphism $\overline{x} : \operatorname{Spa}(C, C^+) \to X$ from an algebraically closed perfectoid field $C$ with open bounded valuation subring $C^+$, the natural map : $$ \cF_{\overline{x}} \to \cF_{\overline{x^\circ}}$$ is an isomorphism, where $\overline{x^\circ}$ denotes the natural composition $\operatorname{Spa}(C, C^\circ) \to \operatorname{Spa}(C, \O_C) \xrightarrow{\overline{x}} X$, where $C^\circ$ denotes the subring of bounded elements in $C$.

We let $\mathcal{D}(X_{\e t}, \F_p)^{oc} \subset \cD(X_{\e t}, \F_p)$ be the full subcategory of overconvergent étale $\F_p$-sheaves on $X$, and denote $\mathcal{D}_{\e t}(X, \F_p)^{oc}$ its left completion. It is stable under pullbacks and (derived) tensor product.
\end{defi}


Let us now introduce the Riemann-Hilbert and $\varphi$-invariants functor, as well as a few of their properties. Both of these functors are defined by descent, so that we start from the local situation. Let us recall the construction from \cite[Prop. 3.9.8]{mann2022-6F}.

\medskip
Let $X = \operatorname{Spa}(A, A^+)$ be a strictly totally disconnected perfectoid space with pseudo-uniformizer $\pi \in A^+$. On such a space, étale covers split, so that is a natural equivalence of categories : $$\cD_{\e t}(X, \F_p)^{oc} \cong \cD(\F_p(X))$$ with the derived category of $\F_p(X)$-modules, where $\F_p(X)$ denotes the ring of continuous maps $|X| \to \F_p$.

The natural morphism of étale sheaves $(\F_p)_X \to \O_X^+/\pi$ induces a map : $$\F_p(X) \cong H^0(X, (\F_p)_X) \to H^0(X, \O_X^+/\pi) \cong A^+/\pi$$ where the last isomorphism follows from the vanishing of higher étale cohomology groups on strictly totally disconnected perfectoid spaces. This map induces a morphism of analytic rings $\F_p(X) \to A^{+,a}/\pi$.

\begin{defi}
    \label{local_RH} 
    The Riemann-Hilbert functor $\cD(\F_p(X))_\omega \to \Dbs(X, A^+/\pi)^\varphi$ is defined as the base change $(-) \otimes_{\F_p(X)} A^{+a}/\pi$ along the morphism of analytic rings constructed above, in the sense of \cite[Defi 2.12.23.(c)]{mann2022-6F}, with the $\sigma$-structure on $A^{+,a}/\pi$ given by the Frobenius morphism.
\end{defi}

This local construction glues to the following : 

\begin{defi}
    \label{RH_def}
    Let $K$ be a complete extension of $\Q_p$, with pseudo-uniformizer $\pi$ such that $\pi | p$. For $X$ an untilted locally spatial diamond over $\operatorname{Spd}(K)$, Mann defines \cite[Defi 3.9.21]{mann2022-6F} : 
    \begin{enumerate}
        \item A "Riemann-Hilbert" functor $- \otimes\O^{+a}_X/\pi : \cD_{\e t}(X, \F_p)^{oc} \to \Dbs(X, \O_X^{+ a}/\pi)^\varphi$
        \item Its left adjoint, a functor of $\varphi$-invariants : $(-)^\varphi : \Dbs(X, \O^+_X/\pi)^\varphi \to \cD_{\e t}(X, \F_p)^{oc}$
    \end{enumerate}
    
    When $X = \operatorname{Spa}(K, \O_K)$ is the spectrum of an complete non-archimedean field with pseudo-uniformizer $\pi$, we simply denote the Riemann-Hilbert functor as $- \otimes \O^{+a}_K/\pi$. 
\end{defi}

\begin{rk}
    This Riemann-Hilbert functor can be seen as an analogue of the one introduced by Bhatt-Lurie in \cite{bhatt2017riemannhilbertcorrespondencepositivecharacteristic}, relating $\F_p$-étale sheaves over a fixed scheme $X$ of characteristic $p$ with quasicoherent modules over $X$ equipped with a Frobenius.
\end{rk}

The key fact for our purposes is as follows :

\begin{pro}(\cite[Thm 3.9.23]{mann2022-6F}) The Riemann-Hilbert functor is fully faithful.
\end{pro}

Moreover, when restricting both sides of the equivalence to the full subcategories of dualizable objects, L.Mann proves that the Riemann-Hilbert functor induces an equivalence of category. In \cite{mann2022-6F}, this is the result used by Mann to establish his form of Poincaré duality. Dualizable objects of $\cD_{\e t}(X, \F_p)^{oc}$ can be seen to correspond to \textit{perfect complexes}, i.e. the ones that are locally quasi-isomorphic to bounded complexes of finite dimensional.

In this paper, we're interested in the étale cohomology of non-proper spaces, for which cohomology groups may be infinite dimensional, and, hence, non dualizable. We'll nonetheless prove that the Poincaré duality holds on the essential image of the Riemann-Hilbert functor.

\subsection{Poincaré Duality on the essential image of Riemann-Hilbert}
    \label{TechLemma}

For any $\cF \in \cD_{\e t}(X, \F_p)$, we define its dual $\cF^{\vee} := \underline{\operatorname{RHom}}_{\cD_{\e t}(X, \F_p)}(\cF, (\F_p)_X)$, where $(\F_p)_X$ denotes the constant sheaf associated to $\F_p$ on $X_{\e t}$. We establish the following : 

    \begin{thm}
    \label{PCCSPD}
    Let $K$ be a complete extension of $\Q_p$, $\O_K$ be its ring of integers, and let $\pi$ be a pseudo-uniformizer, such that $\pi | p$. Let $f : X \to \operatorname{Spa}(K, \O_K)$ be a smooth locally noetherian analytic adic space that is pure of dimension $d$. Let $\cL \in \cD_{\e t}(X, \F_p)^{oc}$.

    Assume that there exists $\cF_\cL \in \cD_{\e t}(\operatorname{Spa}(K, \O_K), \F_p)$ and an isomorphism in $\Dbs(\operatorname{Spa}(K, \O_K), \O^+_K/\pi)$ : $$f^{M}_{!} (\cL \otimes\O^{+ a}_X/\pi) \simeq \cF_\cL \otimes \O^{+a}_K/\pi.$$ 
    Then, there is an isomorphism in $\cD_{\e t}(\operatorname{Spa}(K, \O_K), \F_p)$ : $$Rf_{*}^{H} \cL^{\vee}(d)[2d] \simeq \cF_\cL^\vee. $$
    \end{thm}

    Note that, if $f$ is additionally proper and $\cL$ is a perfect (e.g. a local system), the primitive comparison theorem \cite[Corro 3.9.24]{mann2022-6F} applies, so that the assumption is automatically verified with $\cF_\cL = Rf_*^{M} \cL$ (this is overconvergent since $f$ is quasi-compact). In that setup, this recovers \cite[Corro 3.10.22]{mann2022-6F}.

    \begin{proof}
    By taking stalks from $\operatorname{Spa}(\overline{K}, \O_{\overline{K}})$, we may freely identify $\cD_{\e t}(\operatorname{Spa}(K, \O_K), \F_p)$ with the derived $\infty$-category of $\F_p$-vector spaces equipped with a $\operatorname{Gal}(\overline{K}/K)$-action, with its natural symmetric monoidal structure. In particular, all étale sheaves on $\operatorname{Spa}(K, \O_K)$ are necessarily overconvergent.
    
    The Riemann-Hilbert functor is symmetric monoidal and commutes with pullbacks. Hence, for any $A \in \cD_{\e t}(\operatorname{Spa}(K, \O_K), \F_p) = \cD_{\e t}(\operatorname{Spa}(K, \O_K), \F_p)^{oc}$, we have :
    \begin{align*}
        f^{M}_!\left((\cL \otimes^\L_{\F_p} Rf^*_{H} A) \otimes \O^{+ a}_X/\pi \right) & \cong f^{M}_!\left( (\cL \otimes \O^{+ a}_X/\pi) \otimes^M (Rf^*_{H} A \otimes \O^{+ a}_X/\pi) \right) \\
        & \cong f^{M}_! \left((\cL \otimes \O^{+ a}_X/\pi) \otimes^M f^*_{M} (A \otimes \O^{+ a}_K/\pi)\right)  \\
        & \cong \left(f^{M}_! (\cL \otimes \O^{+ a}_X/\pi)\right) \otimes^M (A \otimes \O^{+ a}_K/\pi) \\ 
        & \simeq (\cF_\cL \otimes \O^{+a}_K/\pi) \otimes^M (A \otimes \O^{+ a}_K/\pi) \\ 
        & \cong (\cF_\cL \otimes^\L_{\F_p} A) \otimes \O^{+ a}_K/\pi.
    \end{align*}

    where the third isomorphism follows from the projection formula (cf. fifth point of \ref{6ffprop}), and the fourth one follows from the hypothesis. From there, we compute :
    \begin{align*}
        & \operatorname{Hom}_{\cD_{\e t}(\operatorname{Spa}(K, \O_K), \F_p)}\left(A, Rf_{*}^{H} \cL^\vee (d)[2d]\right) & \\ 
        & \cong \operatorname{Hom}_{\cD_{\e t}(X, \F_p)}\left(Rf^{*}_{H} A, \underline{\operatorname{RHom}}_{\F_p} (\cL, \F_p) \, (d)[2d]\right) & \text{By adjunction}\\
        & \cong \operatorname{Hom}_{\cD_{\e t}(X, \F_p)^{oc}}\left(Rf^{*}_{H} A \otimes^\L_{\F_p} \cL, \F_p \, (d)[2d]\right) & \text{By adjunction} \\ 
        & \cong \operatorname{Hom}_{\Dbs(\O^{+a}_X/\pi)^\varphi}\left( (Rf^{*}_{H} A \otimes^\L_{\F_p} \cL) \otimes \O^{+a}_X/\pi, \O^{+a}_X/\pi \, (d)[2d]\right) & \text{By full faithfulness of $(-) \otimes \O^{+a}_X/p$} \\
        & \cong \operatorname{Hom}_{\Dbs( \O^{+a}_X/\pi)^\varphi}\left((Rf^{*}_{H} A \otimes_{\F_p}^\L \cL) \otimes \O^{+a}_X/\pi, f^!_{M} \O^{+a}_K/\pi\right) & \text{By the last point of \ref{6ffprop}} \\
        & \cong \operatorname{Hom}_{\Dbs( \O^{+a}_X/\pi)^\varphi}\left(f_!^{M}\left( (Rf^{*}_{H} A \otimes_{\F_p}^\L \cL) \otimes \O^{+a}_X/\pi\right), \O^{+a}_K/\pi \right) & \text{By adjunction}\\
        & \simeq \operatorname{Hom}_{\Dbs( \O^{+a}_X/\pi)^\varphi}\left((\cF_\cL \otimes_{\F_p}^\L A) \otimes \O^{+ a}_K/\pi,  \O^{+a}_K/\pi \right) & \text{By the above} \\ 
        & \cong \operatorname{Hom}_{\cD_{\e t}(X, \F_p)^{oc}}\left(\cF_\cL \otimes_{\F_p}^\L A, \F_p\right) & \text{By fully faithfulness $(-) \otimes \O^{+a}_K/p$} \\
        & \cong \operatorname{Hom}_{\cD_{\e t}(X, \F_p)}\left(A, \underline{\operatorname{RHom}}_{\F_p} (\cF_\cL, \F_p)\right). & \text{By adjunction}
    \end{align*}
    
    We conclude using Yoneda's lemma.
\end{proof}
    
The following corollary will be the version used in practice.
    
\begin{corro}
    \label{Corro_PCCS_PD}
    Let $C$ be a complete algebraically closed extension of $\Q_p$, and $\O_C$ be its ring of integers. Let $f : X \to \operatorname{Spa}(C, \O_C)$ be a smooth locally noetherian analytic adic space that is pure of dimension $d$. 
    
    Let $\cL$ be an overconvergent étale sheaf of $\F_p$-vector spaces on $X$. Assume that there exists some $\cF_\cL \in \cD_{\e t}(\operatorname{Spa}(C, \O_C), \F_p)$ and an isomorphism in $\Dbs(\operatorname{Spa}(C, \O_C), \O^+_C/p)$ : $$f^M_! (\cL \otimes\O^+_X/p) \simeq \cF_\cL \otimes \O^{+a}_C/p.$$
    
    Then, for all $0 \le k \le 2d$, there exists an isomorphism : $$H^k_{\e t}(X, \cL^\vee)(d) \simeq \Hom_{\F_p}(\pi_{k - 2d}(\cF_\cL), \F_p).$$ 
\end{corro}


\begin{proof}

    We apply \ref{PCCSPD}, and then take homotopy groups $\pi_{k -2d}$ on both sides. Since $\cL$ is concentrated in degree zero, so is $\cL^\vee$, and we compute : $$ \pi_{2d-k} \left( Rf_*^{H} (\cL^\vee)(d)[2d]\right) \cong R^{k} f_*^{H} \cL^\vee (d)  \cong H^k_{\e t}(X, \cL^\vee)(d)$$
    
    Moreover, $\F_p$ is an injective object of $\cD_{\e t}(\operatorname{Spa}(C, \O_C), \F_p) \cong \cD(\F_p)$ the derived category of $\F_p$-modules, so that we get :  $$\pi_{2d-k}(\operatorname{RHom}_{\cD(\F_p)} (\cF_\cL, \F_p)) = \operatorname{Hom}_{\F_p}(\pi_{k-2d}(\cF_\cL), \F_p)$$ This yields the desired result.
\end{proof}

\begin{rk}
    \label{Functoriality}
    The above computation can be made functorially in $\cL$, provided that $\cF_{\cL}$ can be chosen functorially in $\cL$. In particular, if a group $G$ acts analytically on $X$ or $Y$ such that : 
    \begin{itemize}
        \item $\cL$ is $G$-equivariant.
        \item For all $g \in G$, there is a natural isomorphism morphism $\phi_{g, M}^* \cF_{\cL} \cong \cF_{\phi_{g, H}^* \cL}$, compatible with the action of $G$, where $\phi_g : x \in X \mapsto g \cdot x$.
        \item The isomorphism in the hypothesis is $G$-equivariant.
    \end{itemize}
    Then the Poincaré duality isomorphism is $G$-equivariant.
\end{rk}





\subsection{Standard almost $\O^+/\pi$-étale sheaves and discrete objects in $\Dbs(X, \O_X^+/\pi)$} 

This section follows \cite[Section 3.3]{mann2022-6F}.
\subsubsection{Standard almost $\O^+/\pi$-étale sheaves} Let $X$ be an untilted locally spatial diamond over a complete extension $K$ of $\Q_p$, with pseudo-uniformizer $\pi$. We will define a left-complete category\footnote{This is denoted $\operatorname{Shv}^\wedge(X_{\e t}, \O_X^{+a}/\pi)$ in \cite{mann2022-6F}.} $\cD_{\e t}(X, \O_X^{+ a}/\pi)$ of almost quasi-coherent étale sheaves of $\O^+_X/\pi$-modules over $X_{\e t}$, that compatible with pushforwards and pullbacks. 

Let $f : X \to Y$ be a morphism of such untilted locally spatial diamonds. While the étale pushforward $f_*^S$ of an $\O_X^+/\pi$-module is naturally an $\O_Y^+/\pi$-module, the pullback $f^*_S$ does not preserve $\O^+/\pi$-modules, so that we should consider a different notion of pullback of \textit{quasi-coherent} sheaves. 

Let us first recall the setup in the non-almost case, following \cite[\href{https://stacks.math.columbia.edu/tag/03A4}{Chapter 03A4}]{stacks-project}.

\begin{defi} 
    \label{qcoh_pullback}
    Let $X$ be an untilted locally spatial diamond over a complete extension $K$ of $\Q_p$, with pseudouniformizer $\pi$. 
    
    We may view the sheaf $\O^+_{X}/\pi$ as a sheaf of rings on $X_{\e t}$, and consider the ringed topos $(\operatorname{Sh}(X_{\e t}), \O^+_X/\pi)$ associated to the ringed site $(X_{\e t}, \O_X^+/\pi)$. 
    
    Every morphism of locally spatial diamonds $X \to Y$ induces a morphism of ringed sites $(X, \O^+_X/p) \to (Y, \O^+_Y/p)$. This defines a pullback functor :\footnote{Note that, what we denote $f^*$ here is usually denoted $f^{-1}$ in a quasi-coherent context, while what we denote $f^{*, qcoh}$ is usually just referred to as $f^*$.} $$f^{*, qcoh}_S \cF = \O^+_X/\pi \otimes_{f^*_S \O^+_Y/\pi} f^*_S \, \cG.$$
    
    We let $\operatorname{Mod}(\O^+_X/\pi)$ denote the abelian (1)-category of sheaves of $\O^+_X/\pi$-modules in $X_{\e t}$, and $\cD(X_{\e t}, \O^+_X/\pi)$ be its derived $\infty$-category. Any morphism $f : X \to Y$ induces an adjoint pair : $$ R f_*^{S} : \cD_{\e t}(X, \O^+_X/\pi) \leftrightarrows \cD_{\e t}(Y, \O^+_Y/\pi) : R f_{S}^{*, qcoh} .$$
\end{defi}

Defining the almost version $\cD_{\e t}(X, \O^{+ a}_X/\pi)$, it is a bit more tricky, as the object "$\O^{+a}_X/\pi$" does not exist as a sheaf of classical rings in $\operatorname{Sh}(X_{\e t}, \F_p)$. Instead, we may define a classical sheaf $\cF$ of $\O^{+ a}_X/\pi$-modules on $X$ as a sheaf on $X_{\e t}$, valued in the category of almost $\O^+_X/\pi(X)$-modules, together with the compatible structure of an $\O^+_X/\pi(U)$-module on $\cF(U)$ for all $U \in X_{\e t}$. 

We then let $\cD(X_{\e t}, \O^{+ a}_X/\pi)$ be the derived $\infty$-category of the above category, and $\cD_{\e t}(X, \O_X^{+a}/\pi)$ be its left completion, as in \cite[Defi 3.3.12]{mann2022-6F}. See also \cite[Paragraph 5.5.1]{Gabber_AlmostRing} for the classical theory of almost étale sheaves over schemes. 

Here are the main properties that will be useful to us (cf. \cite[Lemma 3.3.13]{mann2022-6F}). 

\begin{pro}
    \label{Pullback}
    Let $f : X \to Y$ be a morphism of locally spatial diamonds over a complete nonarchimedean field $K$, with pseudo-uniformizer $\pi$. Then, there exists an adjoint pairs of functors : $$ Rf_*^{S} : \cD_{\e t}(X, \O^{+a}_X/\pi) \leftrightarrows \cD_{\e t}(Y, \O^{+a}_Y/\pi) : Rf^{*,qcoh}_{S}.$$

    Moreover, there is a natural t-exact almostification functor : $$(-)^a : \cD_{\e t}(X, \O^{+}_X/\pi) \to \cD_{\e t}(X, \O^{+ a}_X/\pi)$$ It is computed as $\cF^a(U) = \cF(U)^a$, where the second "$a$" denotes the almostification for derived $\O^+_X/p(X)$-modules. The functor $(-)^a$ is compatible with the formation of pushforward and pullbacks.
\end{pro}

\begin{rk}
    \label{ClassicalTensProd}
    Note that there is a natural morphism $- \otimes_{\F_p}^{\L} \O^+_X/p : \cD_{\e t}(X, \F_p) \to  \cD_{\e t}(X, \O_X^+/\pi)$, obtained by left-completing the derived tensor product. Composing this with the almostification functor $\cD_{\e t}(X, \O_X^+/\pi) \to \cD_{\e t}(X, \O_X^{+a}/\pi)$ defines a functor that it closely related to the Riemann-Hilbert functor, cf. \ref{compat_RH}.
\end{rk}

The same definition as in \ref{Overconv} allows us to consider overconvergent objects in $\cD_{\e t}(X, \O^{+ a}_X/\pi)$, and  we denote the full subcategory of overconvergent étale sheaves of $\O^{+ a}_X/\pi$-modules as $\cD_{\e t}(X, \O^{+ a}_X/\pi)^{oc}$ (cf. \cite[Defi 3.3.14]{mann2022-6F} for more detail).

Note that this is really a notion of overconvergence for \textit{almost} sheaves, in particular, the sheaf $\O^+/p$ itself is not overconvergent when viewed as an étale sheaf of $\F_p$-modules. It is, however, \textit{almost} overconvergent. 

Let us now explain how to identify this category with the subcategory of discrete objects in $\Dbs(X, \O^+_X/\pi)$.


    


\subsubsection{Identification of discrete objects}


Recall that there is a fully faithful inclusion of rings into condensed rings (viewed as sheaves of rings on the category of profinite spaces), given by constant sheaves, which factors through solid objects. The essential image of this functor is the full subcategory of \textit{discrete} rings.

For any strictly totally disconnected affinoid perfectoid space $\operatorname{Spa}(A, A^+)$, we may consider the subcategory $\Dbs(A^+/\pi)_\omega \subset \Dbs(A^+/\pi)$ of discrete (derived) $A^+/\pi$-modules. This construction descends to a subfunctor of discrete objects $\Dbs(\O^+_X/\pi)_\omega \subset \Dbs(\O^+_X/\pi)$ for any untilted v-stack $X$, defined as follows.

\begin{defi}(cf. \cite[Def 3.2.17]{mann2022-6F}) \\
    \label{Disc_def}
    Let $X$ be a locally spatial diamond, with pseudo-uniformizer $\pi$. 
    
    We let $\Dbs(X, \O^+_X/\pi)_\omega \subset \Dbs(X, \O^+_X/\pi)$ be the full subcategory formed by objects $\cF$ such that, for any map $f : Y \to X$ from a totally disconnected affinoid perfectoid space, the pullback $f^*_{M} \cF$ is discrete.
\end{defi}

This subcategory satisfies the following : 

\begin{pro}
    \label{Discretization_def}
    For any untilted small v-stack $X$, the subcategory $\Dbs(\O^+_X/\pi)_\omega$ is stable under colimits, and the inclusion $\Dbs(\O^+_X/\pi)_\omega \subset \Dbs(\O^+_X/\pi)$ admits a right adjoint, called the discretization and denoted $(-)_\omega$.
\end{pro}


Discrete objects identify with overconvergent standard almost étale sheaves, in the following sense.

\begin{pro}
    \label{Comparison_Discobj}
    Let $X$ be an untilted locally spatial diamond, with pseudo-uniformizer $\pi$. Assume moreover that $X$ admits a map to an affinoid perfectoid space.
    
    Then, there is an equivalence of $\infty$-categories :
    $$J_X : \Dbs(X, \O^+_X/\pi)_\omega \cong \cD_{\e t}(X, \O^{+a}_X/\pi)^{oc}$$

    Moreover, for all morphism $f : X \to Y$ of untilted locally spatial diamonds over an affinoid perfectoid space, there is a canonical isomorphism of functors $\Dbs(Y, \O^+_Y/\pi)_\omega \to \cD_{\e t}(X, \O^{+,a}_X/\pi)^{oc}$ : $$J_X \circ f^{*}_{M} \cong Rf^{qcoh, *}_{S} \circ J_Y.$$
\end{pro}

\begin{proof}
    This is \cite[Prop 3.3.16]{mann2022-6F}. Note that the hypothesis that $X$ admits a map to an affinoid perfectoid space allows implies that $X$ admits "enough pseudouniformizers", in the sense of \cite[Defi 3.2.6]{mann2022-6F}.
\end{proof}

If $X = \operatorname{Spa}(C, \O_C)$, we denote for simplicity $J_C := J_{\operatorname{Spa}(C, \O_C)}$.

At this point, we constructed tools for two natural functors $\cD_{\e t}(X, \F_p)^{oc} \to \Dbs(X, \O^+_X/\pi)$ ; the Riemann-Hilbert functor, and the almostification of the tensor product of classical étale sheaves. We show that they actually coincide (this was informally remarked in \cite[Paragraph 6.2]{Zavyalov_Abstract_PD}).

\begin{pro}
    \label{compat_RH}
    Let $X$ be an untilted locally spatial diamond, with pseudo-uniformizer $\pi$. Suppose that $X$ admits a map to an affinoid perfectoid space. 
    
    Let $\cF \in \cD_{\e t}(X, \F_p)^{oc}$. Then, the following square commutes :
    \begin{center}
        \begin{tikzcd}[column sep=large]
            & \cD_{\e t}(X, \F_p)^{oc} \ar[r, "- \otimes \O^{+a}_X/\pi"] \ar[d, "\left(- \otimes^\L_{\F_p} \O^+_X/\pi\right)^{a}", swap] & \Dbs(X, \O^+_X/\pi)^{\varphi}_{\omega} \ar[d] \\ 
            & \cD_{\e t}(X, \O^{+a}_X/\pi)^{oc} \ar[r, "J_X^{-1}"] & \Dbs(X,\O^+_X/\pi)_{\omega}
        \end{tikzcd}
    \end{center}
    
    Here, the top horizontal arrow is the Riemann-Hilbert functor of \ref{RH_def}, the rightmost vertical arrow is the restriction to discrete objects of the forgetful functor described in \ref{Forget_phiinv}, the leftmost vertical arrow is the functor described in \ref{ClassicalTensProd} (it preserves overconvergent objects, since taking stalks commutes with tensor products), and the lower horizontal arrow is the identification of \ref{Comparison_Discobj}.
\end{pro}

\begin{proof}
    It follows from \cite[Lemma 3.9.18.(i)]{mann2022-6F}, the proof of \cite[Prop 3.3.16]{mann2022-6F},  \cite[Remark 3.9.11]{mann2022-6F} and \cite[Lemma 3.2.19.(i)]{mann2022-6F} that the four categories satisfy pro-étale descent, so that it suffices to check the commutation on a basis of strictly totally disconnected affinoid perfectoid spaces.
    
    For such spaces, using the identifications from \cite[Lemma 3.3.15.(ii)]{mann2022-6F} \cite[Lemma 3.9.18.(iii)]{mann2022-6F}, this is simply the definition of the Riemann-Hilbert functor, recalled in \ref{local_RH}.
\end{proof}

Moreover, since the equivalence $J_C$ is compatible with pullback, it is also stable by the right adjoints of such maps. This give the following compatibility with pushforwards :

\begin{corro}
    \label{Commute_pushforward_J}
    Let $C$ be a complete algebraically closed extension of $\Q_p$, and let $f : X \to \operatorname{Spa}(C, \O_C)^\diamond$ be a locally spatial diamond. Let $\cF \in \Dbs(X, \O^+_X/p)_\omega$. 
    
    Then, there is a natural isomorphism : $J_C \left( \displaystyle(f_{*}^{M} \cF)_\omega \right) \cong Rf_{*}^{S} (J_X(\cF))$.
\end{corro}

\begin{proof}
    By \ref{Comparison_Discobj}, the two functors $f_{M}^* \circ J_C^{-1}$ and $J_X^{-1} \circ Rf^{qcoh, *}_{S}$ define two naturally isomorphic functors $\cD_{\e t}(\operatorname{Spa}(C, \O_C/p), \F_p) \to \Dbs(X, \O^+_X/\pi)_\omega$ . We'll deduce the result above by passing to right adjoints. Since $J_X$ and $J_C$ are equivalences of categories, they are naturally adjoint to their inverse.

    Recall that $f^{*}_{M}$ preserves discrete objects by definition. By composition of right adjoints, the pullback $f^*_{M} : \Dbs(\operatorname{Spa}(C, \O_C), \O^+_C/\pi)_\omega \to \Dbs(X, \O^+_X/\pi)_\omega$ restricted to discrete objects admits the right adjoint $(f^{M}_{*} -)_\omega$, where $(-)_\omega$ is the discretization functor described in \ref{Discretization_def}.
    
    All étale sheaves over $C$ are overconvergent, and $Rf_{S}^{qcoh, *} : \D_{\e t}(\operatorname{Spa}(C, \O_C), \O^{+ a}_C/\pi) \to \D_{\e t}(X, \O^{+ a}_X/\pi)^{oc}$ admits the right adjoint $Rf^S_*$.
    
    Since naturally isomorphic functors admit naturally isomorphic adjoints, there is a natural isomorphism $J_C \, \circ \, ((f_*^{M} -)_\omega) \cong Rf_*^{S} \circ J_X$, which concludes. 
\end{proof}

In particular, we have the following : 

\begin{corro}
    \label{Disc_OC_over_field}
    Let $C$ be a complete algebraically closed extension of $\Q_p$, and let $f : X \to \operatorname{Spa}(C, \O_C)^\diamond$ be a locally spatial diamond. Let $\cL$ be an overconvergent étale sheaf of $\F_p$-vector spaces over $X$. 
    
    Then, there is a natural isomorphism : $J_C \left( \displaystyle(f_{*}^{M} (\cL \otimes \O^{+a}_X/\pi))_\omega \right) \cong \left(Rf_{*}^{S} (\cL \otimes_{\F_p} \O^+_X/\pi) \right)^a$.
\end{corro}

\begin{proof}
    This follows directly from \ref{compat_RH} and \ref{Commute_pushforward_J} applied to $\cF = \cL \otimes \O^{+a}_X/\pi$, together with the commutation of $(-)^a$ with pushforwards.
\end{proof}

\begin{rk}
    In the setup of the corollary above, assume furthermore that $f$ is quasi-compact, so that $f^{M}_*$ preserves discrete objects by \cite[Lemma 3.3.10.(ii)]{mann2022-6F}.
    
    Letting $\cL = \F_p$, we get $f_*^{M} \O^{+,a}_X/p \cong J_C^{-1} \left(\left( Rf_*^{S} \O^{+}_X/p\right)^a\right)$, which recovers a variant of \cite[Example 3.3.17]{mann2022-6F}. Informally, this means that $f_*^{M} \O^{+,a}_X/p$ computes, up to an almost isomorphism, $R\Gamma_{\e t}(X, \O^+_X/p)$.
\end{rk}

\section{About the topology of locally spatial diamonds}
\label{Sec3}

This section contains many technical results that will be useful in section \ref{Sec4}. 

\subsection{The closed complement of an open sub variety}

A standard problem in p-adic geometry is the fact that, unlike for schemes, the complement of an open subset of a rigid analytic variety needs not admit a structure of an rigid analytic variety, nor of an analytic adic space. Indeed, every morphism of analytic adic space is generalizing (in a topological sense, cf. \cite[\href{https://stacks.math.columbia.edu/tag/0060}{Section 0060}]{stacks-project}), so that the complement of a non-specializing open subset never admits an analytic structure.



However, this is not the only problem. Consider the natural inclusion of the open unit disk in the closed unit disk over $\C_p$, viewed as adic spaces. The topological complement contains a single point of rank 2, and hence does not admit an analytic structure.

In the realm of rigid-analytic varieties, one often uses the following : 

\begin{lemma}
    \label{Open_Open_dec}
    Let $U$ be an quasi-compact admissible open subset of a quasi-separated rigid analytic variety $X$. Then, the complement $X \setminus U$ is an admissible open subset of $X$.
\end{lemma}

\begin{proof}
    This is \cite[Prop 5.6.5]{Huber_Ecoh} or \cite[Lemma 1.5]{VANDERPUT1992219}.
\end{proof}

Recall that the topology on rigid-analytic spaces is a Grothendieck topology, so we may decompose a connected rigid-analytic variety as the disjoint union of two admissible open subsets, as long as the covering is not admissible.

Fix $K$ a complete extension of $\Q_p$, and let $r$ be the functor from the category of rigid-analytic varieties over $K$ to the category of analytic adic spaces over $\operatorname{Spa}(K, \O_K)$ constructed in \cite[Paragraph 1.11]{Huber_Ecoh}. 

\begin{lemma}
    \label{r(U)}
    Let $X$ be a quasi-separated rigid-analytic variety over a complete non-archimedean field $K$, and $U \subset X$ a quasi-compact admissible open subset.
    
    Then, $r(X \setminus U) = (r(X) \setminus r(U))^\circ$, where $T^\circ$ denotes the interior of a topological space $T$.
\end{lemma}

\begin{proof}
    Recall that, since $X$ is quasi-separated, by \cite[Paragraph 1.1.11]{Huber_Ecoh}, $r$ induces an increasing bijection between quasi-compact admissible open subsets of $X$ and quasi-compact open subsets of $r(X)$. Moreover, for any open $U$ of $X$, $r(U) \cap X = U$, so that $r$ is a isomorphism of posets.

    Since $r$ is order-preserving, $r(X) = r((X \setminus U) \sqcup U) \subset r(X \setminus U) \cup r(U)$, so that $r(U \setminus X) \subset r(X) \setminus r(U)$. By \ref{Open_Open_dec}, the complement $U \setminus X$ is an admissible open of $X$, so that $r(X \setminus U)$ is open in $r(X)$, and we have $r(X \setminus U) \subset (r(X) \setminus r(U))^\circ$.

    Since $r(X)$ admits a basis of quasi-compact open subsets, we may write $(r(X) \setminus r(U))^\circ$ as the union of all quasicompact open subsets contained in $r(X) \setminus r(U)$. Let $V$ be such a quasicompact open. By the properties of $r$, we may write $V = r(W)$ for some quasicompact admissible open $W$ of $X$, so that $r(W) \cap r(U) = \emptyset$. Since $r$ is order-preserving, we have $r(W \cap U) \subset r(W) \cap r(U) = \emptyset$, so that $W \cap U = \emptyset$ since $r$ is injective, and $W \subset X \setminus U$, so that $V = r(W) \subset r(X \setminus U)$. 
    
    Hence $(r(X) \setminus r(U))^\circ \subset r(X \setminus U)$, which concludes.
\end{proof}

Huber, in \cite{Huber_Ecoh}, defines a theory of \textit{pseudo}-adic spaces, so that any "nice" (convex and pro-constructible) subset of an analytic adic space admits a pseudo-analytic structure, and, in turn, an étale site. In this paper, we'll bypass the use of pseudo-adic spaces, using locally spatial diamonds instead.

Note that by \cite[Prop 11.15]{Scholze_ECohD}, open subsets of locally spatial diamonds naturally admit a structure of a diamond, which is then a locally spatial diamond by \cite[Prop 11.20]{Scholze_ECohD}.

We use some notations from \cite{Scholze_ECohD}. In particular, for $X$ a diamond, we denote $|X|$ its underlying topological space, and, for $T$ a topological space, we let $\underline{T}$ denote the v-sheaf given by $\underline{T}(X) := \operatorname{Cont}(|X|, T)$. We will use the following lemma :

\begin{lemma}
    \label{ComplementLSD}
    Let $X$ be a locally spatial diamond, and $F \subset X$ be a generalizing closed subset of $X$. Then $F$ is the underlying space of a sub-locally spatial diamond, still denoted $F$, given by the fiber product of v-sheaves $\displaystyle F = \underline{F} \times_{\underline{|X|}} X$.
\end{lemma}

\begin{proof}
    By \cite[Lemma 2.7]{anschütz2022padictheorylocalmodels}, the above defines a closed sub v-stack of underlying topological space $F$. It is also a v-sheaf, since the fiber product can be taken in the category of v-sheaves.

    By \cite[Prop 11.10]{Scholze_ECohD}, it is naturally a diamond, and, by \cite[Prop 11.20]{Scholze_ECohD}, it is locally spatial.
\end{proof}

Hence, the complement of a specializing open subset of an analytic adic space needs not admit a structure of an analytic adic space always admit the structure of a locally spatial diamond. In that setting, we furthermore have a localization (or excision) sequence for étale sheaves, as developed in the following.

\subsection{Localization for locally spatial diamonds}


For a partially proper open immersion $j$ of locally spatial diamonds, we consider the proper pushforward $Rj_!^{S}$ as taken in the sense of \cite[Defi 19.1]{Scholze_ECohD}, as the left adjoint to $Rj^*_{S}$. It coincides with the functor defined by Huber in \cite{Huber_Ecoh}, when $j$ is defined between locally noetherian analytic adic spaces by \ref{Huber_Scholze_OI}.

\begin{pro}
    \label{excision_diamond}
    Let $X$ be a locally spatial diamond whose underlying topological space decomposes as the union $|X| = |U| \sqcup |F|$ of an open and a closed subset respectively, such that $|F|$ is generalizing. Let $j : U \to X$ and $i : F \to X$ be the associated open (resp. closed) immersion of diamonds. 
    
    For any sheaf $\cF$ of abelian groups on $X_{\e t}$, there is an exact sequence in $\cD_{\e t}(X, \Z)$, functorial in $\cF$ :
    $$R j_!^{S} j_{S}^* \cF \to \cF \to R i_*^{S} i_{S}^* \cF$$
\end{pro}

\begin{proof}
    To simplify the notations, we consider implicitly every functor to be in the sense of Scholze.

    
    The morphisms come from (co)unit of the adjunctions. The étale site of $X$ has enough points, cf. \cite[Prop 14.3]{Scholze_ECohD}, so that it suffices to check the exactness on stalks. Let $\overline{x} : \operatorname{Spa}(C, C^+) \to X$ be a quasi-pro-étale geometric point of $X$, where $C$ is algebraically closed, $C^+ \subset C$ is an open and bounded valuation subring, and $\overline{x}$ maps  the unique closed point of $\operatorname{Spa}(C, C^+)$ to some $x \in X$. 

    We compute the stalks using the following lemma, which immediately concludes :
    \begin{lemma}
        \label{stalkpf}
        Under the hypothesis of \ref{excision_diamond}, we have :
    $$(Rj_! \circ j^* \cF)_{\overline{x}} \cong 
    \begin{cases}
        \cF_{\,\overline{x}} & \text{ if } x \in |U| \\
        0 &\text{ otherwise}
     \end{cases} \hspace{.2cm} \text{ and }\hspace{.2cm} (R i_* \circ i^* \cF)_{\overline{x}} \cong 
     \begin{cases}
        \cF_{\,\overline{x}} &\text{ if } x \in |F| \\
        0 &\text{ otherwise}
     \end{cases}$$
    \end{lemma}

    We apply base change (\cite[Prop 19.1]{Scholze_ECohD}, resp. \cite[Corro 16.10.(ii)] {Scholze_ECohD}) to the diagrams :
    
    \begin{center}
        \begin{tikzcd}
            & U \times_X \operatorname{Spa}(C, C^+) \ar[d, "\overline{x}'"] \ar[r, "j'"] & \operatorname{Spa}(C, C^+) \ar[d, "\overline{x}"] \\ 
            & U \ar[r, "j"] & X
        \end{tikzcd}
        and 
        \begin{tikzcd}
            & F \times_X \operatorname{Spa}(C, C^+) \ar[d, "\overline{x}'"] \ar[r, "i'"] & \operatorname{Spa}(C, C^+) \ar[d, "\overline{x}"] \\ 
            & F \ar[r, "i"] & X
        \end{tikzcd}
    \end{center}

    It follows that : $$(R j_! j^* \cF)_{\overline{x}} = R\Gamma(R j_!' \circ \overline{x}'^{\, *} j^*\cF) = R\Gamma(Rj_!' j'^* \, \overline{x}^{\, *} \, \cF) \text{ and }(R i_* i^* \cF)_{\overline{x}} = R\Gamma(R i'_* \, \overline{x}'^{\, *} \,  i^*  \, \cF) = R\Gamma(Ri'_* i'^{\, *} \, \overline{x}^{\, *}\,  \cF),$$ where $R\Gamma$ denotes the derived global section functor on $\operatorname{Spa}(C, C^+)$.
    
    From \cite[Prop 12.10]{Scholze_ECohD}, we know that there are surjections : $$|U \times_{X} \operatorname{Spa}(C, C^+) |  \twoheadrightarrow |U| \times_{|X|} |\operatorname{Spa}(C, C^+) |\text{ and } |F \times_{X} \operatorname{Spa}(C, C^+)|  \twoheadrightarrow |F| \times_{|X|} |\operatorname{Spa}(C, C^+)|.$$ Here, the fiber product on the right hand side is taken as topological spaces. Topologically $|\operatorname{Spa}(C, C^+)|$ is a totally ordered chain of specializations. It admits a unique closed point, which is a specialization of every other points. since continuous maps preserves specializations, the image of $|\operatorname{Spa}(C, C^+)|$ in $|X|$ is contained in the set of generalizations of $x$. 
    
    \vspace{.1cm}
    Assume that $x \not \in |U|$. Since $|F| = |X| \setminus |U|$ is a generalizing subset of $|X|$ by hypothesis, the image of $\overline{x}$ is contained in $|F|$, and hence disjoint from $|U|$, so that the fiber product $|U| \times_{|X|} |\operatorname{Spa}(C, C^+)|$ is empty. Hence, $U \times_{X} \operatorname{Spa}(C, C^+)$ maps surjectively to the empty set, hence is empty, so that $x'^{\, *}$ is zero, and hence $(R j_! j^* \cF)_{\overline{x}} = 0$, as desired.

    \vspace{.1cm}
    The exact same argument holds for $x \not \in |F|$, since $|U|$ is generalizing (as it is open), so that, for $x \not \in F$, $(Ri_*i^* \cF)_{\overline{x}} = \emptyset$.

    \vspace{.1cm}
    Now, assume that $x \in |U|$. We'll show that $R j'_! j'^{\, *} \cF \cong \cF$. The topological space $|U \times_X \operatorname{Spa}(C, C^+)|$ identifies with an open subset of $|\operatorname{Spa}(C, C^+)|$ containing the unique closed point $x \in \operatorname{Spa}(C, C^+)$, so that, $|U \times_X \operatorname{Spa}(C, C^+)| = |\operatorname{Spa}(C, C^+)|$. Since $j'$ is an open immersion, it's an injection of perfectoid spaces, so that the combined \cite[Prop 5.3]{Scholze_ECohD} and \cite[Prop 5.4]{Scholze_ECohD} prove that $j$ is an isomorphism, hence $R j'_! j'^{\, *} \cong id$, which gives the desired stalk.

    \vspace{.1cm}
    Finally, if $x \in |F|$, which is generalizing, the image of $|\operatorname{Spa}(C, C^+)|$ is included in $|F|$. Since, by \ref{ComplementLSD}, $F = \underline{F} \times_{|\underline{X}|} X$, the morphism $x :\operatorname{Spa}(C, C^+) \to X$ then factors through $F$, so that there exists a diagonal map $x_F$ making the diagram commute :

    \begin{center}
        \begin{tikzcd}
            & F \times_X \operatorname{Spa}(C, C^+) \ar[d, "\overline{x}'"] \ar[r, "i'"] & \operatorname{Spa}(C, C^+) \ar[d, "\overline{x}"] \ar[dl, dashed, swap, "x_F"] \\ 
            & F \ar[r, "i"] & X
        \end{tikzcd}
    \end{center}

    By the universal property of fiber products, this diagonal map $x_F$ induces a splitting of $i'$, which is hence surjective. Finally, a surjective closed immersion of perfectoid spaces is an isomorphism (cf. \cite[Lemma 5.3 and 5.4]{Scholze_ECohD}), so that $i'$ induces an isomorphism of diamonds. 
    
    Hence $Ri'_* i'^{\, *} \cong id$, which concludes.
\end{proof}


\begin{rk}
    In the above setup, we may rewrite the \textit{fake} base change results \cite[Thm 19.2 and Rk 19.3]{Scholze_ECohD} as a real base change result, without needing to go through pseudo-adic spaces.
\end{rk}

Note that this also holds when we consider quasicoherent pullbacks, as defined in \ref{qcoh_pullback}.

\begin{corro}
    \label{excision_diamond_qcoh}
    Let $X$ be an untilted locally spatial diamond with pseudo-uniformizer $\pi$, whose underlying topological space decomposes as the union $|X| = U \sqcup F$ of an open and a closed subset respectively, such that $F$ is generalizing. Let $j : U \to X$ and $i : F \to X$ be the associated open (resp. closed) immersion of diamonds. 
    
    For any sheaf $\cF$ of $\O^+_X/\pi$-modules on $X_{\e t}$, there is an exact sequence in $\cD_{\e t}(X, \O^+_X/\pi)$ :
    $$R j_!^{S} j_{S}^{*, qcoh} \cF \to \cF \to R i_*^{S} i_{S}^{*, qcoh} \cF$$

    This construction is moreover functorial in $\cF$.
\end{corro}

This directly follows from the following lemma :

\begin{lemma}
    \label{Qcoh_Pullback_qproet}
    Let $f : X \to Y$ be a quasi-pro-étale morphism of locally spatial diamonds. 
    
    Then there is an natural isomorphism of functors $f^{qcoh, *}_{S} \cong f^{*}_S$.
    
    This is in particular the case if $f$ is an open immersion or a generalizing closed immersion.
\end{lemma}

\begin{proof}
    Recall that, for any sheaf $\cF$ of $\O^+_Y/\pi$-modules, the quasicoherent pullback $f^{qcoh, *}_S$ is defined as the tensor product $f^*_S \cF \otimes_{f^{*}_S (\O^+_Y/\pi)} \O^+_X/\pi$. But, for a quasi-pro-étale $f$, it follows from \ref{Pullback_qproet}, $f^*_S \, \O^+_Y/\pi \cong \O^+_X/\pi$, which yields the desired result.
    
    Finally, open immersions are quasi-pro-étale by definition, and generalizing closed immersions are quasi-pro-étale by \cite[Corro 10.6]{Scholze_ECohD}.
\end{proof}

\subsection{A convenient open covering}


In this section, we construct two intertwined covering of any partially proper analytic adic space : one consisting of qcqs open, and the other one consisting of partially proper open subsets. Such a covering will play an important technical role in the following section. 

Before anything else, let us first establish two topological lemmas.

\begin{lemma}
    \label{NiceCov}
    Let $X$ be a topological space that admits a basis of quasi-compact quasi-separated (qcqs for short) open subsets and such that, for any quasi-compact open $U \subset X$, the closure $\overline{U}$ is quasi-compact.
    
    Then, there exists a filtered poset $(I, \le)$ and a decomposition of $X$ as an increasing union of qcqs open subspaces $(X_i)_{i \in I}$, such that : 
    \begin{enumerate}[itemsep=0pt]
        \item For any $i \in I$, there are only finitely many $j \in I$ such that $j \le i$.
        \item For any $i < i'$, $\overline{X_i} \subset X_{i'}$.
    \end{enumerate}
\end{lemma}

\begin{proof}
    Since $X$ admits a basis of qcqs open subsets, we may chose a covering $X = \bigcup_{\lambda \in \Lambda} X_\lambda$ by some qcqs open $X_\lambda$. Let us denote $\mathcal{P}_\omega(\Lambda)$ the set of finite subsets of $\Lambda$.
    
    For any $J \in \mathcal{P}_\omega(\Lambda)$, let $X_J = \bigcup_{\lambda \in J} X_\lambda$, so that the $(X_J)_{J \in \mathcal{P}_\omega(\Lambda)}$ form a filtered covering of $X$ by qcqs open subsets. For any $J \in \cP_\omega(\Lambda)$, consider the open cover $\overline{X_J} \subset \bigcup_{J \subset J'} X_{J'}$. Since $\overline{X_J}$ is quasi-compact, we can extract a finite subcover, so that there exists some $m(J) \in \mathcal{P}_\omega(\Lambda)$ such that $J \subset m(J)$ and $\overline{X_{J}} \subset X_{m(J)}$.
    
    By construction, for any $J \in \mathcal{P}_\omega(\Lambda)$, $J \subset m(J)$ and $\overline{X_J} \subset X_{m(J)}$.
    
    Now, define a partial order on $\mathcal{P}_\omega(\Lambda)$, denoted $\le'$, by : $$A \le' B \Leftrightarrow A = B \text{ or } m(A) \subset B.$$ By construction, it satisfies both desired properties. 
    
    It is moreover filtering since, for any $A, B \in \mathcal{P}_\omega(\Lambda)$, both $A \le' m(A) \cup m(B)$ and $B \le' m(A) \cup m(B)$.
\end{proof}

\begin{lemma}  
    \label{incr_qcqs}
    Let $X$ be a topological space that admits a basis of quasi-compact open subsets.
    
    Let $F$ be a quasi-compact closed subset of $X$, and $U$ an open subset such that $F \subset U$. Then there exists a quasicompact open $V$ of $X$ such that $F \subset V \subset U$.
\end{lemma}

\begin{proof}
    Since the quasi-compact opens form a basis of the topology, every $x \in F$ admits a quasi-compact neighborhood contained in $U$, denoted $U_x$. Then, $F \subset \bigcup_{x \in F} U_x$. Since $F$ is quasi-compact, one may bextract a finite sub-covering $F \subset \bigcup_{i=1}^n U_{x_i}$. Then $V = \bigcup_{i=1}^n U_{x_i}$ satisfies the desired property.  
\end{proof}

We may now state our result. Note that it can be seen as a more general variant of the covering used in the proof of \cite[Prop 6.2.5]{ACLB_Ansch_Mann_6FFonFF}. Compare also with the covering constructed in \cite[Corro 2.29]{Achinger_Niziol_cpctsupport}, or with \cite[Lemma 5.3.3]{Huber_Ecoh}.

\begin{pro}
    \label{VNiceCov}
    Let $f : X \to Y$ be a partially proper morphism of locally noetherian analytic adic spaces, with taut $Y$. Then, there exists a filtered poset $(I, \le)$, together with two open coverings of $X$ :
    \begin{itemize}
        \item An increasing open cover $X = \bigcup_{i \in I} X_i$ by qcqs open subsets.
        \item An increasing open cover $X = \bigcup_{i \in I} U_i$ by partially proper open subsets.
    \end{itemize}
    
    such that : 
    
    \begin{enumerate}[topsep=0pt]
        \item For any $i \in I$, there are only finitely many $i' < i$.
        \item For all $i \in I$, $X_i \subset U_i$ and for any $i < j$, $U_i \subset X_j$.
        \item For any $i < j$, $\overline{X_i} \subset U_{j}$, and the morphism $X_i \to X_j$ factors through the universal adic compactification $\overline{X_i}^Y$ of $X_i$ over $Y$, as constructed in \cite[Thm 5.1.5]{Huber_Ecoh}.
    \end{enumerate}  
\end{pro}

\begin{proof}
    Since $f$ is a partially proper, it is taut, so that, since $Y$ is taut, $X$ is also taut. Moreover, $X$ admits a basis of affinoid open subsets, so that the assumptions of \ref{NiceCov} apply. Let $(I, \le)$ and $(X_i)_{i \in I}$ be a filtered open cover of $X$ by qcqs open subsets, satisfying the conditions of \ref{NiceCov}.

    For any $i \in I$, let $Y_i = \overline{X_i}$ be the closure of $X_i$ in $X$. Since $X$ is taut, $Y_i$ is quasi-compact in $X$. Let $Z_i = Y_i \setminus X_i$, which is a quasi-compact closed subset of $Y_i$. The above lemma \ref{incr_qcqs}, applied to $Z_{i} \subset (Y_{i} \setminus \bigcup_{j<i} Y_{j})$ in $Y_i$, yields a quasi-compact open subset $V_i$ of $Y_i$ such that $Z_i \subset V_i \subset Y_{i} \setminus \bigcup_{j < i} Y_{j}$.

    Since $Z_i \subset V_i$, $X_i \cap V_i = \emptyset$, and $X_i$ is open, so that $X_i \cap \bar{V_i} = \emptyset$. Let $U_{i} = Y_{i} \setminus \overline{V_{i}}$, so that, for all $i < j$; $X_i \subset U_i \subset X_{j}$ using the second hypothesis from \ref{NiceCov}. We claim that the map $U_i \to X$ is partially proper. It suffices to show that $U_i$ is specializing, hence that $\overline{V_i}$ is generalizing. $V_i$ is quasi-compact, hence retrocompact, so, by \cite[Ch 0, Corro 2.2.27]{fujiwara2017foundationsrigidgeometryi}, $\overline{V_i}$ is simply the closure under specializations of $V_i$. Since $X_i$ is an analytic adic space, all specializations are vertical, and, by \cite[Lemma 1.1.10.(i)]{Huber_Ecoh}, the set of generalizations (or specializations) of any given point forms a chain, so that $\overline{V_i}$ is stable under generalization. Hence the $U_i$ are specializing open subsets of $X$, hence partially proper.
    
    It suffices to prove the last property. For any $i < j$, there is a natural inclusion $X_i \subset U_i \subset X_j$. Since $X_i \to U_i \to Y$ is the composition of an open immersion and a partially proper morphism, it defines a compactification of $X_i$ over $U_i$, so that, by \cite[Thm 5.1.1]{Huber_Ecoh}, the morphism $X_i \to U_i$ factors through the universal compactification $\overline{X_i}^Y$. Hence, it is also the case for the morphism $X_i \to X_j$.
\end{proof}

In \cite[Section 6.2]{ACLB_Ansch_Mann_6FFonFF}, it is crucial for some arguments that a covering as constructed above can be made to be \textit{countable}. We show that this is equivalent to a the more standard paracompactness assumption. 

Let's say that a rigid-analytic variety is \textit{paracompact} if it admits a locally finite covering by admissible affinoids, cf. \cite[Definition 2.5.6]{deJong1996}.

\begin{pro}
    Let $X$ be a connected partially proper rigid-analytic variety over $\operatorname{Spa}(C, \O_C)$\footnote{The argument would work mutatis mutandis for locally noetherian analytic adic spaces over a more general base, but we stick to the generality used in the litterature.}, where $C$ is a complete algebraically closed extension of $\Q_p$. Then, the following are equivalent : 
    \begin{enumerate}[topsep = 0pt]
        \item $X$ is paracompact.
        \item $X$ is countable at infinity in the sense of \cite[Defi 6.2.5]{ACLB_Ansch_Mann_6FFonFF}, i.e. we may write may write $X = \bigcup_{n \in \N} X_n$ as a countable increasing union of qcqs open subspaces $X_n$ such that, for each $n$, the inclusion $X_n \subset X_{n+1}$ factors through the adic compactification $Y_n$ of $X_n$ over $C$. 
    \end{enumerate}
\end{pro}

Note that the "connected" hypothesis is purely technical - otherwise it is possible that $X$ admits uncountably many connected components, and one may not expect to find a countable covering.

\begin{proof}
    $(i) \implies (ii)$ Since $X$ is connected, we may choose a \textit{countable} admissible locally finite affinoid covering in \ref{NiceCov}, using \cite[Lemma 2.5.7]{deJong1996}. Then, the construction of \ref{VNiceCov} produces such a result.
    
    $(ii) \implies (i)$. Write $X = \bigcup_{n \in \N} X_n$ as a strictly increasing countable union of qcqs open subsets, as in the definition of countability at infinity, and let $Y_n$ be the adic compactification of $X_n$ over $\operatorname{Spa}(C, \O_C)$. Consider $Y_n \setminus X_{n-1}$ as a closed subset of $X_{n+1} \setminus Y_{n-2}$. By the separation lemma \ref{incr_qcqs}, for all $n \in \N$, there exists qcqs open subsets $U_n$ in $X_{n+1}$ such that $Y_{n} \setminus X_{n-1} \subset U_n \subset X_{n+1} \setminus Y_{n-2}$. By construction, each $U_n$ might intersect $U_{n \pm 1}$, but not any $U_k$ for $|k - n| > 1$.
    
    Each $U_n$ is a quasi-compact open subset of $X$, which can hence be written as a finite union of affinoids. Such collection of affinoids is locally finite by construction, so that $X$ is paracompact.
\end{proof}

\begin{rk}
    Note that the above applies to period domains, as they are paracompact by \cite[Lemma A.3]{Hartl_2013}. 
\end{rk}

\section{Compatibility between Mann's and Huber's proper formalisms}
\label{Sec4}

The goal of this section is to identify Mann's $f_!^{M}$ and Huber's $Rf_!^{H}$ (whose construction will be shortly recalled in paragraph \ref{Rf!H}) for partially proper morphisms $f$. The precise statement is as follows : 

\begin{pro}
    \label{Identification_lower_shriek}
    Let $f : X \to \operatorname{Spa}(C, \O_C)$ be a partially proper locally noetherian analytic adic space over a complete algebraically closed extension $C$ of $\Q_p$. Let $\cL$ be an overconvergent étale sheaf of $\F_p$-modules on $X$. 
    
    There is a natural isomorphism in $\cD_{\e t}(\operatorname{Spa}(C, \O_C), \O_C^{+, a}/p)^{oc}$ : $$ J_C\left(f_{!}^{M} (\cL \otimes \O^{+,a}_X/p)\right) \cong  \left(Rf_!^{H}(\cL \otimes_{\F_p} \O^+_X/p)\right)^{a},$$ where $J_C$ denotes the identification functor from \ref{Comparison_Discobj}.
\end{pro}

Here is a sketch of the proof. 

\begin{foreshadow} \label{foreshadow} Fix an open covering of $X$ by partially proper open subsets $(U_i)_{i \in I}$, as constructed in \ref{VNiceCov}. We will proceed in four steps, as follows : 
\begin{enumerate}[itemsep=.1cm]
    \item Prove that, for any standard étale sheaf $\cF$ of $\O^+_X/p$-modules on $X$, there is an exact sequence in $\cD_{\e t}(\operatorname{Spa}(C, \O_C), \O_C^+/p)$ : 
    $$ R f_{!}^{H} \cF \to R f_{*}^{H} \cF \to \rlim_{i \in I} Rf^{H}_{*} \circ i^{X \setminus U_i, S}_{*} \circ i^{X \setminus U_i, *}_S \cF ,$$
    were, $i^{X \setminus U_i}$ denotes the closed immersion of locally spatial diamonds $X \setminus U_i \to X$.
    
    \item Prove that $f^{M}_!$ preserves discrete objects. 

    \item Establish, for any $\cF \in \Dbs(X, \O^+_X/p)_\omega$, the following exact sequence in $\Dbs(\operatorname{Spa}(C, \O_C), \O^+_C/p)_\omega$ :
    $$ f_{!}^{M} \cF \to \left(f_{*}^{M}(\cF)\right)_\omega \to \rlim_{i \in I}  \left( f^{M}_* \circ i^{X \setminus U_i, M}_{*} \circ i^{X \setminus U_i, *}_M \cF \right)_\omega,$$
    
    where the notations are the same as in the first point.
    
    \item Identify the second and third terms of the two exact sequences above, and deduce an isomorphism between the first terms. 
\end{enumerate}
\end{foreshadow}

Let us do so in order. The first point is an essentially standard result (cf. \cite[Corro 2.17]{Achinger_Niziol_cpctsupport}), but is usually written for rigid analytic varieties, where the colimit is taken along a covering by quasi-compact open subsets, and the complement is taken in the rigid-analytic sense, as in \ref{Open_Open_dec}. We will recover this result in \ref{corro_Fib_H}, by viewing the complement as a locally spatial diamond, as in \ref{ComplementLSD}.

\subsection{Proper pushforward for almost quasi-coherent classical étale sheaves}

\label{Rf!H}



First, let us recall Huber's definition of  direct image with compact support $R f^H_{!}$ \cite[Defi 5.2.1]{Huber_Ecoh} for a partially proper morphism $f$. 

For any separated morphism of locally noetherian adic spaces $f : X \to Y$, we may define a left exact functor $f^{H}_! : X_{\e t}^\sim \to Y_{\e t}^\sim$  by : $$ f^{H}_! \cF : U \in Y_{\e t} \mapsto \{ s \in \Gamma(U \times_Y X, \cF), \, \mathrm{supp}(s) \text{ is proper over } U\}.$$

Whenever $f$ is étale, $f^{H}_!$ is exact, and if $f$ is furthermore separated, $f^{H}_!$ is left adjoint to $f_H^*$ by \cite[Lemma 2.7.6.(i) and Prop 5.2.4]{Huber_Ecoh}. When $f$ is proper, $f_!^H = f_*^H$. 

For a general morphism $f$, we should not expect the derived functors of $f^{H}_!$ to have good properties - rather, for morphisms $f$ that can be written as a composition $f = p \circ i$ of a proper morphism $p$ and an open immersion $i$, we should define the derived pushforward $R f^{H}_!$ as $R p^{H}_* \circ i^{H}_!$, where $i_!^H$ implicitely denotes the extension of the exact functor to the derived category.

However, for partially proper morphism $f$ of analytic adic spaces, the derived functors of $f^{H}_!$ admit reasonable properties. In \cite{Huber_Ecoh}, Huber considers a right derived functor $R^+ f^{H}_! : D^+(X_{\e t}, \Z) \to D^+(Y_{\e t}, \Z)$ at the level of the bounded below derived $1$-categories, and additionally constructs a fully derived functor $R f^{H}_! : D(X_{\e t}, \Z) \to D(Y_{\e t}, \Z)$ under some finite dimensionality assumption on $f$. However, using \cite[\href{https://stacks.math.columbia.edu/tag/070K}{Lemma 070K}]{stacks-project} and the discussion at the beginning of \cite[\href{https://stacks.math.columbia.edu/tag/07A5}{Section 07A5}]{stacks-project}, we may define a right derived functor $R f_{!}^H : D(X_{\e t}, \F_p) \to D(Y_{\e t}, \F_p)$ for any partially proper morphism $f$, as the (unbounded) derived functors of $f_!$. Indeed, the categories of modules on a site form Grothendieck abelian categories (cf. \cite[\href{https://stacks.math.columbia.edu/tag/07A5}{Section 07A5}]{stacks-project}).

This construction naturally admits an $\infty$-categorical enhancement (cf. \cite[Section 1.3.5]{Lurie_HA}), i.e. there exists a canonical morphism of $\infty$-categories $Rf_{!}^H : \cD(X_{\e t}, \F_p) \to \cD(Y_{\e t}, \F_p)$, that induces the morphism above at the level of homotopy categories, and it passes to left completions to a morphism $Rf_{!}^H : \cD_{\e t}(X, \F_p) \to \cD_{\e t}(Y, \F_p)$.

Likewise, we can check that, for any partially proper morphism $f : X \to Y$ and any sheaf $\cF$ of $\O_X/p$-modules on $X$, $f^{H}_! \cF$ admits a natural structure of an $\O_Y/p$-module, viewing every $f_!^{H}\cF(U)$ as a submodule of $\Gamma(U \times_Y X)$.

Hence, the derived functor $Rf^{H}_!$ induces a functor $\cD_{\e t}(X, \O_X^+/p) \to \cD_{\e t}(Y, \O_Y^+/p)$. 

\begin{rk}
    \label{Huber_Scholze_OI}
    Note that this approach is not the one of \cite[Defi 22.13]{Scholze_ECohD}. In loc. cit, Scholze defines, for $j$ an open immersion, $j^{S}_!$ as the left adjoint of $j^*_{S}$, and the definition of $f_!$ for a general $f$ is defined by Kan extension from the setup where one may write $f = p \circ j$ as the composition of a proper map and an open immersion (with $p^S_! = p_*^S$).
    
    It is however clear from definitions that, if $j : X \to Y$ is a separated open immersion between locally noetherian adic spaces, $j^*_H \cong j^*S$ so that, by uniqueness of left adjoints, $j_!^H \cong j_!^S$.
\end{rk} 


\begin{lemma}
    \label{pseudo_excision}
    Let $f : X \to Y$ be a partially proper morphism of locally noetherian analytic adic spaces, with taut $Y$. Let $(U_i)_{i \in I}$ be an increasing covering of $X$ by partially proper open subsets as in \ref{VNiceCov}.
    
    Then, for any sheaf $\cF$ of abelian groups over $X_{\e t}$, there is a natural isomorphism in $\cD_{\e t}(X, \Z)$ : $$\rlim_{i \in I} Rf^{H}_* \circ j_{U_i \, !}^{H} \circ j_{U_i, H}^{*} \cF \cong R f_{!}^{H} \cF.$$
\end{lemma}

\begin{rk}
   For technical reasons, the result is stated with a colimit alongside partially proper open, but the colimit is the same when taken along qcqs open subsets, as one can take intertwining cofinal coverings as in \ref{VNiceCov}. 
\end{rk}

    
\begin{proof}
    For any $i \in I$, there is a natural morphism $Rj^{H}_{U_i !} \, j^{*}_{U_i, H} \cF \to \cF$ induced by the counit of the adjunction. Everything lives in the bounded below subcategory $\cD^+(X_{\e t}, \F_p) \subset \cD_{\e t}(X, \F_p)$, so that we may forget about the left completions.
    In order to check that it is an isomorphism in the derived category, it suffices to show that it is a quasi-isomorphism, i.e. to check that it is an isomorphism on homotopy groups. By the proof of \cite[Prop 1.3.5.21]{Lurie_HA}, since we're working over the derived category of a Grothendieck abelian category, the functors $\pi_k$ commute with filtered colimits, so that : $$\pi_{-n} \left( \rlim_{i \in I} R f_*^{H} \circ j_{U_i !}^{H} \circ j_{U_i, H}^{*} \cF \right) \cong \rlim_{i \in I} R^n f_{*}^{H} \circ j_{U_i !}^{H} \circ j_{U_i, H}^{*} \cF.$$ 
    
    We easily check that the covering $(U_i)_{i \in I}$ satisfies the conditions of \cite[Lemma 5.3.3.(ii)]{Huber_Ecoh}, so that the result follows from \cite[Lemma 5.3.3.(iii)]{Huber_Ecoh}.
\end{proof}

    
We may now establish the announced exact sequence.
        
\begin{pro}
    \label{Fiber_H}
    Let $f : X \to Y$ be a partially proper morphism of locally noetherian analytic adic spaces, with taut $Y$. Let $(U_i)_{i \in I}$ be an increasing covering of $X$ by partially proper open subsets as constructed in \ref{VNiceCov}.
    
    For any $i \in I$, we let $Z_i$ be the complement of the partially proper $U_i$, viewed as a closed locally spatial subdiamond of $X^\diamond$ by \ref{ComplementLSD}. Let $i_{Z_i} : Z_i \to X^{\diamond}$ be the associated closed immersion. Then, for any sheaf $\cF$ of abelian groups over $X_{\e t}$, there is a natural exact sequence in $\cD_{\e t}(Y, \Z)$ : $$R f_{!}^{H} \cF \to R f_{*}^{H} \cF \to \rlim_{i \in I} R f_{*}^{H} \circ i_{Z_i \, *}^{S} \circ i_{Z_i, S}^{*} \cF.$$ 
\end{pro}

Note that we implicitly identify the étale site of $X$ with the étale site of $X^\diamond$, so that the composition $R f_{*}^{H} \circ i_{Z_i, *}^S$ makes sense, even though $Z_i$ only admits the structure of a locally spatial diamond.


\begin{proof}
    Using the localization sequence established in Proposition \ref{excision_diamond}, there is a short exact sequence of sheaves of abelian groups on $Y_{\e t}$ : $$ 0 \to j^{S}_{U_i \, !} \, j_{U_i, S}^{*} \cF \to \cF \to i_{Z_i \, *}^{S}\,  i_{Z_i, S}^{*} \cF \to 0.$$ 
    
    This induces an exact sequence in $\cD_{\e t}(X, \Z)$. Applying the exact functor $R f^H_*$ yields the exact sequence : $$ R f^{H}_* \circ j^{S}_{U_i \, !} \circ j_{U_i, S}^{*} \cF \to R f^{H}_* \cF \to R f^{H}_* \circ i^{S}_{Z_i *} \circ i_{Z_i, S}^{*} \cF. $$

    We then form the colimit over $I$. Since the the morphisms from passing from $i$ to some $j > i$ are induced by (co)units of the adjunctions, the diagram is homotopy coherent, and, since colimits of exact sequences remain exact (as they can be viewed as cofiber sequences), we have :
    $$\rlim_{i \in I} R f^{H}_* \circ R j_{U_i \, !}^S \circ j_{U_i, S}^* \cF \to R f^H_* \cF \to \rlim_{i \in I} R f^H_* \circ  i_{Z_i \, *}^S \circ i_{Z_i, S }^* \cF.$$

    By \ref{Huber_Scholze_OI}, the functors $j_{U_i}^*$ and $j_{U_i, !}$ coincide in Huber and Scholze's formalism, so that the first term identifies to $Rf^{H}_! \cF$ by \ref{pseudo_excision}, which concludes.   
\end{proof}

Now, we prove a version for sheaves of $\O^+/p$-modules. 

\begin{corro}
    \label{Exact_Seq_Hub_Qcoh}
    Consider the setup of \ref{Fiber_H} over $K$, and let $\cF$ be an étale sheaf of $\O^+_X/p$-modules on $X$. Then, there is an exact sequence : 
    $$ R f_{!}^{H} \cF \to R f^{H}_* \cF \to \rlim_{i \in I} R f^{H}_* \circ i_{*, H}^{Z_i} \circ i_{Z_i, H}^{*, qcoh} \cF.$$ 
\end{corro}

\begin{proof}
    By \ref{Qcoh_Pullback_qproet}, the quasi-coherent pullback $i_{Z_i, H}^{*, qcoh}$ coincides with the pullback $i_{Z_i, H}^*$ as étale sheaves, so that this directly follows from the above result, since the natural forgetful functor $\cD_{\e t}(X, \O^+_X/p) \to \cD_{\e t}(X, \Z)$ respects exact sequences.
\end{proof}

We conclude this paragraph by noting that we can deduce an alternative writing of \ref{Fiber_H}, that remains in the realm of analytic adic spaces, in the spirit of \cite[Corro 2.17]{Achinger_Niziol_cpctsupport}. Following \ref{r(U)}, we replace the $Z_i$ by the interior of the complement of some qcqs open subspaces. This is what is done in the following :

\begin{corro}
    \label{corro_Fib_H}
    Let $f : X \to Y$ be a partially proper morphism of locally noetherian analytic adic spaces, with taut $Y$. Let $(X_i)_{i \in I}$ be an increasing covering of $X$ by quasicompact open subsets. 

    For any $i \in I$, let $Y_i = (X \setminus X_i)^\circ = X \setminus \overline{X_i}$, which admits an analytic structure as an open subset of $X$. Let $j_{Y_i} : Y_i \to X$ denote the associated open immersion. Then, for any étale sheaf $\cF$ of abelian groups on $X_{\e t}$, there is an exact sequence : $$ R f_{!}^{H} \cF \to R f^{H}_* \cF \to \rlim_{i \in I} R f^{H}_* \circ j^{H}_{Y_i \, *} \circ j_{Y_i, H}^{*} \cF.$$ 
\end{corro}

\begin{proof}
    Up to reducing $I$, we may assume that, for any $i$, there is only finitely many $j$ such that $j < i$. By the proof of \ref{VNiceCov}, we may construct partially proper open subsets $(U_i)_{i \in I}$ such that $X_i \subset U_i \subset X_j$ for all $i < j$, satisfying the properties of \ref{VNiceCov}. 
    
    We reuse the notations from the previous proposition. Hence, it suffices to check that there is a canonical isomorphism $\rlim_{i \in I} R f^H_* \circ  i_{F_i \, *}^H \circ i_{F_i, H}^* \cF \to \rlim_{i \in I} R  f_*^H \circ j_{Z_i \, *}^H \circ j_{Z_i, H}^* \cF$.

    From the properties of the covering, for any $i < j$ in $I$, we have $Y_j \subset Z_j \subset Y_i$, so that both families of subsets are cofinal in one another. Since the category $\cD(X_{\e t}, \Z)$ admits all colimits, the colimit alongside alternating $Z$ and $Y$'s exists, and can can be computed alongside one family or the other, so that both two colimits coincide.
\end{proof}

\subsection{Mann's proper pushforward preserves discrete objects}

Let us proceed with the established program in \ref{foreshadow}. This subsection is dedicated to the proof of the following :

\begin{pro}
    \label{preservesdisc}
    Let $f : X \to Y$ be a partially proper morphism of locally noetherian analytic adic spaces over a complete extension $K$ of $\Q_p$. 
    
    Then, for any $\cF \in \Dbs(X, \O^+_X/p)_\omega$, $f^{M}_! \cF$ is discrete in $\Dbs(Y, \O^+_Y/p)$.
\end{pro}

To avoid busy notations, until the end of this paragraph, all functors are implicitly in the sense of Mann, and we drop the superscript "M". We moreover fix such a field $K$.

We will reduce the the case where $f$ factors as the composition of a partially proper open immersion and a proper morphism. Note that, the pushforward by a proper morphism preserves discrete objects by \cite[Lemma 3.3.10.(ii)]{mann2022-6F}, so that it will then suffice to deal with open immersions. We need with the following lemma :

\begin{lemma} \label{colim_union}
    Let $(I, \le)$ be a filtered partial order, and $(U_i)_{i \in I}$ be an increasing open covering of a locally spatial diamond $X$ over $K$. For any $i \in I$, let $g_i : U_i \to X$ denote the open immersion. 
        
    Then, there is a natural isomorphism $\cF \cong \rlim_{i \in I} g_{i \, !} \, g_i^{*} \, \cF$  for any $\cF \in \cD^a_{\square}(X, \O^+_X/p)$.
\end{lemma}

\begin{proof}
    For any $j \in I$, $g_j^{*} = g_j^{!}$ admits the left adjoint $g_{j  \, !}$. The counit of the adjunctions defines a morphism $\rlim_{i \in I} g_{i \, !} \, g_i^{*} \cF \to \cF$.
        
    By (hyper)v-descent \cite[Thm 1.2.1]{mann2022-6F}, there is an isomorphism : $$\Dbs(X, \O^+/p) \cong \llim_{n \in \Delta} \Dbs(Y_n, \O^+/p)\text{, where } Y = \bigsqcup_{i \in I} U_i \text{ and } Y_n = Y^{\times n} = \bigsqcup_{i_1, \dots, i_n \in I}\displaystyle U_{i_1} \times_X \dots \times_X U_{i_n}.$$
        
    Hence, it suffices to check the isomorphism after pulling back to all such $U_{i_1} \times_X \dots \times_X U_{i_n}$, for fixed $n \in \N$ and $(i_1, \dots, i_n) \in I$. Since all the $U_i$ are open subsets of $X$, the fiber product is simply the intersection, and, since $I$ is filtered, every such morphism $U_{i_1} \times_X \dots \times_X U_{i_n} \to X$ factors through $U_{j}$ for some large enough $j \in I$. Hence, it suffices to prove that, for all $j \in J$, the map : $$g_j^* \left(\rlim_{i \in I} g_{i \, !} g_i^* \cF\right) \to g_j^* \cF $$ is an isomorphism. Since $g_j^*$ admits a right adjoint, it commutes with colimits, so that it suffices to prove that : $$\rlim_{i \in I} g_j^* \, g_{i \, !} \, g_i^* \, \cF \cong g_j^* \, \cF$$
        
    Since the covering is increasing, the fiber product $U_i \times_X U_j$ is $U_{min(i,j)}$. Since $I$ is filtered, the subset  $J = \{ i \in I, i \ge j\}$ is cofinal in $I$, so that we may compute the limit alongside $J$. For any $i \in J$, consider the following cartesian diagram : 
        
    \begin{center}
        \begin{tikzcd}[cramped]
            & U_j \ar[r] \ar[d] & U_i \ar[d] \\
            & U_j \ar[r] & X
        \end{tikzcd}
    \end{center}

    By proper base change (cf. \cite[Thm 1.2.4]{mann2022-6F}), we see that, for all $i \in J$, $g_j^* g_{i\, !} g_i^* \cF = g_j^* \cF$. Hence, the colimit along $J$ is trivial, which concludes.
\end{proof}

\begin{rk}
    Under the setup of \ref{colim_union}, a similar argument shows that and $\cF \cong \llim_{i \in I} g_{i \, *}^{M} \, g_{i \, M}^{*} \, \cF$ for any $\cF \in \Dbs(X, \O^+_X/p)$.
    
    Moreover, using a technical lemma on $(\infty, 2)$-categories \cite[D.4.7]{Mann_Heyer_Smooth}, we may improve the adjunction triangle $g^{M}_{i \,!} \dashv g_{i, M}^{*} \dashv g_{i \, *}^{M}$ to the expected adjunctions between the following three functors : 
    \begin{itemize}
        \item The map $(\cF_i)_{i \in I} \in \llim \Dbs(U_i, \O^+_{U_i}/p)  \mapsto \rlim_{i \in I} g_{i \, !} \, \cF_i \in \Dbs(X, \O^+/p)$ 
        \item A right adjoint $ \cF \in \Dbs(X, \O^+/p) \mapsto (\cF \mapsto (g_i^* \cF))_{i \in I} \in \llim_{i \in I} \cD_\square^a(U_i, \O^+_U/p)$
        \item A right adjoint of the above map $(\cF_i)_{i \in I} \in \llim_{i \in I} \Dbs(U_i, \O^+_{U_i}/p)  \mapsto \llim_{i \in I} g_{i \, *} \, \cF \in \Dbs(X, \O^+/p) $ 
    \end{itemize}
    Here, the limits $\llim_{i \in I} \Dbs(U_i, \O^+_{U_i}/p)$ is taken along pullbacks. The above proposition shows that the counit of the first adjunction and the unit of the second one are isomorphisms, both these statements being equivalent (\cite[\href{https://kerodon.net/tag/02FF}{Corollary 02FF}]{kerodon}), to the full faithfulness of the middle arrow.
\end{rk}

\smallskip
Let us now prove that partially proper open immersions preserve discrete objects. We start by doing so for totally disconnected perfectoid spaces. Note that the proof follows roughly \cite[Lemma 6.2.1 and 6.2.2]{Zavyalov_Abstract_PD}, which moreover assumes $U$ to be a Zariski open subset.

    \begin{lemma}
        \label{oi_totdisc_discrete}
        Let $X$ be a totally disconnected perfectoid space over $K$ and $j : U \to X$ be a partially proper open immersion. Then $j_!$ preserves discrete objects.
    \end{lemma}
    
    \begin{proof}
    The map $j$ is partially proper, hence specializing (by the valuative criterion), in the sense that the image $j(U)$ is stable under specialization. We let $V = |j(U)|$.
    
    Consider the natural quotient map $\pi : |X| \to \pi_0(X)$. Since $X$ is a totally disconnected perfectoid, by definition, every connected component of $|X|$ admits a unique closed point, to which any other point in the connected component specializes. Since $V$ is open and stable under specializations, $V$ contains the connected component of every $x \in U$. Hence, $V \cong \pi^{-1} (\pi(V))$. Since $\pi$ is a quotient map, it follows that $\pi(V)$ is open in $\pi(X)$.
    
    By \cite[\href{https://stacks.math.columbia.edu/tag/0906}{Lemma 0906}]{stacks-project}, $\pi_0(X)$ is profinite, so that it admits a basis of clopen subsets. Thus, $\pi(V)$ can be written as a filtered union of quasi-compact clopen sets ; so that the same holds for $V = \pi^{-1}(\pi(V))$. Hence, there exists an filtered poset $I$ such that $V = \bigcup_{i \in I} V_i$ is a filtered union of clopen subsets of $|X|$. For all $i \in I$, let $U_i = j^{-1}(V_i) \subset U$, which admits an analytic structure as an open subset of $U$.

    Consider the map $j_i : U_i \to X$, which is the inclusion of a clopen subset of $X$. Let $\cF \in \Dbs(U, \O^+_U/p)_\omega$.

    By the lemma \ref{colim_union}, $\cF \cong \rlim_i \, j_{i \, !} \, j_i^* \cF$, so that $j_! \cF \cong j_! \rlim_i \, j_{i \, !} \, j_i^* \, \cF \cong \rlim_i \, (j \circ j_i)_!\, j_i^* \, \cF$,\footnote{Note that, by unraveling the proof of \cite[Lemma 3.6.2]{mann2022-6F}, this is actually the definition of $j_!$ for arbitrary open immersions.} since $j_!$ commutes with colimits (as it admits a right adjoint).

    By definition, pullback preserve discrete objects, and $j \circ j_i$ is the closed immersion of an open subset of $Y$, so that $(j \circ j_i)_! = (j \circ j_i)_*$, which, by \cite[3.3.10.(ii)]{mann2022-6F}, preserves discrete objects. Finally, by \cite[lemma 3.2.19.(ii)]{mann2022-6F}, (small) colimits of discrete objects remain discrete, so that $j_! \cF$ is discrete.
    \end{proof}

Let us now establish the result for general partially proper open immersions.

\begin{lemma}\label{openimm+pp}
    Let $j : U \to X$ be a partially proper open immersion of locally spatial diamonds over $K$. Then $j_{!}$ preserves discrete objects. 
\end{lemma}

\begin{proof}  
    Fix a discrete object $\cF \in \cD_\square^a(U, \O_U^+/p)$. We will reduce to the situation where $X$ is a totally disconnected perfectoid space.
    
    As discreteness is defined by descent (cf \cite[Def 3.2.17]{mann2022-6F}), we need to show that, for every morphism $f : \widetilde{X} \to X$ from a totally disconnected perfectoid space $\widetilde{X}$, the sheaf $f^* j_! \cF$ is discrete in $\Dbs(\widetilde{X}, \O^+_{\widetilde{X}}/p)$. Consider the cartesian square below :
    \vspace{-.9cm}
        \begin{center} 
        \[\begin{tikzcd}[cramped]
	\widetilde{U} & \widetilde{X} \\
	U & X
	\arrow["{j'}", from=1-1, to=1-2]
	\arrow["{f'}"', from=1-1, to=2-1]
	\arrow["f", from=1-2, to=2-2]
	\arrow["j", from=2-1, to=2-2]
    \end{tikzcd}\]
        \end{center}

    By proper base change \cite[Thm 1.2.4.(iv)]{mann2022-6F}, $f^* j_! \cF = j'_! f'^* \cF$. 

    By definition \cite[Defi 10.7]{Scholze_ECohD}, $j'$ is a separated open immersion. Note that the valuative criteria from \cite[Defi 18.4]{Scholze_ECohD} is automatically stable under base change, so that $j'$ remains partially proper, and, by the above lemma \ref{oi_totdisc_discrete}, $j'_!$ preserves discrete objects.

    Finally, pullbacks preserve discrete objects by definition, so that $j'_! f'^{*} \cF$ is discrete, which concludes.
\end{proof}

We may now prove the main result \ref{preservesdisc}.

\begin{proof} \label{proofcountable} (of \ref{preservesdisc})
    Using \ref{VNiceCov}, we write $X = \bigcup_{i \in I} X_i$ as an increasing filtered cover by qcqs open subsets, such that, for any $i < j$, the inclusion $X_i \subset X_j$ factors through a partially proper $U_i$ and such that the map $X_i \to X_j$ factors through Huber's universal compactification $\overline{X_i}^Y$ of $X_i$ over $Y$. Denote $f_i : X_i \to X$ and $g_i : U_i \to X$ the associated open immersions.

    \vspace{.05cm}
    By lemma \ref{colim_union}, there is an isomorphism $\cF \cong \rlim_{i \in I} g_{i \,!} \, g_i^* \, \cF$.
    
    Hence, $f_! \, \cF = f_! \left( \rlim_{i \in I} g_{i\, !} \, g_i^* \cF \right) = \rlim_{i \in I} (f \circ g_i)_! \, g_i^* \cF$, since left adjoints commute with colimits. By \cite{mann2022-6F}, Lemma 3.2.19.(ii), a (small) colimit of discrete objects remains discrete, so that, as pullbacks preserve discrete objects, it suffices to show that every $(f \circ g_i)_!$ preserves discrete objects.

    \vspace{.05cm}
    Here, $g_i : U_i \to X$ is the inclusion of a partially proper open subset in $X$. Moreover, the partially proper map $f \circ g_i$ factors through $X_i$, then, for any $i < j$, through the adic compactification $\bar{X}_{j}^Y$ of $X_j$ over $Y$ by the universal property of universal compactifications \cite[Defi 5.1.1]{Huber_Ecoh}. We then have the following diagram, where $U_i \to Y$ is partially proper : 
    \begin{center}
    \begin{tikzcd}[row sep = large]
        & U_i \arrow[rr, "g_i", bend left]\arrow[drr, swap] \arrow[r, hookrightarrow] & \bar{X}_{j}^Y \arrow[r, hookrightarrow] \arrow[dr]& X \arrow[d, "f"]
    \\ & & & Y
    \end{tikzcd}
    \end{center} 
    By \cite[Corro 5.1.6]{Huber_Ecoh}, the map $\bar{X}_{j}^Y \to Y$ is proper, so that, by the two-out-of three property for partially proper morphisms \cite[Lemma 1.10.17.(vi)]{Huber_Ecoh}, the open immersion $U_i \to \bar{X}_{j}^Y$ is additionally partially proper.
    
    Hence, the composition $g_i \circ f$ factors as the composition of a partially proper open immersion and a proper morphism, and both $g_i$ and $f$ preserve discrete objects by the lemma \ref{openimm+pp} and \cite[Lemma 3.3.10.(ii)]{mann2022-6F} respectively. This concludes.
\end{proof}

\subsection{Localization sequence in $\Dbs(-, \O^+/p)$}

We will use the following result, which is an analogue, in Mann's formalism, of the excision result proven in \ref{excision_diamond}.

\begin{pro}
    \label{Excision_M}
    Let $X$ be an untilted locally spatial diamond over $K$, together with a partially proper open immersion $j : U \to X$. By \ref{ComplementLSD}, the complement $F = |X| \setminus j(|U|)$ admits a structure of a locally spatial diamond, and we let $i : F \to X$ be the corresponding closed immersion. 

    Let $\cF \in \Dbs(X, \O^+_X/p)$. Then, there is a natural exact sequence in $\Dbs(X, \O^+/p)$ : $$ j_!^{M} j^*_{M} \cF  \to \cF \to i_*^{M} i^*_{M} \cF .$$
\end{pro}

Note that the proof resembles the one of\cite[Lemma 6.2.2]{Zavyalov_Abstract_PD}, which only considers Zariski open subsets.

\begin{proof}
    The morphisms are naturally (co)-unit morphism, so that it suffices to prove the exactness of the sequence. This can be checked on a basis, i.e. after pullback to any totally disconnected perfectoid space. By proper base change, we may assume that $X$ admits a map to a perfectoid space\footnote{This hypothesis is needed to use \ref{Comparison_Discobj}.}. 
    
    Since $\O^{+}_X/p$ is the unit of Mann's tensor product $\otimes^M$ on $X$, we may write $\cF = \cF \otimes^M \O_X^{+}/p$ and, using the projection formula (cf. \ref{6ffprop}), the sequence reduces to : $$ \cF \otimes^M j_!^{M} j_{M}^* \O^+_X/p \to \cF \to \cF \otimes^M i_*^{M} i_{M}^* \O^+_X/p.$$

    Since tensoring with $\cF$ is exact, we reduced to the case $\cF = \O_X^+/p$, which is discrete. By the above proposition \ref{preservesdisc} and \cite[Lemma 3.3.10]{mann2022-6F}, all objects appearing in the triangle are discrete. Recall that there is a functorial equivalence of categories  $\Dbs(X, \O^+_X/p)_\omega \cong \cD_{\e t}(X, \O_X^{+, a}/p)^{oc}$, as recalled in \ref{Comparison_Discobj}, so that we can check that it is exact as a sequence of almost étale sheaves, which are étale sheaves valued in an abelian category.
    
    Since the étale site $X_{\e t}$ admits enough points by \cite[Prop 14.3]{Scholze_ECohD}, it suffices to check the isomorphism after pullback to geometric points of the form $\operatorname{Spa}(C, C^+)$, for some complete algebraically closed $C$, with open bounded valuation subring $C^+$. Here, the pullback is to be taken in the quasi-coherent sense, and it coincides with $f_{M}^*$ by \ref{Comparison_Discobj} (we implicitely use the identification of \ref{Comparison_Discobj}).
    Let $f : \operatorname{Spa}(C, C^+) \to X$ be such a morphism.
    
    By proper base change applied to the same diagram as in the proof of \ref{excision_diamond}, we reduce to the case where $X$ is a totally disconnected perfectoid space, and $U$ is a partially proper open subset of $X$, and $F$ is its complement. Following the proof of \ref{excision_diamond}, $U$ is either empty or $U = X$, so that, in both cases, the exactness is trivial.
\end{proof}

Let us now construct the announced exact sequence.

\begin{lemma}
    Let $f : X \to \operatorname{Spa}(K, \O_K)$ be a partially proper analytic adic space. Let $(U_i)_{i \in I}$ be an open covering of $X$ by partially proper open subsets, as in \ref{VNiceCov}.
    
    For any $i \in I$, we let $Z_i$ be the complement of the partially proper $U_i$, viewed as a closed locally spatial subdiamond of $X^\diamond$ by \ref{ComplementLSD}. Let $i_{Z_i} : Z_i \to X^{\diamond}$ be the associated closed immersion. Then, for $\cF \in \Dbs(X, \O^+_X/p)$, there is a natural exact sequence :
    $$f_!^{M} \cF \to f^{M}_* \cF \to \rlim_{i \in I} f_*^{M} \, i^{M}_{Z_i,*} \, i_{Z_i}^{*, M} \cF$$ 
\end{lemma}

For improved readability of the following proof, everything is implicitely taken in the sense of Mann.

\begin{proof} 
    Let $j_{U_i} : U_i \to X$ denote the associated open immersions. Starting from the excision triangle from the above proposition \ref{Excision_M}, we apply $f_*$ and take colimits, as in the proof of \ref{Fiber_H}. This yields the following exact sequence : $$\rlim_{i \in I} f^{}_* \circ j^{}_{U_i \, !} \circ j_{U_i}^{*} \cF \to  f^{}_* \cF \to \rlim_{i \in I} f^{}_* \circ  i_{Z_i \, *}^{} \circ i_{Z_i}^{*} \cF$$
    
    Let $(X_i)_{i \in I}$ be a covering of $X$ by quasi-compact opens subsets, associated to $(U_i)_{i \in I}$ as in \ref{VNiceCov}. 
    
    Let us now establish that $f_* \circ j_{U_i \, !} \cong f_! \circ j_{U_i \, !}$. 
    
    \vspace{-.1cm}
    Consider the following diagram for any $i < j \in I$, where we let $Y := \operatorname{Spa}(K, \O_K)$, and $\overline{X_j}^Y$ denote the universal adic compactification of $X_j$ over $Y$ :
    \begin{center}
    \begin{tikzcd}
        & U_i \ar[r, "j_{U_i}"] \ar[dr, swap, "k_{i,j}"] & X \ar[r, "f"] & Y \\
        & & \overline{X_j}^Y \ar[u, "\pi_j"] \ar[ur, "p_j", swap] 
    \end{tikzcd}
    \end{center}
    By \cite[Corro 5.1.6]{Huber_Ecoh}, since $X_j$ is quasi-compact, $p_j$ is proper, and so is $p_j$ by \cite[Lemma 1.10.17.vi)]{Huber_Ecoh} since $f$ is separated. Hence : $$f_! \circ j_{U_i \, !} = p_{j\, !} \circ k_{i,j \, !} = p_{j\, *} \circ k_{i,j \, !} = f_* \circ \pi_{j \, !} \circ k_{i,j \,!} = f_* \circ j_{U_i \, !}.$$

    Hence, using \ref{colim_union} and the fact that $f_!$ commutes with colimits (as it admits a right adjoint), we have $$\rlim_{i \in I} f_* \,  j_{U_i \, !} \, j_{U_i}^{*} \cF \cong \rlim_{i \in I} f_! \, j^{}_{U_i \, !} \,  j_{U_i}^{*} \cF \cong f_! \left( \rlim_{i \in I} j^{}_{U_i \, !} \, j_{U_i}^{*} \cF \right) \cong f_! \cF$$  
    So that the first term of the above exact sequence identifies with $f_! \cF$, which concludes.
\end{proof}

Let us now apply the discretization functor, to get the sequence announced in \ref{foreshadow}.

\begin{corro}
    \label{Mann_Discrete_exactseq}
    In the situation above, there is a natural exact sequence :   
    $$f_!^{M} \cF \to \left(f^{M}_* \cF\right)_\omega \to \rlim_{i \in I} \left(f^{M}_* i^{M}_{Z_i,*}  i_{Z_i}^{*, M} \cF\right)_\omega$$
\end{corro}

\begin{proof}
    Apply the exact functor $(-)_\omega$ from \ref{Discretization_def} to previous result, which commutes with colimits by the lemma \ref{Commute_Colimits} below. The object $f_!^{M} \cF$ is discrete by \ref{preservesdisc}, which concludes.
\end{proof}

\begin{rk}
    As in \ref{corro_Fib_H}, we may rewrite the right hand side as a colimit alongside any complement of a reasonable covering - in particular a quasi-compact one, using a cofinality argument. 
\end{rk}

In the proof above, we use the following lemmas. They are purely technical, and can be skipped in first lecture.

\begin{lemma}
    \label{Notalmost_commutes_colimits}
    Let $\cC$ be an $\infty$-category that admits colimits and such that finite products coincide with finite coproducts (for example, an additive category). Let $\kappa$ be a strong limit cardinal.
    
    Let $\operatorname{Cond}(\cC)$ be the category of condensed objects in $\cC$ as constructed in \cite{Scholze_Condensed}, and $(-)_\kappa :\operatorname{Cond}(\cC) \to \operatorname{Cond}(\cC)_\kappa$ be the functor computed by restricting to $\kappa$-small extremally disconnected sets. 
    
    The functor $(-)_\kappa$ commutes with colimits.
\end{lemma}

\begin{proof}
    Let $(\cF_i)_{i \in I} \in \operatorname{\operatorname{Cond}(\cC)} = \rlim_{\lambda} \operatorname{Cond}(\cC)_\lambda$, where the colimit is computed along strict limit cardinals. Recall from \cite[Section 2]{Scholze_Condensed} that objects of $\operatorname{Cond}(\cC)$ can be seen as $\cC$-valued sheaves from the category of extremally disconnected compact Hausdorff sets, and colimits are given by sheafifying the pointwise colimit. 
    
    The sheaf condition for a presheaf on such extremally disconnected rewrites as mapping finite disjoint unions to finite products. Since finite products coincide with finite direct coproducts in $\cC$, the pointwise colimit of sheaves remains a sheaf, and there is no need to sheafify.  
    
    Hence, for all $\kappa$-small extremally disconnected set $S$, the evaluation $\cF \in \operatorname{Cond}(\cC) \mapsto \cF(S)$ commutes with colimits, so that $(-)_{\kappa}$ commutes with colimits.
\end{proof}

\begin{lemma}
    \label{Commute_Colimits}
    Let $(A, A^+)$ be a perfectoid pair with pseudo-uniformizer $\pi$, and $\kappa$ be a strong limit cardinal. Let $\Dbs(A^+/\pi)$ be Mann's category of solid almost $A^+/\pi$-modules, and $(-)_{\kappa, a} : \Dbs(A^+/\pi) \to \Dbs(A^+/\pi)_\kappa$ be the restriction to $\kappa$-small objects constructed in \cite[Defi 3.2.20]{mann2022-6F}.
    
    Then $(-)_{\kappa, a}$ commutes with colimits.
\end{lemma}

\begin{proof}
    Since $A^+/\pi$ is discrete, it follows from \cite[Defi 2.3.28]{mann2022-6F} that the almostification functor $(-)^a : \cD_\square(A^+/\pi) \to \Dbs(A^+/\pi)$ restricts to a functor $(-)^{a,\kappa} : \cD_\square(A^+/\pi)_\kappa \to \Dbs(A^+/\pi)_\kappa$, that is essentially surjective by definition \cite[Defi 2.3.3]{mann2022-6F}. Consider the following commutative diagram :
    \begin{center}
    \begin{tikzcd}[column sep = small]
        & \cD_\square(A^+/\pi)_\kappa \ar[d,swap, "(-)^{a,\kappa}"] \ar[r, "I_\kappa"]& \cD_\square(A^+/\pi) \ar[d, "(-)^a"] \\ & \Dbs(A^+/\pi)_\kappa \ar[r, "I^a_\kappa"] & \Dbs(A^+/\pi) 
    \end{tikzcd}.
    \end{center}
    
    Here, $I_\kappa$ and $I_\kappa^a$ are the inclusion of full subcategories, and the functor $(-)_{\kappa, a}$ that we are interested in is defined as the right adjoint of $I_\kappa^a$. Recall that the vertical arrows are localization functors, and that $(-)^a$ admits both a fully faithful left adjoint $(-)_!$ and a fully faithful right adjoint $(-)_*$. Likewise, $(-)^{a, \kappa}$ admits fully faithful left and right adjoints denoted respectively $(-)^\kappa_!$ and $(-)^\kappa_*$. 
    
    The morphism $I_\kappa$ admits a right adjoint, that we denote $(-)_\kappa^\square$, defined as the restriction to the full subcategories of solid objects of the functor $(-)_\kappa : \operatorname{Cond}(\cD(A^+/\pi)) \to \operatorname{Cond}(\cD(A^+/\pi))_\kappa$ considered in \ref{Notalmost_commutes_colimits}, that commutes with colimits by said lemma (here, $\cD(A^+/\pi)$ denotes the $\infty$-derived category of the abelian category of $A^+/\pi$-modules). Hence $(-)_\kappa^\square$ also commutes with colimits. 
    
    Since $(-)^{a, \kappa}$ admits the fully faithful left adjoint $(-)^\kappa_!$, the unit morphism $id \to (-)^{a, \kappa} \circ (-)_!^\kappa$ is a natural isomorphism. It then follows from the commutation of the diagram that $I_\kappa^a \cong (-)^a \circ I_\kappa \circ (-)^\kappa_!$, where all three functors in the composition admit a right adjoint, so that we may write $(-)_{\kappa, a} \cong (-)^{a, \kappa} \circ (-)_\kappa^\square \circ (-)_*$.
    
    Since the functor $(-)^\square_\kappa$ respects almost isomorphisms (as this is a pointwise condition), it follows from the universal property of localizations \cite[Prop 5.2.7.12]{Lurie_HTT} that there exists a functor $J : \Dbs(A^+/\pi) \to \Dbs(A^+/\pi)_\kappa$ such that : $J \circ (-)^a \cong (-)^{a, \kappa} \circ (-)_\kappa^\square$.
    
    As the right adjoint $(-)_*$ of $(-)^a$ is fully faithful, the counit morphism $(-)^a \circ (-)_*\to id$ is a natural isomorphism, so that $J \cong J \circ (-)^a \circ (-)_* \cong (-)^{a, \kappa} \circ (-)_\kappa^\square \circ (-)_*$, hence there is a natural isomorphism $J \cong (-)_{\kappa, a}$.  Moreover, the co-unit $(-)^a \circ (-)_! \to id$ is a natural isomorphism, so that : $$(-)_{\kappa, a} \cong J \cong J \circ (-)^a \circ (-)_! \cong (-)^{a, \kappa} \circ (-)_\kappa^\square \circ (-)_!.$$ The functors $(-)^{a, \kappa}$ and $(-)_!$ commutes with colimits as they are left adjoints, and $(-)^\square_\kappa$ commutes with colimits by the argument earlier in the proof. Hence $(-)_{\kappa, a}$ commutes with colimits.
\end{proof}
 
\subsection{Identification of Mann and Huber's proper pushforwards}

Let us finally prove \ref{Identification_lower_shriek}. Recall the statement :

\begin{pro_no}
    Let $f : X \to \operatorname{Spa}(C, \O_C)$ be a partially proper locally noetherian analytic adic space over a complete algebraically closed extension $C$ of $\Q_p$.
    
    Let $\cL$ be an overconvergent étale sheaf of $\F_p$-modules on $X$. Then, there is a natural isomorphism : $$J_C\left( f_{!}^{M} (\cL \otimes \O^{+,a}_X/p) \right) \cong \left( Rf_!^{H}(\cL \otimes_{\F_p} \O^+_X/p) \right)^a $$ in $\cD_{\e t}(\operatorname{Spa}(C, \O_C), \O^{+,a}_C/p)$, where $J_C$ denotes the identification from \ref{Comparison_Discobj}.
\end{pro_no}

\begin{proof}
    Let $(X_i)_{i \in I}$ be a covering of $X$ by partially proper open subsets, as in \ref{VNiceCov}. For any $i$, let $Z_i$ be the complement of $X_i$ in $X$, viewed as a locally spatial diamond, and let $i_{Z_i} : Z_i \to X$ be the associated closed immersion. Applying the colimit-preserving exact functor $J_C$ to the exact sequence from \ref{Mann_Discrete_exactseq} yields : 

    $$J_C\left(f_{!}^{M} (\cL \otimes \O^{+,a}_X/p)\right) \to J_C\left(f^{M}_* (\cL \otimes \O^{+,a}_X/p)_\omega\right) \to \rlim_{i \in I} J_C\left(\left(f^{M}_* \circ i^{M}_{Z_i \, *} \circ i_{Z_i, M}^{*} (\cL \otimes \O^{+,a}_X/p)\right)_\omega\right).$$ 

    Applying the colimit-preserving exact functor $(-)^a : \cD_{\e t}(\operatorname{Spa}(C, \O_C), \O^+_C/p) \to \cD_{\e t}(\operatorname{Spa}(C, \O_C), \O^{+,a}_C/p)$ to the exact sequence constructed in \ref{Exact_Seq_Hub_Qcoh}, we get the following exact sequence : 
    
    $$ \left( R f_{!}^{H} (\cL \otimes_{\F_p} \O^+_X/p)\right)^a \to \left( R f^{H}_* (\cL \otimes_{\F_p} \O^+_X/p)\right)^a \to \rlim_{i \in I} \left( R f^{H}_* \circ i^{S}_{Z_i \, *} \circ i_{Z_i, S}^{*} (\cL \otimes_{\F_p} \O^+_X/p)\right)^a.$$ 
    
    By \ref{Disc_OC_over_field}, there is a natural isomorphism between the middle terms of both sequences. 
    
    We will now prove that, for all $i \in I$, there is a natural isomorphism : $$J_C\left(\left(f^{M}_* \circ i^{M}_{Z_i \, *} \circ i_{Z_i, M}^{*} (\cL \otimes \O^{+,a}_X/p)\right)_\omega\right) \cong \left( R f^{H}_* \circ i^{S}_{Z_i \, *} \circ i_{Z_i, S}^{*} (\cL \otimes_{\F_p} \O^+_X/p)\right)^a.$$
    
    We may compute the left hand side as : 
    \begin{align*}
        J_C & \left(\left(f^{M}_* \, i^{M}_{Z_i \, *} \, i_{Z_i, M}^{*} (\cL \otimes \O^{+,a}_X/p)\right)_\omega \right) & \text{} \\ 
        & \cong Rf_*^S \, i_{Z_i, *}^S \, J_{Z_i} \left( i_{Z_i, M}^* (\cL \otimes_{\F_p} \O^{+a}_X/p)\right)  & \text{By \ref{Commute_pushforward_J}} \\
        & \cong Rf_*^S \, i_{Z_i, *}^S \, i_{Z_i, S}^{*, qcoh} \left( J_C (\cL \otimes_{\F_p} \O^{+a}_X/p)\right) &\text{By \ref{Comparison_Discobj}} \\ 
        & \cong Rf_*^S \, i_{Z_i, *}^S \, i_{Z_i, S}^{*, qcoh}  \left( \left( \cL \otimes_{\F_p} \O^{+}_X/p \right)^a \right) &\text{By \ref{compat_RH}} \\ 
        & \cong \left(  Rf_*^S \, i_{Z_i, *}^S \, i_{Z_i, S}^{*, qcoh}  \left( \cL \otimes_{\F_p} \O^{+}_X/p \right) \right)^a &\text{By \ref{Pullback}} 
    \end{align*}
    
    Moreover, $Rf_*^S \cong Rf_*^H$ since $f$ is a morphism between locally analytic analytic adic spaces, and, by \ref{Qcoh_Pullback_qproet}, $i_{Z_i, S}^{*, qcoh} \cong i_{Z_i, S}^{*}$. Hence, there is a canonical isomorphism between the rightmost terms of the exact sequences : $$ \rlim_{i \in I} J_C\left(\left(f^{M}_* \circ i^{M}_{Z_i \, *} \circ i_{Z_i, M}^{*} (\cL \otimes \O^{+,a}_X/p)\right)_\omega\right) \cong  \rlim_{i \in I} \left( R f^{H}_* \circ i^{S}_{Z_i \, *} \circ i_{Z_i, S}^{*} (\cL \otimes_{\F_p} \O^+_X/p)\right)^a$$
    
    that is compatible with the isomorphism between the middle terms. 
    
    There is then a natural isomorphism between the fibers, which concludes.
    \end{proof}

\section{Primitive Comparison with Compact Support}
\label{Sec5}

In this section, we finally define our notion of primitive comparison with compact support, and prove that it implies Poincaré duality. We start this section by recalling various improvements around the primitive comparison theorem, from the litterature.

\subsection{Primitive comparison theorem for Zariski-constructible coefficients}

Let us start by recalling the notion of Zariski-constructible complexes from \cite{BH_Sixfunctors_Zariskiconstr}.

\begin{defi} 
    Let $X$ be a rigid analytic space over a nonarchimedean field $K$. A sheaf $\cF$ of $\F_p$-vector spaces on $X_{\e t}$ is Zariski-constructible if $X$ admits a locally finite stratification $X = \bigcup_{i \in i} X_i$ into Zariski locally closed subsets $X_i$ such that, for all $i \in I$, $\cF_{| X_i}$ is locally constant and admits finite stalks.
    
    We let $D^b_{zc}(X_{\e t}, \F_p) \subset D(X_{\e t}, \F_p)$ be the full subcategory formed by bounded complexes $\cF$ whose homotopy sheaves $\pi_n(\cF)$ are Zariski-constructible for all $n \in \Z$.
\end{defi}

Note that any Zariski-constructible sheaf is automatically overconvergent. 

\begin{rk}
    \label{Morph_Proj}
Recall that, for any finite dimensional\footnote{A morphism is said to be finite-dimensional if it is of finite transcendence degree. This hypothesis is used in \cite{Huber_Ecoh} to ensure that complexes are bounded in the derived category, and is probably not needed.} partially proper morphism of locally noetherian analytic adic spaces $f : X \to Y$ and sheaf $\cF$ of $\F_p$-modules on $X$, there is a natural morphism : 
     $$ Rf_!^{H} \cF \otimes_{\F_p}^\L \O^+_Y/p \cong Rf_!^{H} (\cF \otimes_{\F_p}^\L f_{H}^* \O^+_Y/p) \to Rf^{H}_! (\cF \otimes_{\F_p} \O^+_X/p) $$
 Where the isomorphism is the projection formula of \cite[Prop 5.5.1.(iv)]{Huber_Ecoh} and the second map is induced by the natural morphism $f^*_H \O^+_Y/p \to \O^+_X/p$.
\end{rk}

In the following, we will need the following version of the primitive comparison theorem :

\begin{thm}
    \label{PrimComparison_Galois}
    Let $K$ be a complete extension of $\Q_p$, $\overline{K}$ an algebraic closure of $K$, and $C$ be the completion of $\overline{K}$. Let $f : X \to \operatorname{Spa}(K, \O_K)$ be a proper rigid analytic variety, and let $\cF \in D_{zc}^b(X, \F_p)$.
    
    Then, for all $n \in \N$, the natural morphism : $$(R^n f^H_* \cF) \otimes_{\F_p} \O^+_{K}/p \to R^n f^H_* (\cF \otimes_{\F_p} \O_{X}^+/p)$$
    is an almost isomorphism. In particular, for any Zariski-constructible sheaf $\cL$ of $\F_p$-vector spaces on $X$, the natural morphism : 
    $$H^n_{\e t}(X_C, \cF) \otimes_{\F_p} \O_C/p \to H^n_{\e t}(X_C,  \cL \otimes_{\F_p} \O_{X_C}^+/p)$$
    is a $\operatorname{Gal}(\overline{K}/K)$-equivariant almost isomorphism, where the Galois action on the left hand side is diagonal.
\end{thm}

\begin{proof}
    By \cite[Prop 3.6]{BH_Sixfunctors_Zariskiconstr}, the category $\cF \in D_{zc}^b(X, \F_p)$ is generated by objects of the form $g_* \F_p$ for finite morphism $g : Y \to X$, and, by the arguments of the proof of \cite[Lemma 7.3.6]{zavyalov2024coherentmodulescoherentsheaves}, it suffices to check the isomorphism for $\mathcal{F} = g_* \F_p$. Using \cite[Corro 7.2.8]{zavyalov2024coherentmodulescoherentsheaves}, we may replace $f$ by $g \circ f$, and it suffices to prove the statement for $\cF = \F_p$. It suffices to check the isomorphism on stalks $\operatorname{Spa(C', C'^+)}$, which reduces to $\operatorname{Spa(C', \O_{C'})}$ following the beginning of the proof of \cite{Scholze_padHTforRAV}. The result then follows from \cite[Thm 3.17]{Scholze_Perfspace_Survey}.
    
    The second statement follows from taking stalks with respect to $\operatorname{Spa}(C, \O_C) \to X$, and using proper base change as in \cite[Thm 16.10]{Scholze_ECohD}, together with the identification of pullbacks of the structure sheaf $\O^+/p$ by quasi-pro-étale morphisms as in \ref{pro_et_pullback_is_restr}.
\end{proof}

\begin{rk}
    Note that primitive comparison theorem does not apply to proper locally noetherian analytic adic spaces, see \cite[Example 5.2.5.(b)]{ACLB_Ansch_Mann_6FFonFF} for a counterexample.
\end{rk}

\subsection{Primitive comparison with compact support and Poincaré duality}

\label{Comparison_Morphism_PCCS}

Let us define our notion of primitive comparison with compact support. 

\begin{defi}
    \label{Criterion_PCCS}
    Let $f : X \to \operatorname{Spa}(C, \O_C)$ be a partially proper locally noetherian analytic adic space of finite dimension over a complete algebraically closed extension of $\Q_p$, and $\cL$ be an overconvergent sheaf of $\F_p$-vector spaces on $X$.
    
    We say that $X$ satisfies primitive comparison with compact support with respect to $\cL$ if the morphism from \ref{Morph_Proj} is an almost isomorphism, or equivalently if, for all $i \in \N$, the induced map : $$H^i_{\e t, c}(X, \cL) \otimes_{\F_p} \O_C/p \to H^i_{\e t, c}(X, \cL \otimes \O_X^+/p)$$ is an almost isomorphism of $\O_C/p$-modules.
\end{defi}

Our main result is as follows : 

\begin{thm}
    \label{PCCS_PD}
    Let $C$ be a complete algebraically closed extension of $\Q_p$, and $f : X \to \operatorname{Spa}(C, \O_C)$ be a partially proper locally noetherian analytic adic space over $C$, that is smooth of pure dimension $d$. 
    
    Let $\cL$ be an overconvergent sheaf of $\F_p$-modules on $X$, and assume that $X$ satisfies primitive comparison with compact support with respect to $\cL$.
    
    Then, for any $0 \le k \le 2d$, there exists a natural isomorphism : $$H^k_{\e t}(X, \cL^\vee(d)) \simeq \operatorname{Hom}_{\F_p}(H^{2d-k}_{\e t, c}(X, \cL), \F_p).$$
\end{thm}

\begin{proof}
    Since $X$ satisfies primitive comparison with compact support, the map : $$\left(Rf_!^{H} \cL \otimes_{\F_p}^\L \O_C/p \right)^a \to \left(R f_!^{H} (\cL \otimes_{\F_p} \O_X^+/p) \right)^a$$ is an isomorphism in $\cD_{\e t}(\operatorname{Spa}(C, \O_C), \O_C^{+a}/p)$. 
    
    By \ref{Identification_lower_shriek}, the second term canonically identifies with $J_C \left( f_{!}^{M} (\cL \otimes \O^{+,a}_X/p) \right)$.
    
    By \ref{compat_RH}, there is an isomorphism $J_C^{-1} \left((Rf^{H}_! \cL \otimes_{\F_p}^\L \O_C/p)^a \right) \cong \left( Rf^{H}_! \cL \right) \otimes \O^{+,a}_{C}/p$.
    
    Combining the above yields an isomorphism $f_{!}^{M} (\cL \otimes \O^{+,a}_X/p) \cong \left(Rf^{H}_! \cL \otimes \O_C/p \right)^a$. 
    
    We may then apply corollary \ref{Corro_PCCS_PD} with $\cF_\cL = Rf^{H}_! \cL$, which gives isomorphisms for all $0 \le k \le 2d$ : $$H^k_{\e t}(X, \cL^\vee)(d) \simeq \Hom_{\F_p}\left(\pi_{k - 2d}\left(Rf^{H}_! \cL\right), \F_p\right).$$
    
    Finally, $\pi_{k-2d} \left(Rf_!^{H} \cL\right) \cong R^{2d-k} f_!^{H} \cL \cong H^{2d-k}_{\e t, c}(X, \cL)$, which concludes.
\end{proof}

We can also make the following Galois-equivariant version.

\begin{thm}
    \label{PCCS_PD_Gal}
    Let $K$ be a complete extension of $\Q_p$ and $f : X \to \operatorname{Spa}(K, \O_K)$ be a partially proper locally noetherian analytic adic space over $K$ that is smooth of pure dimension $d$. Let $\cL$ be an overconvergent sheaf of $\F_p$-modules on $X$. 
    
    Let $C$ be the completion of an algebraic closure of $K$, and $X_C$ denote the base change of $X$ to $C$, and $\cL_C$ denote the base change of $\cL$ to $X_C$.
    
    Assume that $X_C$ satisfies primitive comparison with compact support with respect to $\cL_C$.
    
    Then, for any $0 \le k \le 2d$, there exists a $\operatorname{Gal}(\overline{K}/K)$-equivariant natural isomorphism : $$H^k_{\e t}(X_C, \cL_C^\vee(d)) \cong \operatorname{Hom}_{\F_p}(H^{2d-k}_{\e t, c}(X_C, \cL_C), \F_p).$$
\end{thm}

\begin{proof}
    The isomorphism follows from \ref{PCCS_PD} above, so that it suffices to establish Galois equivariance.

    By construction, $\cL_C$ is equivariant, and the morphism appearing in the primitive comparison with compact support is Galois-equivariant.
    
    Any $g \in \operatorname{Gal}(\overline{K}/K)$ acts on $C$ and $X_C$ via a morphism that we denote $\psi_g$, and the induced action on étale sheaves is then given by the pullback $\psi_{g, H}^*$. Using base change with respect to zero-dimensional morphisms \cite[Thm 5.4.6]{Huber_Ecoh}, it follows that $\psi_{g, H}^{*} Rf^H_! \cL \cong  Rf^H_! \psi_{g, H}^* \cL.$
    
    The equivariance now follows from \ref{Functoriality}.
\end{proof}

\subsection{Poincaré duality for almost proper varieties}

Checking primitive comparison with compact support behaves nicely with respect to excision. In particular, we have the following technical criterion.

\begin{lemma}
    \label{Test_Comple}
    Let $f : X \to \operatorname{Spa}(K, \O_K)$ be a proper rigid-analytic variety over a complete extension $K$ of $\Q_p$ with fixed algebraic closure $\overline{K}$, and let $C$ be the completion of $\overline{K}$. Let $F$ be a generalizing closed subset of $X$, endowed with the induced structure of a locally spatial diamond as in \ref{ComplementLSD}.
    
    Denote $i : F \to X$ be the closed immersion, and $h = f \circ i$. Let $U = X \setminus F$ be the complementary open subset with the induced analytic structure, and $j : U \to X$ the associated open immersion. Denote by $X_C, U_C$, and $F_C$ the base change of $X, U$ and $F$ to $C$, and $i_C, j_C$ the associated immersions. Let $f_C : X_C \to \operatorname{C}$ be the structural morphism of $X_C$, and $h_C = f_C \circ i_C$.
    
    Let $\cL$ be an overconvergent étale sheaf of $\F_p$-vector spaces on $X$, and $\cL_C$ its pullback to $X_C$. Assume that the natural morphisms : 
    \begin{align}
        \label{first_iso} \left( Rf^H_{C \, *} \cL_C \right) \otimes_{\F_p}^\L \O_C/p & \to Rf_{C *}^H (\cL_C \otimes_{\F_p} \O^+_{X_C}/p) \\
        \label{second_iso} 
         \left( R h^H_{C \, *}  i_{C, S}^* \cL_C \right) \otimes_{\F_p}^\L \O_C/p  & \to R h^H_{C \, *} (i_{C,S}^* \cL_C \otimes_{\F_p} \O^+_{F_C}/p)
    \end{align}
    
    are almost isomorphisms. Assume furthermore that $X$ is smooth of pure dimension $d$.
    
    Then, for any $0 \le k \le 2d$, there exists a natural $\operatorname{Gal}(\overline{K}/K)$-equivariant isomorphism : $$H^k_{\e t}(U_C, j^*_{C, H} \cL_C^\vee(d)) \cong \operatorname{Hom}_{\F_p}(H^{2d-k}_{\e t, c}(U_C, j_{C, H}^* \cL_C), \F_p).$$
\end{lemma}

\begin{proof}
    Recall that $j_{C, H}^* \cong j_{C, S}^*$. We can consider every pullbacks and (proper) pushforwards to be taken in the sense of Scholze, and we drop the "S" superscripts to alleviate notations. 
    
    From the localization sequence \ref{excision_diamond} applied to $\cL_C$, we get the following exact sequence in $\cD_{\e t}(X_C, \F_p)$ : 
    $$ R j_{C \, !} \, j_C^* \cL_C \to \cL_C \to R i_{C\, *} \, i_C^* \cL_C$$ 
    Likewise, by applying the localisation sequence to $\cL_C \otimes \O^+_{X_C}/p$, we get the following exact sequence in $\cD_{\e t}(X_C, \O^+_{X_C}/p)$ : $$R j_{C \, !} \, (j_{C}^{*} \cL_C \otimes_{\F_p} \O^+_{X_C}/p) \to \cL_C \otimes_{\F_p} \O^+_{X_C}/p \to Ri_{C \, *} \, i_{C}^{*} (\cL_C \otimes_{\F_p} \O^+_{X_C}/p)$$
    
    We may then apply the exact functor $Rf_{C \, *} = Rf_{C \, !}$ to both sequences, and extend the scalars of the first one to $\O_C/p$. By functoriality of localization sequences, and using the functorial morphism constructed in \ref{Morph_Proj}, we obtain the following diagram with exact rows : 
    \begin{center}
    \begin{tikzcd}[column sep = small]
        R (f_C \circ j_C)_! \, j_C^* (\cL_C \otimes_{\F_p} \O^+_{X_C}/p) \ar[r] 
        & R f_{C \,*} (\cL_C \otimes_{\F_p} \O^+_{X_C}/p) \ar[r] 
        & R (f_C \circ i_C)_* \, i_C^* (\cL_C \otimes_{\F_p} \O^+_{X_C}/p)
        \\ 
        \left( R(f_{C} \circ j_C)_! \, j_C^* \cL_C \right) \otimes^\L_{\F_p} \O_{C}/p \ar[r] \ar[u]
        & Rf_{C *} \cL_C \otimes^\L_{\F_p} \O_{C}/p \ar[r] \ar[u]
        & (R(f_C\circ i_C)_* \, i_C^* \cL_C) \otimes^\L_{\F_p} \O_{C}/p \ar[u]
    \end{tikzcd}
    \end{center}
    
    By the hypothesis (\ref{first_iso}), the middle vertical arrow induces an almost isomorphism. 
    
    Using (\ref{second_iso}), the rightmost vertical arrow is an almost isomorphism. Indeed, pullbacks commute with tensor products, and, by \ref{Pullback_qproet} (since generalizing closed immersions are quasi-pro-étale), we have a canonical isomorphism : $$R h_{C *} \, i_C^* (\cL_C \otimes \O^+_{X_C}/p) \cong R h_{C \, *} \left( i_C^* \cL_C \otimes \O^+_{F_C}/p\right).$$

    Hence the leftmost vertical arrow is an almost isomorphism. By the same argument as the one used for the right hand side, this corresponds to the fact that the following morphisms : $$H^i_{\et, c}(U_C, j^* \cL_C \otimes \O^+_{X_C}/p) \to  H^i_{\et, c}(U_C, j^* \cL_C) \otimes \O_C/p)$$ are almost isomorphisms for all $i \in \N$. We conclude by theorem \ref{PCCS_PD_Gal}.
\end{proof}

We may deduce Poincaré Duality for \textit{almost proper} (also known as Zariski compactifiable) rigid analytic varieties. Doing so, we recover a result by Zavyalov-Li-Reinecke \cite[Corro 1.1.2]{Zavyalov&co}.

Note that the approach of loc.cit. uses different tools, as it barely uses perfectoid methods, and, rather, relies on a study of formal models, and of trace formalisms. 

\begin{corro}
    \label{Almost_proper}
    Let $K$ be a complete extension of $\Q_p$, $\overline{K}$ a fixed algebraic closure, and $C$ be the completion of $\overline{K}$. Let $X$ be a smooth proper rigid analytic variety over $K$, and $U \subset X$ be a Zariski open subset. Let $j : U \to X$ denote the associated open immersion.
    
    Then, for any étale $\F_p$-local system $\cL$ on $U$, $U_C$ satisfies Poincaré duality with respect to $\cL_C$, i.e there is a natural $\operatorname{Gal}(\overline{K}/K)$-equivariant isomorphism : $$H^k_{\e t}(U_C, \cL_C^\vee(d)) \cong \operatorname{Hom}_{\F_p}(H^{2d-k}_{\e t, c}(U_C, \cL_C), \F_p)$$
\end{corro}

\begin{proof}
    Let $F = X \setminus U$, endowed with the reduced structure, making it a proper rigid analytic variety. Let $\cL_C^X := j_!^H \cL_C$. Since $j$ is a Zariski open immersion, it follows from \cite[Prop 3.11]{BH_Sixfunctors_Zariskiconstr} that $j^H_{!} \cL_C$ is Zariski constructible on $X_C$. 
    
    Moreover, it follows from base change \cite[Thm 19.1]{Scholze_ECohD} that $j_{C, H}^* j_{C \, !}^H \cL_C \cong \cL_C$.
    
    We will use the criterion from \ref{Test_Comple} (and the associated notations) to the overconvergent sheaf $j_!^H \cL$. By base change, $i^*_{C, S} j_{C, !}^H \cL = 0$, so that the second almost isomorphism is trivial. It then suffices to show that the following morphism is an almost isomorphism :
    \begin{align*}
        \left( Rf^H_{C \, *} \, j^H_{C \, !} \, \cL_C \right) \otimes_{\F_p}^\L \O_C/p  \to Rf_{C *}^H (j^H_{C\,!} \, \cL_C \otimes_{\F_p}^\L \O^+_{X_C}/p).
    \end{align*}

    This follows directly from the primitive comparison theorem as stated in \ref{PrimComparison_Galois}.
\end{proof}

The criterion \ref{Test_Comple} also allows us to consider other setups, such as Drinfeld's upper half plane. Let us start with the following lemma : 

\begin{lemma}
    Let $K$ be a complete extension of $\Q_p$, and $C$ be the completion of its algebraic closure. Let $X$ be a smooth proper rigid analytic variety over $K$. Let $S$ be a profinite topological space, viewed as a locally spatial diamond over $K$, together with a closed immersion $S \to X$. Let $U = X \setminus S$, endowed with the induced analytic structure.
    
    Let $U_C$ be the base change of $U$ to $C$. Then $U_C$ satisfies primitive comparison with compact support with respect to the constant local system $\F_p$.
\end{lemma}

\begin{proof}
    Recall that, by definition of the locally spatial diamond associated to a profinite space, the base change of $S$ to $C$ admits the structure of an affinoid perfectoid given by $S_C = \operatorname{Spa}(C^0(S, C), C^0(S, \O_C))$. 
    
    Write $S$ as a projective limit $S = \llim_{i \in I} S_i$ of finite discrete spaces. By \cite[Prop 14.9]{Scholze_ECohD}, there is a natural isomorphism $H_{\e t}^j(S_C, \cF) \cong \rlim_{i \in I} H_{\e t}^j(S_{i,C}, \cF_i)$ for any étale sheaf $\cF$ on $S_C$, where $\cF_i$ denotes the pullback of $\cF$ to $S_{i, C}$, viewed as a closed subspace of $X_C$.
    
    For any $i \in I$, $S_{i,C}$ is a strictly totally disconnected affinoid perfectoid space (it is a finite disjoint union of some $\operatorname{Spa}(C, \O_C)$), so that étale cohomology coincides with sheaf cohomology of the associated topological space. Since the $S_{i}$ are finite topological spaces, there is no higher cohomology, and $H^0$ is simply given by global sections. Hence,$H^j_{\e t}(S_C, \cF) = 0$ for $j > 0$.
    
    Moreover $H^0_{\e t}(S_C, \F_p) \cong \rlim_{i \in I} \operatorname{Cst}(S_{i}, \F_p) \cong \cC^0(S, \F_p)$ where $\operatorname{Cst}$ denotes the ring of constant functions,  and $$H^0(S_C, \O^+_{S_C}/p) \cong \rlim_{i \in I} \operatorname{Cst}(S_{C,i}, \O_C/p) \cong \cC^0(S_C, \O_C/p)$$ since $\O_C/p$ is discrete.
    
    The natural map $\cC^0(S_C, \F_p) \otimes \O_C/p \to \cC^0(S_C, \O_C/p)$ is naturally an isomorphism, which concludes.
\end{proof}

As a corollary, we deduce Poincaré duality for Drinfeld's upper half plane.

\begin{corro}
    Let $C = \C_p$ be the completion of an algebraic closure of $\Q_p$, and $\Omega_{C} = \P^1_C \setminus \P^1(\Q_p)$ denote Drinfeld's upper half plane over $C$.
    
    Then $\Omega_{C}$ satisfies Poincaré duality, i.e., for $0 \le k \le 2$, there is a natural $\operatorname{Gal}(\overline{\Q_p}/\Q_p)$-equivariant isomorphism : $$H^k_{\e t}(\Omega_C, \F_p(1)) \cong \operatorname{Hom}(H^{2-k}_{\e t, c}(\Omega_C, \F_p), \F_p).$$
\end{corro}

\begin{proof}
    This follows directly from the above proposition and \ref{PCCS_PD_Gal}.
\end{proof}

Note that one can easily check that the above morphism is also $\operatorname{GL}_2(\Q_p)$-equivariant using the functoriality established in \ref{Functoriality}. 

In the rest of this paper, we prove that such a Poincaré duality result also holds for Drinfeld's symmetric spaces of higher dimension, as well as general period domains.

\section{Period domains}
\label{PerDom}
\label{Sec6}

In this section, we give a precise definition of period domains, for which we claim no originality. We start by specializing the theory for $\GL_n$. Such a context is sufficient to define Drinfeld's symmetric spaces, is necessary for the general theory, and contains most of the necessary intuition.

The original construction was developed  in \cite{RappZink}, see also a short introduction by the first author in \cite{RappIntro}. Concise introductions, somewhat similar to the one below, can be found in \cite{Orlik_local}, \cite{cdhn}. The book \cite{dat_orlik_rapoport_2010} presents an extensive version of the theory, whose structure we broadly follow. 

A more modern interpretation using the Fargues-Fontaine curve is laid out in \cite{Fargues_Shen_Chen}, however we do not expand on this point of view.



\subsection{Filtered isocrystals, and period domains for $\GL_n$}


Period domains for $\GL_n$ are classifying spaces of \textit{weakly admissible} isocrystals in the sense of Fontaine, which are viewed as \textit{semi-stable objects of slope zero} in an adequate Harder-Narasimhan formalism on the category of filtered isocrystals. In what follows, we fix a field $K$.

In the following, we assume some familiarity with Tannakian formalism, and refer the result to \cite{Deligne1982_Tannaka} or \cite[Section 4]{dat_orlik_rapoport_2010}.
\subsubsection{Filtered vector spaces}
    \label{filtered}
    Let $V$ be a finite dimensional $K$-vector space. A $\Q$-filtration $\cF V$ on $V$ is the datum of a map $\mathcal{F} : x \in \Q \mapsto \mathcal{F}^x V \in \{ K\text{-linear subspaces of V}\}$ such that :
    \begin{itemize}[topsep=0pt]
        \item $\cF$ is decreasing, i.e, for all $x < y$, $\cF^y \subset \cF^x$.
        \item There exists $x<y$ such that $\cF^x V = V$ and $\cF^y V = \{0\}$.
        \item For all $x$, $\mathcal{F}^x = \bigcap_{y < x} \cF^y$.
    \end{itemize}


    We define the \textbf{type} of the filtration $\cF$ to be $\nu(\cF) = \left( x_1^{n(1)}, \dots, x_n^{n(i)}\right)$\footnote{Here, the exponent $n_i$ indicates that the corresponding $\nu(i)$ is repeated $n_i$ times.}, where $x_1 > \dots > x_n$ are the jumps of the filtration, i.e. $\{ x_1, \dots, x_n\} = \left\{ x \in \R, gr^{x_i}_{\cF} \neq \{0\}\right\}$, and $n_i = \operatorname{dim}_V(gr^{x_i}_\cF)$ is the multiplicity at $x_i$.

    The datum of $\left( x_1^{n(1)}, \dots, x_n^{n(i)}\right)$, where $x_1 > \dots > x_n$ and $\sum n_i = dim(V)$ is called the possible type of a filtration on $V$.
        
\begin{defi}
    \label{flag_gln}
    Let $V$ be a finite dimensional $K$-vector space, and let $\nu$ be the possible type of a filtration on $V$. The flag variety (cf. \cite[1.31]{RappZink}) $\mathscr{F}(\GL(V), \nu)$\footnote{ The more standard notation (cf. \cite{dat_orlik_rapoport_2010}) is $\mathscr{F}(V, \nu)$. We prefer denoting $\GL(V)$ so that our notation coincides with the general case, where $\GL(V)$ will be replaced by a reductive algebraic group.} parameterizes filtrations of type $\nu$ over $K$, i.e. it is the projective $K$-scheme such that, for any field extension $K'$ over $K$ : $$\mathscr{F}(\GL(V), \nu)(K') \cong \{ \text{filtrations on } V \otimes_K K' \text{ of type } \nu \}.$$

    It can be constructed as $\mathscr{F}(\GL(V), \nu) = \mathbb{GL}(V)/P(\nu)$, where $\GL(V)$ is the reductive group defined by $R \mapsto \operatorname{GL}(V \otimes_K R)$ for all $K$-algebras $R$, and $P(\nu)$ is the parabolic subgroup of $\GL(V)$ stabilizing a given filtration of type $\nu$ on $V$.
    
    It is naturally endowed with an action of the linear group $GL(V) =\GL(V)(K)$. 
\end{defi}

    For $K/k$ a field extension, we let $\Q-\operatorname{Fil}_k^K$ denote the category of pairs $(V, \cF_K)$, where $V$ is a $k$-vector space, and $\cF_K$ is a $\Q$-filtration on the $K$-vector space $V \otimes_k K$.
    A morphism in $\Q-\operatorname{Fil}_k^K$ is a $k$-linear map $f : V \to W$ such that, for all $x \in \Q$, $f \otimes id_K(V^x) \subset W^x$.

    The category $\operatorname{Fil}_k^K$ is an additive $k$-linear $\otimes$-category. It is not abelian (as a change of the filtration cannot be read on the kernel nor cokernel), but it is \textit{quasi-abelian}, and is in particular an exact category. 
    


\subsubsection{Isocrystals} \label{Isoc_sec}

\label{isoc}
    Let $L$ be a perfect field of characteristic $p > 0$, $K_0 = \operatorname{Frac}(W(L))$, and $\sigma$ be the Frobenius endomorphism on $K_0$. 
    
    An isocrystal\footnote{This is sometimes referred to as F-isocrystals, where the  "F" reminds us of the presence of a Frobenius morphism. We simply call them Isocrystals for clarity of notation. Dear reader, please remind yourself at all times that isocrystals come endowed with a Frobenius morphism.} over $K_0$ is a pair $N = (V, \varphi)$, where $V$ is a finite-dimensional vector space over $K_0$, and $\varphi : V \to V$ is a $\sigma$-linear\footnote{i.e. $\forall \lambda \in K_0$, $\forall x \in V$, $\varphi(\lambda\cdot x) = \sigma(\lambda) \cdot \varphi(x)$.} bijective automorphism.

    We let $\operatorname{Isoc}(K_0)$ denote the category of such isocrystals, with morphisms given by $K_0$-linear maps $f : V \to V'$ such that $f \circ \varphi_V = \varphi_{V'} \circ f$. It is an abelian category, and, moreover, we have the following (cf. \cite{Manin_1963}) :

\begin{thm}(Dieudonné-Manin classification) \label{DN} \\ 
Assume that $L$ is algebraically closed. Then, $\operatorname{Isoc}(K_0)$ is semi-simple, and the simple objects are of the form :
$$V_{r/s} = \left(K_0^s, \begin{bmatrix}
        0 & 1 &  \\ 
          & \ddots & 1 \\
        p^r &  &   0 & 
    \end{bmatrix} \cdot \sigma\right)$$

for any coprime $(r, s) \in \Z \times \N_{\ge 1}$. 

More explicitly, any $N \in \operatorname{Isoc}(K_0)$ can be uniquely written as $V = \bigoplus_{x \in E} V_x$ for some finite $E \subset \Q$. Using the terminology described above, the decomposition as a sum of simple objects yields defines a grading functor $\omega : \Isoc(K_0)\to \Q-\operatorname{Grad}_{K_0}$, where $\Q-\operatorname{Grad}_{K_0}$ denotes the category of finite dimensional $K_0$-vector spaces equipped with a grading indexed by $\Q$.
\end{thm}

\subsubsection{Filtered isocrystals} \label{FilIC}
    We keep the notation of \ref{Isoc_sec}, and let $K$ be a field extension of $K_0$. A \textbf{filtered} isocrystal $(V, \varphi, \cF V_K)$ over $K_0$ is the datum of an isocrystal $(V, \varphi)$ over $K_0$, together with a $\Q$-filtration $\cF V_K$ on the vector space $V_K = V \otimes_{K_0} K$.

    We let $\operatorname{FilIsoc}_{K_0}^K$ denote the category of such filtered isocrystals. The category $\operatorname{FilIsoc}_{K_0}^K$ is endowed with a Harder-Narasimhan formalism, defined by :   

    \begin{itemize}
        \item The \textbf{rank} of $(V, \varphi, \cF V_K)$ is $dim_{K_0} V$,
        \item The \textbf{degree} of $(V, \varphi, \cF ,V_K)$ is $\displaystyle deg(V, \varphi, \cF V_K) = \sum_{i \in \Z} i \cdot dim_{K} (gr^i_\cF(V_K)) - v_F(det(\varphi))$, 
        \item The \textbf{slope} of $(V, \varphi, \cF V_K)$ is $\displaystyle\nu(V, \varphi, \cF V_K) = \frac{deg(V, \varphi, \cF V_K)}{rk(V, \varphi, \cF V_K)}$,
        \item If $(V, \varphi, \cF_K)$ is a filtered isocrystal, a \textbf{sub}-isocrystal is a triplet $(V', \varphi', \cF'_K)$ where $V'$ is a subspace of $V$ stable under $\varphi$, endowed with the induced structure, i.e. $\varphi' = \varphi_{V'}$, and $\cF'_K = \cF_K|_{V' \otimes_{K_0} K'}$.
    \end{itemize}

    We say that a filtered isocrystal is \textbf{weakly admissible} if it is semi-stable\footnote{Recall that an object is said to be semi-stable if every sub-object has a non-greater slope.} and of slope zero. A nontrivial result of Faltings and Totaro (cf. \cite{Faltings_TPofWA} and \cite{Totaro1996TensorPI}) states the following : 
    
    \begin{thm}
        \label{Totaro}
        The tensor product of semi-stable filtered isocrystals is semi-stable.
    \end{thm}

    \begin{rk}
    \label{HodgeNewton}
    This degree function (and, in turn, the slope function) is the difference of two natural and classical degree functions : one associated to the isocrystal, and the other to the filtration :
        \begin{enumerate}
        \item For $(V, \cF_K V)\in \operatorname{Fil}_k^K$ a filtered vector space, the degree $deg(V, \cF_K V) = \sum_{i} i\cdot  dim_K(gr^i_\cF V_K)$ is the y-axis value of the terminal point of the Hodge polygonal associated to the filtration $\cF V_K$.
        \item For an isocrystal $(V, \varphi)$ over $K_0$, we consider its degree $-v_p(det(A))$, where $A$ is a matrix representing the action of $\varphi$ in any given base. It is the y-axis value of the terminal point of the Newton polygonal associated with the isocrystal $(V, \varphi)$. 
    \end{enumerate}

    In other terms, a filtered isocrystal is weakly admissible if and only if its Newton and Hodge polygonals terminate at the same point, and the Newton polygonal of any sub-isocrystal lies above its Hodge polygonal. This is the original interpretation of Fontaine \cite[Defi 4.4.3]{Fontaine_Original_Art_Faiblement_Adm}.
    \end{rk}

    We may now define p-adic period domains for $\GL_n$ (cf.  \cite[Prop 1.36]{RappZink}).

\begin{prodef}
    Let $(V, \varphi) \in \operatorname{Isoc}(K_0)$, and $\nu$ be the possible type of a filtration on $V$.

    There exists an open subset defining a smooth, partially proper rigid-analytic variety over $K_0$ : $$\displaystyle \breve{\mathscr{F}}^{wa}(\GL(V), \varphi, \nu) \subset \sF(V, \nu)$$ such that, for any complete field extension $K$ of $K_0$, $$\breve{\sF}^{wa}(\GL(V), \varphi,\nu)(K) = \{ \cF \in \mathscr{F}(\GL(V),\nu)(K), (V \otimes_{K_0} K, \varphi \otimes_{K_0} id_{K}, \cF) \text{ is weakly admissible}\}.$$

    We call it the \textbf{period domain} associated to $(\GL(V), \varphi, \nu)$. 
    
    \end{prodef}    

Let us conclude this section by some examples of such period domains.

\subsubsection{Drinfeld's symmetric spaces} Fix $n \in \N_{\ge 1}$, $L$ a perfect field of characteristic $p$, $K_0 = \operatorname{Frac}(W(L))$, $V = K_0^n$, and $\sigma \in \operatorname{Aut}(K_0/\Q_p)$ be the Frobenius endomorphism. Let $\nu = (n-1, -1, \dots, -1)$.\\ For any field extension $K/K_0$, $\mathscr{F}(\GL(V), \nu)(K)$ parametrizes filtrations on $V \otimes_{K_0} K$ with only two jumps : a $1$-dimensional jump at $n-1$, and a $(n-1)$-dimensional jump at $(-1)$. 
    
    Any such filtration $\cF_{K} \in \mathscr{F}(\GL(V), \nu)(K)$ is of the following form, for some $K$-line $\cL_{\cF} \subset V \otimes_{K_0} K$: $$\cF_{K, x} = \begin{cases}
        V \otimes_{K_0} K & \text{ if } \hspace{1.7cm} x \le -1 \\
        \cL_{\cF} &\text{ if } \hspace{.52cm} -1 < x \le n-1 \\ 
        \{0\} &\text{ if } \hspace{.3cm} n-1 < x
    \end{cases} .$$
    Such filtrations hence correspond one-to-one to lines $\mathcal{L}_\cF$ of $V \otimes_{K_0} K$.

    Hence, $\mathscr{F}(\GL(V), \nu) = \P(V)$ is the $K_0$-variety whose $K$-points are $\P(V)(K) := \P(V \otimes_{K_0} K)$. 

  We may now define Drinfeld's symmetric spaces as open subsets of those flag varieties. Endow $V$ with the the trivial isocrystal structure $\sigma^{\times n}$.    One can check that a filtration $\cF$ is semi-stable if and only if the line $\cL_{\cF}$ is not contained in any $K_0$-rational hyperplane. Hence, $$ \breve{\sF}^{wa}\left(\GL(V), \varphi, \nu \right)(K) = \P^{n}(K) \setminus \bigcup_{H \in \mathcal{H}_{K_0} } H$$ where $\mathcal{H}$ describes all $K_0$-rational hyperplanes of $V$. Such a space will be denoted $\Omega^{n}_{K_0}$, and is called Drinfeld's symmetric space of dimension $n$ over $K_0$.

    It is naturally the base change to $K_0$ of the Drinfeld space over $\Q_p$, given by : $$\Omega^{n}_{\Q_p} = \P^n_{\Q_p} \setminus \bigcup_{H \in \mathcal{H}_{\Q_p} } H$$ where $\mathcal{H}$ describes all $K_0$-rational hyperplanes of $\Q_p^n$. 

\subsection{Period domains for arbitrary reductive groups}

We'll define period domains relative to more general reductive algebraic groups. In the case $G = \GL_n$, we will recover the above. The construction will be mostly the same, except that we replace :

\begin{itemize}[topsep=0pt]
    \item Isocrystals by isocrystals with G-structure.
    \item Filtrations by 1-parameter subgroups.
    \item Type of filtrations with conjugacy classes of 1-parameter subgroups.
\end{itemize}

\subsubsection{Gradings and filtrations} \label{GrFiltrations} Let $k$ be any field, and $K/k$ be a fixed field extension. Let $G$ be a connected algebraic group over $k$.

\begin{defi}
    We let $\operatorname{Rep}_{k}(G)$ be the category of finite dimensional representations of $G$ over $k$. 
    
    An object of $\operatorname{Rep}_{k}(G)$ is of the form $(V, \rho)$ where $V$ is a finite dimensional $k$-vector space, and $\rho$ is a morphism of affine group schemes $G \to \mathbb{GL}_V$, i.e. contains the compatible data of morphisms $\rho_R : G(R) \to GL_n(V \otimes_{k} R)$ for any $k$-algebra $R$. The category $\operatorname{Rep}_{k}(G)$  is naturally equipped with a tensor product, and the fiber functor obtained by forgetting the representation $w^G : (V, \rho) \in  \Rep_{k}(G) \mapsto V \in \operatorname{Vect}_{k}$.
\end{defi}

This makes it a Tannakian category over $k$.

Our goal is to understand filtrations \textit{of the fiber functor} $\omega^G$. Let us first introduce graduations. We let $\Q-\operatorname{Grad}_k$ denote the category of finite dimensional $k$-vector spaces, equipped with a grading indexed by $\Q$, as $V = \bigoplus_{x \in \Q} V_x$. Morphisms of such graded vector spaces are linear maps respecting the gradation. 

Consider the functor $\omega_0 : \Q-\operatorname{Grad}_k \to \operatorname{Vect}_k$ obtained by forgetting the grading : $\displaystyle \omega_0 \left(V = \bigoplus_{x \in \Q} V_x\right) \mapsto V$.

\begin{rk}
\label{grtofil}
This category of $\Q$-grading is related to the category of $\Q$-filtrations introduced in \ref{filtered}. Indeed, 
there exists a functor $ \operatorname{fil} : \Q-\operatorname{Grad}_k \to \Q-\operatorname{Fil}_k^K$, obtained by $$\operatorname{fil} : \displaystyle V = \bigoplus_{x \in \Q} V_x \mapsto \left(V_K, \cF^x V_K = \sum_{y \ge x} V_y \otimes K\right).$$

More generally, if we let $\Q-\operatorname{Grad}_{k}^K$ denote the category of vector spaces $V$ over $k$ together with a $\Q$-gradation on $V_K$, the above definition extends to a functor $\operatorname{fil} : \Q-\operatorname{Grad}_{k}^K \to \operatorname{Fil}_k^K$.

In the other direction, there exists a functor $\operatorname{gr} : \Q-\Fil_k^K \to \Q-\operatorname{Grad}_K$ obtained by $\displaystyle \cF \mapsto \bigoplus_{x \in \Q}gr^x_\cF \, \cF_K$.

Note that, while $\operatorname{fil}$ respects the underlying vector space (i.e. commutes with fiber functors), $\operatorname{gr}$ \textit{does not}, even if $K=k$.
\end{rk}

In order to study the category of $\Q$-gradings, we introduce the slope torus $\mathbb{D}_k =  \llim_{n \in \N} \bG_{m,k}$, as studied in \cite[IV.2]{dat_orlik_rapoport_2010}. We review its major properties in the following.

\begin{prodef}
    We define a functor $D_k : R \in  \operatorname{Ab}^{op} \mapsto \operatorname{Spec}(k[R]) \in \operatorname{AffSch}_k$. It is left adjoint to $\operatorname{Spec}(R) \mapsto R^*$. By construction, $D_k(\Z) = \bG_{m,k}$, so that : $$D_k(\Q) = D_k\left(\rlim_{n \in \N^*} \frac{1}{n} \cdot \Z\right) \cong \llim_{n \in \N^*} D_k(\Z) = \llim_{n \in \N^*} \bG_{m,k}.$$ 
    
    We let $\mathbb{D}_k = D_k(\Q)$ be the slope torus. The multiplication by any $x \in \Q$ induces a map $\chi_x : \D_k \to D_k(\Z) \cong \G_{m, k}$, viewed as a character of $\mathbb{D}_k$.
\end{prodef}

\begin{pro}
    \label{Tannak_Gradingsfromtorus}
    The category $\Q-\operatorname{Grad}_k$ is a neutral Tannakian category over $k$, and $\omega$ is a fiber functor over $k$. Moreover, there is a natural equivalence of category :
    
    \begin{center}
    \begin{tabular}{ccc}         $\operatorname{Rep}_k(\mathbb{D}_k)$ & $\to$ & $\Q-\operatorname{Grad}_k$ \\
         $(V, \rho)$ & $\mapsto$ & $\displaystyle V = \bigoplus_{x \in \Q} V_x$
    \end{tabular}
    \end{center}
    where $V_x = \{ v \in V, \rho(g)v = \chi_x(g) \cdot v, \forall g \in D(\Q)\}$ is the weight space of $V$ associated to $x \in \Q$.
\end{pro}

\begin{proof}
    cf. \cite{Rivano1972}, IV.1.2.1. 
\end{proof}

\subsubsection{Relative gradings and filtrations}
We keep the notations from \ref{GrFiltrations}. Let $G_K = G \times_{\Spec(k)} \Spec(K)$ denote the base change of $G$ to $K$. 

\begin{defi}
    \label{exa}
    A $\Q$-grading of $\omega^G$ over $K$ is a tensor functor $\operatorname{Gr} : \operatorname{Rep}_k(G) \to \Q-\Grad_K$ such that the diagram below commutes :
    \begin{center}
    \begin{tikzcd}
        & \Rep_k(G) \arrow[dr, "\operatorname{Gr}", blue] \arrow[rr, "\omega^G \otimes_k K"]& & \operatorname{Vect}_K \\
        & & \Q-\operatorname{Grad}_K \arrow[ur, "\omega_0"] & 
    \end{tikzcd}.        
    \end{center}
        Note that such functors are automatically exact and faithful.
\end{defi}

\begin{proof}
    The faithfulness is clear. The exactness follows from the exactness of $\omega^G$, hence of $\omega^G \otimes_k K$, and the fact that $\omega_0$ both preserves and reflects exact sequences. 
\end{proof}

\begin{lemma}
    \label{Gradingsfromtorus}
    There is a one-to-one correspondence between $\Q$-gradings of $\omega^G$ over $K$ and homomorphisms $\mathbb{D}_K \to G$. 
\end{lemma}

\begin{proof}
    This follows from Prop. \ref{Tannak_Gradingsfromtorus} and a form of Tannakian duality, cf. \cite[Corro 4.1.19]{dat_orlik_rapoport_2010}.
\end{proof}

Note that such a grading can be seen as the compatible data of a grading on $V \otimes_k K$ for each representation $(V, \rho)$ of $G$ (cf. \cite[Remark 4.2.5.(ii)]{dat_orlik_rapoport_2010} for a more precise formulation).

\begin{defi} \label{Filtr}
    A $\Q$-filtration of $\omega^G$ over $K$ is a tensor functor $\cF : \Rep_k(G) \to \Q-\operatorname{Fil}_k^K$ such that the diagram below commutes : 

    \begin{center}
    \begin{tikzcd}
        & \Rep_k(G) \arrow[dr, blue, "\cF"] \arrow[rr, "\omega^G"]& & \operatorname{Vect}_k \\
        & & \Q-\operatorname{Fik}_k^K \arrow[ur, "\omega_0"] & 
    \end{tikzcd}
    \end{center}
    
    and such that the composition $\operatorname{Rep}_k(G) \xrightarrow{\cF} \Q-\Fil_k^K \xrightarrow{gr} \Q-\Grad_K$ (where $gr$ is from \ref{grtofil}) is exact. 
\end{defi}

\begin{rk}
    \label{Filtr_induces_reps}
    Note that, by definition, such a filtration $\cF$ defined over $K$ induces, for all representation $(V, \rho)$ of $G$, a filtration $\cF^\rho$ on $V \otimes_k K_0$. One can reformulate the above definition as the \textit{compatible data} of filtrations $\cF^\rho$ on $V \otimes_k K$, cf. \cite[Defi 4.2.6]{dat_orlik_rapoport_2010}.
\end{rk}

The exactness condition ensures we can recover the filtration from the grading.

\begin{pro}
    \label{Gradingfromfil}
    Every $\Q$-grading of $\omega^G$ over $K$ naturally induces a $\Q$-filtration of $\omega^G$ over $K$. 
\end{pro}

\begin{proof}
    Let $\operatorname{Gr}$ be such a grading. The morphism $\operatorname{fil}$ from \ref{grtofil} commutes with fiber functors, so that the forgetful functor $\omega_0 : \Q-\operatorname{Grad}_{k,K} \to \operatorname{Vect}_K$ factors through $\Q-\operatorname{Fil}_K^K$. We have : 

    \begin{center}
    \begin{tikzcd}
        & \Rep_k(G) \arrow[ddrrr, start anchor = {south}, end anchor={west}, bend right = 30] \arrow[drr, "{V \mapsto (V, \operatorname{Gr}(V))}", near end, swap]\arrow[drrr, dashed] \arrow[rrrr, "\omega^G"]& & & & \operatorname{Vect}_k \\ 
        & & & \Q-\operatorname{Grad}_{k,K} \arrow[r, "fil", swap]& \Q-\Fil_k^K \arrow[ur, "\omega_0"] \arrow[d, "gr"] & \\
        & & & &\Q-\operatorname{\Grad}_K&
    \end{tikzcd}
    \end{center}

    The dashed arrow will be our filtration. 
    
    It suffices to show that the long curved composition $\operatorname{Rep}_k(G) \to \Q-\operatorname{Grad}_K$ is exact. The map $(V \mapsto (V, \operatorname{Gr}(V))$ is exact by the remark \ref{exa}, and $fil$ and $gr$ can be checked to be exact by an easy linear algebra argument.
\end{proof}

We conclude this paragraph by showing that gradings and filtrations relative to $\GL(V)$ correspond to "standard" grading and filtrations on $V$ (see also \cite[Remark 4.2.11]{dat_orlik_rapoport_2010}).

\begin{ex}
    \label{grfil_gln}
    Let $G = \GL(V)$. Then, $\Q$-gradings of $\omega^G$ over $K$ correspond to graduations on $V \otimes_k K$, and $\Q$-filtrations of $\omega^G$ over $K$ correspond to filtrations on $V \otimes_k K$.

    Indeed, any $\Q$-grading of $\omega^G$ over $K$ induces a $\Q$-grading on $V$ through the standard representation of $\GL(V)$ acting on $V$. Reciprocally, a $\Q$-grading $V = \bigoplus_{x \in \Q} V_x$ induces a morphism $D(\Q) \to \GL(V)$ as : $$g \mapsto \sum_{x \in \Q} i \circ \chi_x(g) \cdot id_{V_x}$$ where $i : \bG_m \to \GL(V)$ denotes the diagonal embedding, and $\chi_x : D(\Q) \to \bG_m$ is the map from \ref{Tannak_Gradingsfromtorus}.

    Likewise, filtrations on $\omega^G$ induce filtrations on any representations, and one can consider the standard representation on $V \otimes_k K$. Reciprocally, consider a filtration on $V \otimes_k K$. As a filtration on a finite dimensional vector space, it is necessarily split (i.e. lies in the image of the morphism $fil$ from \ref{grtofil}). Hence, by the discussion above, it defines a grading on $\omega^G$ over $K$. Then the filtration given by \ref{Gradingfromfil} works.
\end{ex}

\subsubsection{Cocharacters and flag varieties}
We keep the notations of \ref{GrFiltrations}. We fix $k^{sep}$ a separable closure of $k$ containing $K$.

\begin{defi}
    A $\Q$-{one-parameter subgroup} (or $\Q$-cocharacter) of $G$ defined over $K$ is a morphism of algebraic groups $\D_{K} \to G_K$. We let $X_*(G)_\Q^K$ denote the set of such cocharacters. 
\end{defi}

We call $X_*(G)_\Q^{k}$ the set of \textbf{rational} cocharacters, and $X_*(G)_\Q^{k^{sep}}$ the set of \textbf{geometric} cocharacters.

\smallskip
\begin{rk}
    The literature is quite inconsistent on whether a cocharacter denotes a rational cocharacter, or a geometric one.\footnote{This is partly due to the fact that most literature of algebraic groups was written in the language of algebraic varieties over algebraically closed fields.} In this paper, we always specify the field of definition.
\end{rk}
\smallskip

By \ref{Tannak_Gradingsfromtorus}, elements of $X_*(G)_\Q^{K}$ are equivalent to gradings of $\omega^G$ over $K$. In particular, they induce filtrations of $\omega^G$ over $K$, by \ref{grtofil}.

\begin{defi}
    
   We say that two 1-parameter subgroups $\lambda, \lambda' \in X_*(G)_\Q^K$ are \textbf{par-equivalent} if they define the same filtration on $\omega^G$ over $K$.
\end{defi}


For $\lambda \in X_*(G)_\Q^K$, we let $P_\lambda$ be the subgroup of $G_K$ defined, for any $K$-algebra $R$, by : 
\begin{equation*}
    \label{parab}
    P_{\lambda}(R) = \left\{ g \in G(R),\lim_{t \to 0} \lambda(t) \circ g \circ \lambda(t)^{-1} \text{ exists}\right\}.
\end{equation*}

\begin{pro}
    \label{parequivalence}
    Assume that $G$ is reductive. Then :
    \begin{enumerate}
        \item All filtrations of $\omega^G$ over $K$ are splittable, i.e. come from a grading by \ref{Gradingfromfil}.
        \item For any $\lambda$ $\in X_*(G)_\Q^K$, $P_{\lambda}$ is a parabolic subgroup of $G_K$.
        \item Let $\lambda, \lambda' \in X_*(G)_\Q^K$, and $\cF_\lambda$ (resp. $\cF_{\lambda'}$) be the filtration on $\omega^G$ associated to $\lambda$ (resp. $\lambda'$). The following are equivalent : 
        \begin{enumerate}
             \item $\cF_\lambda = \cF_{\lambda'}$, i.e $\lambda$ and $\lambda'$ are par-equivalent.
             \item There exists $g \in P_{\lambda}(K)$ such that $\lambda' = g \circ \lambda \circ g^{-1}$.
        \end{enumerate}
    \end{enumerate}
\end{pro}

\begin{proof}
    cf. \cite[Def 2.3/Prop 2.6]{Mumford_Geometric_Invariant_Theory} and \cite[Thm 4.2.13]{dat_orlik_rapoport_2010}.
\end{proof}

Finally, we turn our attention to conjugacy classes of such $\Q$-parameter subgroups, and the associated flag variety. There are two natural and compatible actions on $X_*(G)_\Q^K$ : 

\begin{itemize}[itemsep=.1cm]
    \item Assume $K/k$ is Galois. The Galois group $\Gal(K/k)$ acts, for $\sigma \in \Gal(K/k)$ and $\lambda \in X_*(G)^K_\Q$, by : $$^\sigma \lambda = \sigma_{G_K} \circ \lambda \circ \sigma_{\mathbb{D}_{K}}^{-1},$$ where $\sigma_G$ (resp $\sigma_{\mathbb{D}_{K}}$) is the natural action of $\Gamma_E$ on $G_K$ (resp. $\mathbb{D}_{K}$).
    \item The group $G(K)$ acts by conjugation, via : $$^g\lambda := g \cdot \lambda \cdot g^{-1}$$ where the right hand side uses the natural action of $G(K)$ on $G_K$.
\end{itemize}

There is a natural inclusion $X_*(G)^{k}_\Q \subset X_*(G)^{K}_\Q$ compatible with the action of $G$, and it follows from Galois descent (cf. \cite[Theorem 6.1.6.a and Paragraph 6.2.B]{BSL_NeronModels_descent}) that $(X_*(G)^{K}_\Q)^{\operatorname{Gal}(K/k)} \cong X_*(G)^{k}_\Q$. 

For simplicity (and following the literature), we simply denote $X_*(G)/G = X_*(G)_\Q^{k^{sep}}/G(k^{sep})$ the group of conjugacy classes of geometric cocharacters. It is a discrete $\operatorname{Gal}(k^{sep}/k)$-module.

We may now define our flag varieties of interest, that will classify filtrations of \textit{fixed type}. This will be replaced by a conjugacy class of geometric cocharacters.

Every conjugacy class $\{ \mu \} \in X_*(G)_\Q/G$, is defined over a finite extension of $k$, called its \textit{reflex field} :
$$E(G, \{ \mu\}) = (k^{sep})^{\Gamma_{\{\mu\}}} \text{ where } \Gamma_{\{\mu\}} = \operatorname{Stab}_{\operatorname{Gal}(k^{sep}/k)} 
(\{ \mu\}) = \{ g \in \operatorname{Gal}(k^{sep}/k), \forall \lambda \in \{ \mu\}, \phantom{}^{g} \lambda(t) \in \{ \mu \} \} $$ 

\begin{prodef}
    \label{FlagVar_Gen}
    Let $G$ be a connected reductive group over $k$, and $\{ \mu \}$ be a conjugacy class of geometric cocharacters of $G$, with reflex field $E$. Fix a separable closure $E^{sep}$ of $E$.

    There exists a unique projective reduced scheme of finite type over $E$, denoted $\mathscr{F}(G, \{ \mu\})$, whose $K$-points, for any intermediary field extension $E^{sep} / K / E$, are : 
    \begin{align*}
        \mathscr{F}(G, \{ \mu\})(K) &= \left\{\lambda \in X_*(G)_\Q^{K}, \lambda \in \{ \mu \}\right\}\big/ \text{ par-equivalence } \\ 
        &= \left\{ \cF \in \operatorname{Fil}_{K}(\omega^G), \exists \lambda \in \{ \mu \},  \cF = \cF_{\lambda} \right\}.
    \end{align*}

    This is called the \textit{flag variety} associated to $(G, \{ \mu\})$. 
\end{prodef}

\begin{proof}    
    The equivalence of the descriptions follows from the definition and first point of \ref{parequivalence}.

    Assume $G$ is quasi-split. Then, a lemma of Kottwitz (cf. \cite[Lemma 1.1.3]{Kottwitz_mu_lemma} (see also the proof of \cite[Prop 6.2]{Kottwitz_isoc_1})) shows that there exists an element $\mu \in \{ \mu\}$ defined over $E$.
    Let $P_{\mu, E}$ be the parabolic subgroup of $G_E$ associated to $\mu$ under \ref{parab}. Then $\mathscr{F}(G, \{ \mu\}) = G_E/P_{\mu, E}$ satisfies the above property, by \ref{parequivalence}. The unicity follows by descent from an algebraically closed field, where reduced schemes of finite types can be viewed as algebraic varieties.
    
    If $G$ is not quasi-split, we refer to the argument of \cite{dat_orlik_rapoport_2010}. 
\end{proof}

In the $\GL(V)$ setup, this coincides with the flag varieties defined in \ref{flag_gln}. We will then define period domains as analytic open subsets of the above flag varieties, corresponding to semi stable isocrystals.

\subsubsection{Isocrystals with additional structure} 

\label{isoc_add}
Let $G$ a connected reductive algebraic group over $\Q_p$, and $K_0 = \breve{\Q}_p$. Recall that the category $\operatorname{Isoc}(K_0)$ defined in \ref{isoc} is naturally endowed with the tensor product $(V_1, \varphi_1) \otimes (V_2, \varphi_2) := (V_1 \otimes V_2, \varphi_1 \otimes \varphi_2)$, which makes it a (non neutral !) $\Q_p$-linear\footnote{It is \textbf{not} $K_0$ linear, as endomorphisms of the unit object $(K_0, id_{K_0})$, are endomorphisms $K_0 \to K_0$ commuting with the Frobenius, so that they correspond with elements of $K_0^\sigma = \Q_p$.} Tannakian category, whose associated fiber functor is the forgetful functor $\omega_{\operatorname{Isoc}} : \operatorname{Isoc}(K_0) \to \operatorname{Vect}_{K_0}$. 

\begin{defi}
    An \textbf{isocrystal with G-structure} is a $\otimes$-exact functor $\Rep_{\Q_p}(G) \to \operatorname{Isoc}(K_0)$.
\end{defi}

It follows from the definition that such a functor commutes with fiber functors, since $G$ is connected (cf. \cite[Remark 9.1.6]{dat_orlik_rapoport_2010}). 

\begin{ex}
    \label{GLnRappRich}
    Let $V$ be a $\Q_p$-vector space of dimension $n$, and $G = \mathbb{GL}(V)$ be the group scheme defined by $R \mapsto GL(V \otimes_{\Q_p} R)$ $($so that there is a non canonical isomorphism $\GL(V) \cong \GL_n)$.

    Through the canonical representation of $\mathbb{GL}(V)$ acting on $V$, we can attach to any isocrystal with $\mathbb{GL}(V)$-structure a (standard) isocrystal over $V \otimes_{\Q_p} K_0$.
    
    It follows from the classification of irreducible representations of $\mathbb{GL}(V)$ via Schur functors that, reciprocally, the structure of an $\mathbb{GL}(V)$- isocrystal is fully determined by its value on the standard representation (cf. \cite[Remark 3.4.(ii)]{Rappoport_Richartz}) - so that $\mathbb{GL}(V)$-isocrystals correspond one-to-one with structure of isocrystals on the vector space $V \otimes_{\Q_p} K_0$.
\end{ex}

Any $b \in G(K_0)$ induces the isocrystal with $G$-structure $N_b : (V, \rho) \mapsto \left(V \otimes_{\Q_p} K_0, \rho_{K_0}(b)(id \otimes \varphi)\right)$.

\begin{pro}
    All isocrystals with $G$-structure are of the form $N_b$, for some $b \in G(K_0)$. \\ Moreover, there is an isomorphism of exact $\otimes$-functors $N_b \cong N_{b'}$ if and only if there exists $g \in G(K_0)$ such that $b' = gb\varphi(g)^{-1}$.
\end{pro}

\begin{proof}
    This relies on Tannakian duality, cf. \cite[Lemma 9.1.4 and Remark 9.1.6]{dat_orlik_rapoport_2010}.
\end{proof}

We define $B(G)$, the \textit{Kottwitz} set of $G$, as the quotient of $G(K_0)$ by the above equivalence relation. Hence, $B(G)$ classifies isocrystals with $G$-structure. For any $b \in G(K_0)$, we denote $[b]$ its associated class in $B(G)$. 

\subsubsection{Period domains} We reuse the notations from \ref{isoc_add}.
\label{Period_Domains}

\begin{defi}
    Let $G$ be an algebraic group over $\Q_p$, and $b \in G(K_0)$, with associated isocrystal with $G$-structure $N_b$. Recall that the datum of a point $\cF \in \sF(G, \{ \mu\})(K) \subset \operatorname{Fil}_K(\omega^G)$ defines a filtration $\cF^\rho$ on $V \otimes_k K$ for all representations $(V, \rho)$ of $G$ by \ref{Filtr_induces_reps}. 
    
    We say that $(b, x)$ is \textbf{weakly admissible} if, for all representation $(V, \rho)$ of $G$, the filtered isocrystal $(N_b(V), \cF^x V \otimes_k K)$ is weakly admissible in the sense of \ref{FilIC}.
\end{defi}

By the tensor product theorem \ref{Totaro}, this is equivalent to asking that   $(N_b(V), \cF^x V \otimes_k K)$ is weakly admissible for \textit{some} faithful representation of $G$ (see also \cite[Defi 1.18]{RappZink}).





We may finally define general period domains. Fix the datum $(G, [b], \{ \mu\})$ of :
    \begin{enumerate}[topsep = 0pt]
        \item A connected quasi-split group $G$ over $\Q_p$. 
        \item An element $[b]$ of the Kottwitz set $B(G)$, i.e. a $\sigma$-conjugacy class in $G(K_0)$.
        \item A conjugacy class $\{ \mu \} \in X_*(G)/G$ of geometric cocharacters of $G$.
    \end{enumerate}
    
\begin{defi}
    \label{def_PD}
    Let $(G, [b], \{ \mu\})$ be a datum as above. Let $E = E(G, \{ \mu\})$ be the associated reflex field, and $\breve{E} = E \cdot K_0$. There exists a unique partially proper open subset : $$\breve{\sF}^{wa}(G, b, \{ \mu\}) \subset \sF(G, \{ \mu\}) \otimes_E \breve{E}$$ such that, for any field extension $K/\breve{E}$, 
    $$\breve{\sF}^{wa}(G, b, \{ \mu\})(K) = \left\{\cF \in \sF(G, \{ \mu\})(K), \text{ the filtered isocrystal }(N_b, \cF) \text{ is weakly admissible}\right\}.$$

    This is called the \textbf{period domain} associated to $(G, b, \{ \mu \})$. Moreover, it can be obtained as the analytification of a Berkovich space.
\end{defi}

\begin{proof}
    cf. \cite[Prop 1.36]{RappZink} in the rigid analytic setup and \cite[Prop 9.5.3]{dat_orlik_rapoport_2010} in the Berkovich setup.
\end{proof}

Note that, up to an isomorphism, the period domain only depends on the class $[b] \in B(G)$. Indeed, if $b' = gb\sigma(g)^{-1}$ for some $g \in G(K_0)$, the map $\lambda \mapsto ^g \hspace{-.10cm}\lambda$ induces an isomorphism $\breve{\sF}^{wa}(G, b, \{\mu\}) \cong \breve{\sF}^{wa}(G, b', \{\mu\})$.

Let $J_b$ denote the automorphims group of $N_b$, defined, for any $\Q_p$-algebra $R$, by : $$ J_b(R) = \{ g \in G(\breve{\Q}_p \otimes_{\Q_p} R), gb\sigma(g)^{-1} = b\}.$$ The group $J_b(\Q_p)$ acts naturally on $\breve{\sF}^{wa}(G, b, \{\mu\})$, via $J_b(\Q_p) \subset G(K_0)$ via conjugation $\lambda \mapsto ^g \hspace{-.10cm}\lambda$.

\subsubsection{Decency and basicness} We reuse the notations from \ref{isoc_add}. For any $b \in G(\Breve{\Q}_p)$, Kottwitz \cite[sec 4.2]{Kottwitz_isoc_1}) defines a \textit{slope} morphism $\nu_b : \mathbb D_{K_0} \to G_{K_0}$ as follows : 

\begin{defi}
    \label{Slope_defi}
    For any $V \in \Rep_{\Q_p}(G)$, the $\Q$-grading on the vector space $V \otimes_{\Q_p} K_0$ associated with the representation of the pro-torus $\mathbb{D}_{K_0} \xrightarrow{\nu_b} G_{K_0} \xrightarrow{\rho_{K_0}} \operatorname{GL}(V \otimes K_0)$ coincides with the Dieudonné-Manin slope decomposition of the isocrystal $N_b$, as recalled in \ref{DN}.
\end{defi}

Note that, if $b' = gb\sigma(g)^{-1}$, then $ \nu_{b'} = \operatorname{Int}(g) \circ \nu_b$. 

\begin{defi}
    We say that an element $b \in G(K_0)$ is \textbf{s-decent} for an integer $s \ge 1$ if :
    \begin{itemize}
        \item The morphism $s \cdot \nu_b : \mathbb{D}_{K_0} \to G_{K_0}$ factors through $\mathbb{G}_{m,{K_0}}$,
        \item The following equality holds in $G(K_0) \rtimes \sigma^\Z$ : $(b\sigma)^s = (s \nu_b)(p) \sigma^s$.
    \end{itemize}
\end{defi}

Since $G$ is connected, every class $[b] \in B(G)$ contains an $s$-decent element, for some $s \ge 1$ (cf. \cite[lemma 9.1.33, remark 9.1.34]{dat_orlik_rapoport_2010}). By \cite[lemma 9.6.19]{dat_orlik_rapoport_2010}, if $G$ is furthermore assumed to be quasi-split, each $[b] \in B(G)$ contains an $s$-decent element $b$ such that $\nu_b$ is defined over $\Q_p$.

\begin{pro}
    \label{Field_def_E_s}
    Let $(G, [b], \{ \mu \})$ be a datum as above. Let $s \ge 1$, and assume that $b$ is $s$-decent, with reflex field $E$. Then, the period domain $\breve{\sF}^{wa}(G, b, \{ \mu\})$ admits a model over $E_s := E \cdot \Q_{p^s}$.
\end{pro}

\begin{proof}
    cf. \cite[Prop 1.36.(ii)]{dat_orlik_rapoport_2010}.   
\end{proof}

\begin{defi}
    We say that an element $b \in B(K_0)$ is \textbf{basic} if the slope morphism $\nu_b$ factors through $Z_{K_0}$, where $Z$ denotes the center of $G$. We let $B(G)_{basic}$ denote the subset of $B(G)$ consisting of basic elements.
\end{defi}

Note that an element $b$ is basic if and only if $J_b$ is an inner form of $G$.





\section{Primitive comparison theorem with compact support for period domains}
\label{Sec7}

The goal of this section is to show that period domains, as defined in the previous paragraph, satisfy primitive comparison with compact support (the main result is \ref{PCCS_Period_Domains}). This computation is based on standard sheaf techniques, together with a non-trivial geometric stratification of the complement $\sF \setminus \sF^{wa}$ of period domains, constructed by S.Orlik in \cite{Orlik_local}. 

In loc.cit., Orlik uses his geometric decomposition to compute the étale cohomology with compact support of period domains with $\Z/\ell^n\Z$ coefficients. In \cite{cdhn}, the authors adapt the computation for $\ell = p$. While the geometric part of the argument does not really depend on the coefficients, working with $p$-adic coefficients adds a lot of representation-theoretic complexity, as one needs to compute some $\operatorname{Ext}^1$ groups between Steinberg representations.

In this section, we compute the $\O^+/p$ cohomology with compact support of period domains. Going from $\F_p$ to $\O^+/p$-coefficients adds no representation theoretic complexity, as only the "geometric" part of the argument is needed.

In this article, we black-box the origin of the geometric decomposition (which uses ideas from geometric invariant theory, cf. \cite{Mumford_Geometric_Invariant_Theory}), and directly study the consequences for period domains, as they will be reviewed in Prop \ref{Shubert}.

Compared to the previous expositions of the geometric decomposition and its consequence on cohomology appearing in the literature, we provide the following technical improvements : 

\begin{enumerate}
    \item Throughout most of the computation, we work with coefficients in a general abelian étale sheaf. Doing so, we remove a superfluous overconvergence hypothesis.
    \item We work directly over the field of definition of the period domains, rather than over a complete algebraic closure. This allows to properly keep track of the Galois action.
    \item We work with the language of locally spatial diamonds throughout, rather than pseudo-adic spaces.
\end{enumerate}



Let us now precisely state our result :

We fix a datum $(G, [b], \{ \mu \})$ as in \ref{def_PD}, where we assume furthermore $b$ to be basic and s-decent for some integer $s \ge 1$, and we use notations of section \ref{Sec6}. 

For simplicity, we simply denote $\sF = \sF(G, b, \{ \mu\})$ and $\sF^{wa} = \sF^{wa}(G, b, \{ \mu\})$, with field of definition $E_s$ (as in \ref{Field_def_E_s}), and $\sF_{E_s}$ be the base change of $\sF$ to $E_s$, and $Y = \sF_{E_s} \setminus \sF^{wa}$, viewed as a locally spectral diamond over $E_s$ by \ref{ComplementLSD}. Let $J := J_b$ be the inner form of $G$ associated to the basic element $b$.

We also let $\O^+_{E_s}/p := \O^+_{\operatorname{Spa}(E_s, \O_{E_s})}/p$.

We fix $\overline{E_s}$ an algebraic closure of $E_s$, $C$ be a completion of $\overline{E_s}$, and $\sF_C$, $\sF^{wa}_C$ and $Y_C$ be the respective base change to $C$.

In the rest of this section, every pullback and pushforward is in the sense of Scholze (and coincides with the one of Huber between analytic adic spaces).

\begin{thm}
    \label{PCCS_Period_Domains}
    In the above setup, let $j : \sF^{wa} \to \sF_{E_s}$ be the open immersion, and $j_C :  \sF^{wa}_C \to \sF_{C}$ be its base change. Let $\cL$ be an étale $\F_p$-local system on $\sF_{E_s}$, and $\cL_C$ be its pullback to $\sF_{C}$. 
    
    Then, $\sF^{wa}_C$ satisfies Poincaré duality with respect to $j^*_C \cL$, i.e, for any $0 \le k \le 2d$, there is a $J(\Q_p) \times \operatorname{Gal}(\overline{E_s}/E_s)$-equivariant isomorphism : $$ H_{\e t}^{k}(\sF_C, j_C^* \cL_C) \cong \operatorname{Hom}_{\F_p}(H^{2d-k}_{\e t, c}(\sF_C, j_C^* \cL_C), \F_p)$$
    
    where $d$ denotes the dimension of $\sF_C$.
\end{thm}


Using the criterion \ref{Test_Comple}, the statement reduces to a establishing primitive comparison on $\sF_{C}$ and on $Y_C$. As $\sF_C$ is smooth and proper, the primitive comparison on $\sF_C$ follows directly from the primitive comparison theorem \ref{PrimComparison_Galois}, so that it suffices to establish the following (applied to the pullback of a local system on $\sF_{E_s}$, and base changing to $C$ using \cite[Prop 16.10]{Scholze_ECohD}) :

\begin{pro}
    \label{complement}
    Let $h : Y \to \operatorname{Spa}(E_s, \O_{E_s})$ denote the structure morphism of $Y$. Let $\cL$ be an étale $\F_p$- local system on $Y$. Then, for all $q \in \N$, the natural map :
    $$R^q h_* \cL \otimes \O^+_{E_s}/p \to R^q  h_* (\cL  \otimes_{\F_p} \O_{Y}^+/p) $$ is an almost isomorphism in the category of étale  $\O^+_{E_s}/p$-sheaves over $\operatorname{Spa}(E_s, \O_{E_s})$.
\end{pro}

Note that this morphism (as well as the primitive comparison morphism on $\sF_C$) is naturally $J(\Q_p)$-invariant, so that the $J(\Q_p)$-invariance in \ref{PCCS_Period_Domains} follows from \ref{Functoriality} by mimicking the proof of \ref{PCCS_PD_Gal}.

\subsection{Orlik's geometric decomposition} Let us fix an étale $\F_p$-local system $\cL$ on $Y$. In this section, we present a geometric decomposition of $Y$, due to Orlik in \cite{Orlik_local} (see also \cite{Orlik_Finite} for the simpler analogue over finite fields).

Fix $S$ a maximal $\Q_p$-split torus of $J$, and let $\Delta$ be a basis of the relative root system of $G$, which is a finite set. The following proposition regroups the properties that we will need :
    
\begin{pro}
    \label{Shubert}
    Let $P_0$ be a fixed minimal parabolic subgroup of $J$ containing $S$. For any $I \subset \Delta$, let $P_I$ be the associated standard parabolic subgroup of $J$.
    
    There exists a finite set of rigid-analytic subvarieties $Y_I \subset Y$ indexed by subsets $I \subset \Delta$ together with a stratification :  $$|Y| = \bigcup_{|\Delta \setminus I| = 1} \bigcup_{g \in J(\Q_p)} g \cdot |Y_I|$$

    such that: 
    \begin{enumerate}[topsep = 0pt]
        \item \label{firstpoint} For any $I \subset \Delta$, $Y_I$ is the analytification of a proper algebraic variety over $E_s$ (that is moreover a Schubert variety).
        \item For any $I, J \subset \Delta$, $Y_I \cap Y_J = Y_{I \cap J}$.
        \item \label{proper} $(J/P_I)(\Q_p) = J(\Q_p)/P_I(\Q_p)$.
        \item \label{Lastpoint} The action of $J(\Q_p)$ on $\sF_{E_s}$ restricts to an action of $P_I(\Q_p)$ on $Y_I$.
    \end{enumerate} 
\end{pro}

\begin{proof}
    Cf. \cite[Section 3]{Orlik_local} or \cite[Section 5.4]{cdhn}.
\end{proof}

From the above stratification, we will construct a functorial resolution $\mathcal{C}(\cF)$ of any sheaf of abelian group $\cF$ on $Y_{\e t}$, relating the cohomology of $\cF$ on $Y$ to the cohomology of $\cF$ restricted to the $Y_I$, for varying $I \subset \Delta$. This exactness was only considered in \cite{Orlik_local} for \textit{overconvergent} étale sheaves $\cF$, but, as we will see, this assumption is not needed, so that it also works for $\cF = \cL \otimes \O_Y^+/p$.

The strategy to prove primitive comparison with compact support for $Y$ is as follows :

\begin{enumerate}
    \item From the functorial resolution $\cC(\cF)$, construct a spectral sequence, whose first page computes the cohomology of $\cF$ on the $Y_I$, and which converges to the cohomology of $\cF$ on the $Y$ - in a way that is functorial in the sheaf $\cF$. 
    \item Use the primitive comparison theorem to compare the first page of the spectral sequence with coefficients in $\cL$ or in $\cL \otimes \O_Y^+/p$\footnote{Here, the primitive comparison theorem will suffice, but one can explicitly compute the $\F_p$-cohomology groups of the $Y_I$ (at least after base change to an algebraically closed base) ; cf. \cite[Corrolary 5.11]{cdhn}}.
    \item Prove that this comparison goes through the convergence of the spectral sequence in order to deduce a primitive comparison result for the cohomology of $Y$.  
\end{enumerate}

For any $I \subset \Delta$, we let $X_I = J(\Q_p)/P_I(\Q_p)$. By the third point of \ref{Shubert}, this is a profinite set - as it identifies to the $\Q_p$-points of a proper algebraic variety (which can then be seen as $\Z_p$-points using the valuation criterion). We will need the following lemma.


\begin{lemma}
    \label{generalizing_union}
    For any compact open subset $T \subset J/P_I(\Q_p)$, the space $Z_I^T = \bigcup_{g \in T} g \cdot Y_I$ admits a natural structure of a closed locally spatial sub-diamond of $Y$.
\end{lemma}

\begin{proof}
    Note, by that \ref{Shubert}, $g \cdot Y_I$ is well defined for $g \in J/P_I(\Q_p) = J(\Q_p)/P_I(\Q_p)$ as the action of $P_I(\Q_p)$ stabilizes $Y_I$. 
    
    Fix a compact open $T \subset X_I$. By \ref{ComplementLSD}, it suffices to show that the underlying topological space of $Z_I^T$ is closed and generalizing. It's closed by the proof of \cite[Lemma 3.2]{Orlik_local}. All of the $g \cdot |Y_I|$ are generalizing since they are Zariski closed subsets, and a union of generalizing spaces remains generalizing.
\end{proof}

For any compact open subset $T \subset X_I$, we let $Z_{I}^T$ be as above and $i_T : Z_I^T \to Y$ be the associated closed immersion. 

\begin{defi}
    \label{Disj}
    We let $\operatorname{Disj}_{X_I}$ be the category of covers of $X_I$ by disjoint nonempty open subsets, with morphisms given by refinements. Such disjoint covers are necessarily finite, as $X_I$ is compact.
    
    Objects of $\operatorname{Disj}_{X_I}$ will be denoted $c = \{T_k\}_{1 \le k \le n(c)}$.

    For any $c \in \operatorname{Disj}_{X_I}$ and any étale sheaf $\cF$ on $Y$, we let $\cF_c := \displaystyle \bigoplus_{k=1}^{n(c)} i_{T_k *} i_{T_k}^* \cF$, and $\cF_I := \rlim_{c \in \operatorname{Disj}_{X_I}} \cF_c$.
\end{defi}

Note that, if $\cF$ is $J(\Q_p)$-equivariant, then $\cF_I$ is also $J(\Q_p)$-equivariant.

For a fixed $c$, $\cF_c$ admits a natural interpretation as the étale sheaf on $X_I$ whose sections are locally constant sections of $\cF$, trivialized on the disjoint covering $c$. Likewise, $\cF_I$ can be thought of as the sheaf of \textit{locally constant} sections of $\prod_{x \in X_I} \cF_{x Y_{I}^{ad}}$. This is made more explicit in \cite[p.538]{Orlik_local}.

\smallskip
\begin{rk}
\label{Morph_Disj}
Let us briefly explicit the transition morphisms appearing in the colimit. 
Let $c = \{ T_k\}_{1 \le k \le n(c)}$ be an object of $\operatorname{Disj}_{X_I}$, and $c' = \{ S_{k,l}\}_{1 \le k \le n(c) ; 1 \le l \le m(k)}$ be a refinement of $c$, in the sense that, for any $1 \le k \le n(c)$, the $(S_{k,l})_{1 \le l \le m(k)}$ form a disjoint compact open covering of $T_k$.\footnote{Note that, since coverings are disjoint, if such a refinement exists, the refinement map is uniquely determined.}

The morphism $\cF_{c'} \to \cF_c$ is then given as a direct sum of the natural morphisms : $i_{S_{k,l} *} i_{S_{k,l}}^* \cF \to i_{T_{k} *} i_{T_{k}}^* \cF$ induced by the unit of the adjunction $\cF \to i_{S *} i_S^* \cF$ together with the natural isomorphism : $$ i_{S_{k,l} *} i_{S_{k,l}}^*\cF \cong i_{S_{k,l} *} i_{S_{k,l}}^* i_{T_k *} i_{T_k}^* \cF$$ 

induced by the unit of the adjunction $\cF \to i_{T_k *} i_{T_k}^* \cF$. The fact that it is indeed an isomorphism can be checked on stalks using \ref{stalkpf}.
\end{rk}

Let us study the sheaves $\cF_I$. We start with the following lemma (which is a version of \cite[Prop 6.11.(a)]{cdhn} over an arbitrary base field, and in the language of diamonds) :

\begin{lemma}
    \label{Stalk}
    Let $\bar{x}$ be a geometric point of $Y_{\e t}$, with underlying topological point $x$, and $\cF$ be an étale sheaf of abelian groups over $Y$. The stalk of $\cF_I$ at $\bar{x}$ identifies with : $$(\cF_I)_{\bar{x}} \cong LC(X_I(x), \cF_{\bar{x}})$$
    
    where $X_I(x) = \{g \in X_I, x \in g \cdot Y_{I}\}$ is a closed subset of $X_I$, and $LC$ denotes the set of locally constant functions on $X_I(x)$ (with the topology induced by $X_I$), with value in the stalk $\cF_{\overline{x}}$.
\end{lemma}

\begin{proof}
    Taking stalks commutes with colimits, so that : $$(\cF_I)_{\bar{x}} \cong \rlim_{c \in \operatorname{Disj}_{X_I}} \bigoplus_{k=1}^{n(c)} (i_{T_k, *} i^{*}_{T_k} \cF)_{\bar{x}}$$

    Using \ref{stalkpf}, since $i_{T_k}$ is a generalizing closed immersion, it follows that, for all $1 \le k \le n(c)$ : $$(i_{T_k*} i_{T_k}^* \cF)_{\bar{x}} = \begin{cases}
        \cF_{\bar{x}} & \text{if }\,  x \in T_k \cdot Y_I, \text{ i.e. if } T_k \cap X_I(x) \neq \emptyset  \\
        0 & \text{otherwise}
    \end{cases}$$ 
    We may then replace all terms in the direct sum by $i_{T_k \cap X_I(x), *} i_{T_k \cap X_I(x)}^* \cF$.
    
    If $c = \{ T_k\}_{1 \le k \le n(c)} \in \operatorname{Disj}_{X_I}$, the $T_k \cap X_I(x)$ form a disjoint compact open cover of $X_I(x)$. Reciprocally, any $c' \in \operatorname{Disj}_{X_I(x)}$ can be refined by a covering of the above form, by the lemma \ref{Tech_Top_Lemma} below. 
    
    Hence, coverings of the form $\{T_k \cap X_I(x)\}_{1 \le k \le n(c)}$ for $\{T_k\}_{1 \le k \le n(c)} \in \operatorname{Disj}_{X_I}$ form a cofinal system in $\operatorname{Disj}_{X_I(x)}$, so that :
    $$(\cF_I)_{\bar{x}} = \rlim_{c \in \operatorname{Disj}_{X_I}} \bigoplus_{k=1}^{n(c)} \left(i_{T_k \cap X_I(x) *} i_{T_k \cap X_I(x)}^* \cF\right)_{\bar{x}} \cong \rlim_{c' \in \operatorname{Disj}_{X_I(x)}} \bigoplus _{k = 1}^{n(c')} \cF_{\overline{x}} \cong LC(X_I(x), \cF_{\bar{x}}).$$

    The last isomorphism follows from identifying a locally constant function to the set of its value over any disjoint trivializing cover (that is necessarily finite since $X_I(x)$ is compact).
\end{proof}

The following lemma was used in the proof : 

\begin{lemma}
    \label{Tech_Top_Lemma}
    Let $X$ be a profinite topological space, and $F \subset X$ be a closed subset. \\ Let $c = \{ T_k\}_{1 \le k \le n(c')} \in \operatorname{Disj}_F$. Then, there exists $\{ S_{k,l} \}_{1 \le k \le n(c), 1 \le l \le m(k)} \in \operatorname{Disj}_X$ such that, for all $k$, the $\{S_{k,l} \cap T\}_{1 \le l \le m(k)}$ form a disjoint clopen union of $T_k$.
\end{lemma}

\begin{proof}
    For any $1 \le k \le n(c)$, $T_k$ is an open subset of $F$, so it can be written as $T_k = U_k \cap F$ for some open subset $U_k$ of $X$. Consider the open cover $X = (X \setminus F) \cup \bigcup_{{1 \le k \le n(c)}} U_k$ of $X$. By \cite[\href{https://stacks.math.columbia.edu/tag/08ZZ}{Lemma 08ZZ}]{stacks-project}, we may find a refinement by a disjoint cover formed by clopen subsets of $X$ : $$ X = S'_0 \sqcup \bigsqcup_{{1 \le k \le n(c), 1 \le l \le m(k)}} S'_{k,l}$$ where $S'_0 \subset X \setminus F$ and, for all $k,l$ in range, $S'_{k,l} \subset U_k$.

    We then let $S_{k,l} := S'_{k,l} \cap T$, which are closed and open subsets of $T$. Pick some arbitrary $(k_0, l_0)$ in range. After replacing $S_{k_0,l_0}$ by $S_{k_0,l_0} \sqcup S'_0$, the family $\{ S_{k,l}\}$ satisfies the desired properties.
\end{proof}

\subsection{Orlik's fundamental complex}

We may now define the fundamental complex, promised in \ref{Shubert}.

\begin{pro}
    \label{fundamental}
    Let $\cF$ be a sheaf of abelian groups on $Y_{\e t}$. The following is complex of étale sheaves on $Y_{\e t}$ is exact : $$ \mathcal{C}(\mathcal{F}) :  0 \to \cF \to  \bigoplus_{|\Delta \setminus I| = 1} \cF_I \to \bigoplus_{|\Delta \setminus I| = 2} \cF_I \to \dots \to  \bigoplus_{|\Delta \setminus I| = |\Delta| - 1} \cF_I \to \cF_\emptyset \to 0$$

    where the morphisms are defined as follows : fix $I \subset I' \subset \Delta$ with $|I' \setminus I| = 1$, The inclusion $P_I \subset P_{I'}$ induces a quotient map $p_{I, I'} : X_I \to X_{I'}$. For any $x \in X_{I'}$, we consider a map $\cF_{x \cdot Y_{I'}} \to \cF_{p_{I, I'}(x) \cdot Y_I}$ induced by the natural inclusion $x \cdot Y_I \subset p_{I, I'}(x) \cdot Y_{I'}$. The collection of such maps defines a map $f_{I', I} : \cF_{I'} \to \cF_I$. 
    
    Finally, recall that the root system $\Delta$ is naturally ordered, so that we may write $I' = \{\alpha_1 < \dots < \alpha_n \}$, and $I = I' \setminus \{ \alpha_i\}$. We then let $d_{I, I'} = (-1)^i p_{I, I'} : \cF_{I'} \to \cF_I$, and consider the direct sum along $I, I'$.
\end{pro}

When $\cF$ is a sheaf of $\Z/n\Z$-modules for some $n$ prime to $p$, the exactness is due to \cite[Thm 3.3]{Orlik_local}. This has been generalized in \cite[Thm 2.1]{Orlik_Cont} (see also \cite[Thm 6.13]{cdhn}) for overconvergent étale sheaves of abelian groups. 

\begin{proof}
    The étale topos of $Y$ has enough points by \cite[Prop 14.3]{Scholze_ECohD}, so that it we may check exactness on stalks. Let $\eta : \overline{x} \to Y$ be a geometric point, and $x \in Y$ be the image of the unique closed point.

    By \ref{Stalk}, the stalk at $\overline{x}$ is : 
    $$ C(\cF)_{\overline{x}} : 0 \to \cF_{\overline{x}} \to \bigoplus_{|\Delta \setminus I| = 1} LC(X_I(x), \cF_{\overline{x}}) \to \dots \to \bigoplus_{|\Delta \setminus I| = |\Delta| - 1} LC(X_I(x), \cF_{\overline{x}}) \to LC(X_\emptyset(x), F_{\overline{x}}) \to 0.$$

    However, by \cite[Thm 5.4]{Scholze_Condensed} (result attributed to Bergman), the modules of locally constant maps on profinite spaces are free, so that $C(\cF)_{\overline{x}} \cong C(\Z)_{\overline{x}} \otimes_{\Z} \cF_{\overline{x}}$. Since free $\Z$-modules are flat, it suffices to show that $C(\Z)_{\overline{x}}$ is exact, which follows from \cite[Thm 2.1]{Orlik_Cont}.
\end{proof}

\begin{rk}
\label{Fundamental_TensorProd}
Note that the above proof also shows that, for any two sheaves $\cF$, $\cG$ of abelian groups on $Y_{\e t}$, we have $\cC(\cF \otimes_\Z \cG) \cong \cC(\cF) \otimes_\Z \cG$. 
\end{rk}

\begin{rk}
    The proof of \cite[Thm 3.3]{Orlik_local} can in fact be adapted to work directly in the above setup. Indeed, it is shown that the above localized complex $C(\cF)_{\overline{x}}$ computes the cohomology of the constant sheaf\footnote{Note that, in that setup, constant sheaf cohomology needs \textit{not} coincide with singular cohomology.} $\cF_{\overline{x}}$ on a topological space obtained as the geometric realization of some simplicial subcomplex of the topologized combinatorial Tits building of $J$, which can be shown to be contractible - and hence has trivial cohomology for arbitrary constant sheaves.
\end{rk}

We may now consider the spectral sequence obtained by applying $R h_*$ to $\cC(\cF)$, (recall that $h : Y \to \operatorname{Spa}(E_s, \O_{E_s})$ is the structure morphism of $Y$). Considering the filtration on $\cC(\cF)$ (minus the $\cF$ term) obtained by truncations, we deduce a spectral sequence : $$E_1^{p,q}(\cF) := R^qh_*\left(\bigoplus_{|\Delta \setminus I| = p+1} \cF_I\right) \implies R^{p+q} h_* \cF =: E_\infty^{p+q}(\cF).$$

The above spectral sequence is clearly functorial in $\cF$. We will show the following :

\begin{pro}
    \label{PCTermSS}    
    Let $\cL$ be a local system of $\F_p$-vector spaces on $Y$. Then, for all $p, q \in \N$, the natural morphism : $$R^qh_* \left(\bigoplus_{|\Delta \setminus I| = p+1}  \cL_I\right) \otimes_{\F_p} \O^+_{E_s}/p\to R^q h_* \left(\bigoplus_{|\Delta \setminus I| = p+1} \left( \cL \otimes_{\F_p} \O_{Y}^+/p \right)_I\right) $$ is an almost isomorphism.
\end{pro}

Before proving this proposition, let us show that it implies \ref{complement}.

\begin{pro}
    Assume that the result of \ref{PCTermSS} holds. Then, for all $q \ge 0$, the natural morphism : $$(R^q h_* \cL) \otimes_{\F_p} \O^+_{E_s}/p \to R^q h_* (\cL \otimes \O^+_{Y}/p)$$

    is an almost isomorphism, so that \ref{complement} holds (and, in turn, so does \ref{PCCS_Period_Domains}).
\end{pro}

\begin{proof}
    There is a natural morphism of complexes of sheaves on $Y_{\e t}$ : $$\cC(\cL) \otimes_{\F_p} h^* \O^+_{E_s}/p \to \cC(\cL) \otimes_{\F_p} \O_{Y}^+/p \cong \cC(\cL \otimes_{\F_p} \O_{Y}^+/p)$$

    using the isomorphism from \ref{Fundamental_TensorProd}. This induces a morphism of complexes on $\operatorname{Spa}(E_s, \O_{E_s})_{\e t}$ : $$R^q h_* \cC(\cL) \otimes_{\F_p} \O^+_{E_s}/p \to R^q h_* \left(\mathcal{C}(\cL) \otimes_{\F_p} h^* \O^+_{E_s}/p \right)\to R^q h_* \mathcal{C}(\cL \otimes_{\F_p} \O^+_Y/p)$$
    
    where the first map is the one of \ref{Morph_Proj}, and the second one is induced by the map $h^* \O^+_{E_s}/p \to \O^+_Y/p$ using \ref{Fundamental_TensorProd}. This induces a morphism of spectral sequences $E_k^{p,q} (\cL) \otimes_{\F_p} \O^+_{E_s}/p \to E_k^{p,q} (\cL \otimes_{\F_p} \O_Y^+/p)$.
    
    All morphisms above commute with the localization to the category of sheaves of almost $\O^+_{E_s}/p$-modules, which is a an abelian category. By \ref{PCTermSS}, the considered morphism of spectral sequences induces an almost isomorphism on $E_1$, so that it also induces an almost isomorphism on $E_{\infty}$. This is the desired result by construction of the spectral sequence.
\end{proof}


\subsection{Proof of Prop \ref{PCTermSS}}

Let us rewrite the terms appearing in the proposition. Fix an étale $\F_p$-local system $\cL$ on $Y$, and let $\cF$ be an arbitrary étale sheaf of abelian groups on $Y$.

Since $Y$ is proper, it necessarily qcqs, so that the morphism $Y \to \operatorname{Spa}(E_s, \O_{E_s})$ is a coherent morphism of algebraic topoi (cf. \cite[Section 8]{Scholze_ECohD}). Then, by \cite[Thm 5.1]{SGA4}, $R^q h_*$ commutes with filtered inductive limits.

    Direct sums and the colimit along $\operatorname{Disj}_{X_I}$ are filtered, so that, for all $p,q \in \N$ : $$R^q h_*\left( \bigoplus_{|\Delta \setminus I| = p+1} \cF_I\right) \cong \bigoplus_{|\Delta \setminus I| = p+1} R^q h_* \cF_I \cong \bigoplus_{|\Delta \setminus I| = p+1} \rlim_{c \in \operatorname{Disj}_{X_I}} \bigoplus_{i=1}^{n(c)} R^q h_* \left(i_{T_k *} i_{T_k}^* \cF\right).$$

    The map $i_{T_k}$ is a closed immersion, hence, by \cite[Lemma 21.13]{Scholze_ECohD}, $R^q h_* \circ i_{T_k *} = R^q (f \circ i_{T_k})_*$.
    
    Hence, we reduced to the following : 

    \begin{lemma}
        \label{LemmaLoc}
        For any $q \in \N$, the morphism : $$\rlim_{c \in \operatorname{Disj}_{X_I}} \bigoplus_{k=1}^{n(c)} R^q h_* (i_{T_k *} i_{T_k}^* \cL) \otimes_{\F_p} \O^+_{E_s}/p \to \rlim_{c \in \operatorname{Disj}_{X_I}} \bigoplus_{k=1}^{n(c)} R^q h_* (i_{T_k *} i_{T_k}^* (\cL \otimes_{\F_p} \O^+_Y/p)) $$ is an almost isomorphism.
    \end{lemma}

    \smallskip
    The key idea will be to view the $U \mapsto R^q h_* (i_{U *} i_U^* \cF)$ as presheaves on the profinite space $X_I$, and $\rlim_{c \in \operatorname{Disj}_{X_I}} \bigoplus_{i=1}^{n(c)} R^q h_* \left(i_{T_k *} i_{T_k}^* \cF\right)$ as a form of \textit{sheafification} of the above presheaves. Proving such an isomorphism of sheaves reduced to proving an isomorphism on stalks, which will compute the étale cohomology of the $Y_I$, for which the primitive comparison theorem \ref{PrimComparison_Galois} applies. 

    \smallskip
    \begin{defi}    
        \label{shfS}
        We let, for all $T \subset X_I$ and $q \ge 0$ : $$\mathcal{S}^q(\cF)(T) := R^q h_* (i_{T *} i_T^{*} \cF) \in \operatorname{Sh}(\operatorname{Spa}(E_s, \O_{E_s})_{\e t}, \operatorname{Ab}).$$
        
        We let $\cB$ denote the set of all clopen subsets of $X_I$. Since $X_I$ is profinite, this forms a basis of the topology on $X_I$ by \cite[\href{https://stacks.math.columbia.edu/tag/08ZZ}{Lemma 08ZZ}]{stacks-project}.
    \end{defi}
    
    Recall the notion of (pre)sheaves on a basis, as in \cite[\href{https://stacks.math.columbia.edu/tag/009H}{Section 009H}]{stacks-project}. In this formalism, the $S^q(\cF)$ are presheaves on $\cB$ valued in the abelian category $\operatorname{Sh}(\operatorname{Spa}(E_s, \O_{E_s})_{\e t}, \operatorname{Ab})$. 
    
    Let us start by two quick lemmas regarding sheaves on profinite spaces.


    \begin{lemma}
        \label{plusshf}
        Let $X$ be a profinite topological space, and $\cB$ be the basis formed by clopen subsets of $X$. Let $\cA$ be an abelian category, and $\cG$ be an $\cA$-valued presheaf on $\cB$.

        Then, the presheaf on $\cB$ defined, for any $U \in \cB$ by $\cG^\sharp(U) = \rlim_{c \in \operatorname{Disj}_U}  \bigoplus_{k=1}^{n(c)} \cG(T_k)$ defines a $\cB$-sheaf on $X$, valued in $\cA$ (in the sense of \cite[\href{https://stacks.math.columbia.edu/tag/009J}{Definition 009J}]{stacks-project}).
    \end{lemma}

    \begin{proof}
        We will apply the criterion of \cite[\href{https://stacks.math.columbia.edu/tag/009L}{Lemma 009L}]{stacks-project}. For any $U \in \cB$, the set of \textit{disjoint} clopen coverings is cofinal\footnote{Here, we use cofinality in the sense of the lemma of loc. cit. We do not claim cofinality in the sense of \cite[\href{https://stacks.math.columbia.edu/tag/04E6}{Definition 04E6}]{stacks-project}, as the second condition needs not hold.} in the coverings of $U$ by elements of $\cB$. 

        Hence, it suffices to check the sheaf condition relative to a disjoint clopen covering $U = \bigsqcup_{i \in I} U_i$, where $I$ is finite. Since the covering is disjoint, the cocycle condition is trivial, so that it suffices to check that $\cG^\sharp(U) \cong \bigoplus_{i \in I} \cG^\sharp(U_i)$.

        For all $c = \{ T_k \}_{1 \le k \le n(c)} \in \operatorname{Disj}_U$, the covering $\{ T_k \cap U_i\}_{(k, i) \in A \times I}$ remains in $\operatorname{Disj}_U$ and forms a refinement of the covering $\{U_i\}_{i \in I}$, so that the subcategory whose objects are : $$\operatorname{Ob} \left(\operatorname{Disj}_{U|U_i}\right) := \left\{ \{ T_k\}_{1 \le k \le n(c)} \in \operatorname{Disj}_U , \forall i \in I, \exists 1 \le k \le n(c), T_k \subset U_i \right\}$$ and whose morphisms are given by refinements, is cofinal in $\operatorname{Disj}_U$. Since the choices of refinements in $\operatorname{Disj}_U$ are canonical, the second condition of \cite[\href{https://stacks.math.columbia.edu/tag/04E6}{Definition 04E6}]{stacks-project} is automatically satisfied. Hence, we can compute the colimit defining $\cG^\sharp$ alongside $\operatorname{Disj}_{U|U_i}$, so that : $$\bigoplus_{i \in I} \cG^\sharp(U_i) = \bigoplus_{i \in I} \rlim_{c_i \in \operatorname{Disj}_{U_i}} \bigoplus_{k=1}^{n(c_i)} \cG(T_{k,i}) = \rlim_{c \in \operatorname{Disj_{U|U_i}}}\bigoplus_{k=1}^{n(c)} \cG(V_\alpha) \cong \rlim_{c \in \operatorname{Disj_{U}}}\bigoplus_{k=1}^{n(c)} \cG(V_\alpha) = \cG^\sharp(U).$$
        This concludes.
    \end{proof}

        
        

    \begin{rk}
        $\cG^\sharp$ can be thought of as the\textit{ $\cB$-sheafification of $\cG$}. The construction above is the analogue, with respect to the basis $\cB$, of the standard "plus" construction (as in \cite[\href{https://stacks.math.columbia.edu/tag/00W1}{Section 00W1}]{stacks-project}), that usually needs to be applied twice to sheafify. Here, once suffices, even for non separated presheaves.
    \end{rk}

    Recall that, for $\cG$ a $\cB$-presheaf on $X$, one can compute its stalk at some $x \in X$ by $\cF_{x} = \rlim_{U \in \cB, x \in U} \cF(U)$. Note that the plus construction of the above lemma is functorial. We have the following :

    \begin{lemma}
        \label{ext_shf}
        Let $X$ be a profinite topological space, with basis $\cB$ of clopen subsets as above. 
        
        Let $\cF$, $\cG$ be two $\cB$-presheaves on $X$, together with a morphism $f : \cF \to \cG$ of $\cB$-presheaves. Assume that $f$ induces isomorphisms on stalks $\cF_x \cong \cG_x$ for all $x \in X$. 
        
        Then, the induced morphism $f^\sharp : \cF^\sharp \to \cG^\sharp$ is an isomorphism of $\cB$-sheaves.
    \end{lemma}

    \begin{proof}
        By \cite[\href{https://stacks.math.columbia.edu/tag/009R}{Lemma 009R}]{stacks-project}, there is a functorial equivalence of category between $\cB$-sheaves on $X$ and sheaves on $X$, denoted $\cF \mapsto \cF^{ext}$, such that, for all $U \in \cB$, $\cF^{ext}(U) = \cF(U)$, and, for all $x \in X$, $\cF^{ext}_x = \cF_x$.

        Hence, the morphism $f^\sharp$ induces a morphism of sheaves $(\cF^\sharp)^{ext} \to (\cG^\sharp)^{ext}$, that induces an isomorphism of stalks, so that it is an isomorphism of sheaves. 
        
        Hence, so that, for all $U \in \cB$ the map $(\cF^\sharp)(U) \cong (\cF^\sharp)^{ext}(U) \to  (\cG^\sharp)^{ext}(U) \cong (\cG^\sharp)(U)$ is an isomorphism. This concludes.
    \end{proof}

    Let us now come back to the proof of \ref{PCTermSS}. Using the notations above, it is equivalent to the fact that : 
    $$ \mathcal{S}^q(\cL)^\sharp \otimes_{\F_p} \O^+_{E_s}/p \acong \mathcal{S}^q(\cL \otimes \O^+_Y/p)^\sharp.$$
    
    Hence, by \ref{ext_shf},  it suffices to show that the morphism $\mathcal{S}^q(\cL) \otimes_{\F_p} \O^+_{E_s}/p \to \mathcal{S}^q(\cL \otimes_{\F_p} \O^+_Y/p)$ of $\cB$-sheaves on $X_I$ induces an almost isomorphism on stalks. 

    \begin{lemma}
        \label{Stalk_of_S}
        For any $x \in X_I$, let $\operatorname{Neigh}_x$ be the category of compact open neighborhoods of $x \in X_I$. Let $\cF$ be a sheaf of abelian groups on $Y_{\e t}$. Then, there is a natural isomorphism $$ \mathcal{S}^q(\cF)_x \cong R^q (h \circ i_{x})_* \, i_{x}^* \, \cF$$

        where, for $x \in X_I$, $i_x$ denotes the closed immersion $x \cdot Y_I \subset X$.
    \end{lemma}

    \begin{proof}
        Fix $x \in X_I$. By \cite[Lemma 4.4]{Orlik_local}, for any family $(W_s)_{s \in \N}$ of compact open neighborhoods of the point $1 \cdot P_I(\Q_p) \in X_I$ such that $\bigcap_{s \in \N} W_s = \{ 1 \cdot P_I(\Q_p)\}$, we have $\bigcap_{s \in \N} Z_I^{W_s} = Y_I$. Fix such a family $W_s$. Hence, for all $x \in X$ the family $\{ x \cdot W_s\}_{s \in \N}$ forms a family of compact open neighborhoods of $x$ such that $\bigcap_{s \in \N} x \cdot W_s = \{ x\}$, and $\bigcap_{s \in \N} Z_I^{x \cdot W_s} = x \cdot Y_I$. 
        
        We let $W_s(x) = x \cdot W_s$. Then $\{ W_s(x)\}_{s \in \N}$ forms a cofinal family in $\operatorname{Neigh}_x$.

        Moreover, since $h$ is qcqs, by \cite[Prop 14.9]{Scholze_ECohD}, $R^q h_*$ commutes with filtered colimits, so that :$$\mathcal{S}^q(\cF)_x = \rlim_{W \in \operatorname{Neigh}_x} R^q (h \circ i_{W})_* i_{W}^* \cF \cong R^q h_* \rlim_{W \in \operatorname{Neigh}_x} i_{W*} i_{W}^* \, \cF$$
        
        Since $x \cdot Y_I \subset Z_I^{W(x)}$ for all compact open neighborhood $W(x)$ of $x$, there is a natural morphism (as in the discussion in \ref{Morph_Disj}) : 
        $$ \rlim_{W \in \operatorname{Neigh}_x}  i_{W*} i_{W}^* \cF \to i_{x*} i_{x}^* \cF$$ It suffices to prove that it is an isomorphism of étale sheaves on $Y$, and we will check this on fibers. Let $\eta : \operatorname{Spa}(C, C^+) \to Y$ be a geometric point, and $y \in Y$ be the image of the unique closed point. 
        
        Then, by \ref{stalkpf} we have, for any $W \in \operatorname{Neigh}_x$ :
$$ (i_{x *} \, i_{x}^* \cF)_{\overline{y} = }\begin{cases}
        \cF_{\overline{y}} &\text{ if } y \in x \cdot Y_I \\
        0 &\text{ otherwise}
     \end{cases} \text{ and } (i_{W*}\,i_W^* \cF)_{\overline{y}} =  \begin{cases}
        \cF_{\overline{y}} &\text{ if } x \in Z_I^W \\
        0 &\text{ otherwise}
     \end{cases}.$$

    Since stalks commute with colimits, it follows that $$\left(\rlim_{W \in \operatorname{Neigh}_x} i_{W*} i_{W}^* \cF \right)_{\overline{y}} \cong \begin{cases}
        \cF_{\overline{y}} &\text{ if }  y \in \bigcap_{W \in \operatorname{Neigh}_x } Z_I^W \\ 
        0 &\text{ otherwise}
    \end{cases}.$$

    This concludes since $\bigcap_{W \in \operatorname{Neigh}_x} Z_I^W = x \cdot Y_I$, as we explicited a cofinal system $W_s(x)$.
    \end{proof}

We are now ready to prove \ref{LemmaLoc}, hence \ref{PCTermSS}.

\begin{proof} \textit{of \ref{LemmaLoc}}. Using the reformulation above, we want to prove that the natural morphism : $$\mathcal{S}^q(\cL)^\sharp \otimes_{\F_p} \O^+_{E_s} \to \mathcal{S}^q(\cL \otimes_{\F_p} \O^+_Y/p)^\sharp$$
is an almost isomorphism. By the combined lemmas \ref{ext_shf} and \ref{Stalk_of_S} , this reduces to proving that, for all $x \in X_I$, the morphism : $$ R^q (h \circ i_x)_* \, i_x^* (\cL) \otimes_{\F_p} \O^+_{E_s}/p \to R^q (h \circ i_x)_* \, i_x^* (\cL \otimes \O^+_Y/p) $$ is an almost isomorphism. Pullbacks commutes with tensor products, and, since $i_x : x \cdot Y_I \to Y$ is a generalizing closed immersion, \ref{pro_et_pullback_is_restr} shows that $i^{*}_{x} \O_Y^+/p \cong \O_{x \cdot Y_I}^+/p$. Moreover, the pullback $i_x^* \cL$ defines an $\F_p$-local system on $x \cdot Y_I$, and $h \circ i_{x}$ is the structure morphism of $x \cdot Y_I$, which admits a structure of a proper rigid analytic variety, by \ref{Shubert}. 

The result now follows from the primitive comparison theorem \ref{PrimComparison_Galois}.
\end{proof}

\end{document}